\documentclass[10pt,a4paper,leqno]{article}
\usepackage{amsmath,amssymb,amscd,latexsym}

\usepackage[all]{xy}
\newcommand{\A}{{\mathbb{A}}}
\newcommand{\C}{{\mathbb{C}}}
\newcommand{\F}{{\mathbb{F}}}
\newcommand{\Ge}{\mathbb{G}}
\newcommand{\Lee}{\mathbb{L}}

\newcommand{\Q}{{\mathbb{Q}}}

\newcommand{\Pro}{{\mathbb{P}}}

\newcommand{\Z}{{\mathbb{Z}}}
\newcommand{\Zl}{{\mathbb{Z}}/\ell}
\newcommand{\Zln}{{\mathbb{Z}}/\ell^{\nu}}

\newcommand{\Ab}{\mathrm{Ab}}
\newcommand{\car}{\mathrm{char}\;}

\newcommand{\cocell}{\mathrm{cocell}}
\newcommand{\et}{\mathrm{\acute{e}t}}

\newcommand{\Ev}{\mathrm{Ev}}
\newcommand{\fib}{\mathrm{fib}}
\newcommand{\Ho}{\mathrm{Ho}}
\newcommand{\id}{\mathrm{id}}
\newcommand{\ok}{\overline{k}}
\newcommand{\Xok}{X_{\ok}}

\newcommand{\hOmega}{\hat{\Omega}}
\newcommand{\pro}{\mathrm{pro}}
\newcommand{\proj}{\mathrm{proj}}

\newcommand{\Spec}{\mathrm{Spec}\,}
\newcommand{\colim}{\mathrm{colim}}
\newcommand{\holim}{\mathrm{holim}}
\newcommand{\Et}{\mathrm{Et}\,}
\newcommand{\hEt}{\hat{\mathrm{Et}}\,}

\newcommand{\hEtsp}{\hat{\mathrm{Et}}_\mathrm{Sp}\,}
\newcommand{\thEtsp}{\tilde{\hat{\mathrm{Et}}}_\mathrm{Sp}\,}
\newcommand{\LhEt}{\mathrm{L}\hat{\mathrm{Et}}\,}
\newcommand{\LthEt}{\mathrm{L}\tilde{\hat{\mathrm{Et}}}\,}
\newcommand{\LhEtsp}{\mathrm{L}\hat{\mathrm{Et}}_\mathrm{Sp}\,}

\newcommand{\compl}{\hat{(\cdot)}}

\newcommand{\Gal}{\mathrm{Gal}}

\newcommand{\Hom}{\mathrm{Hom}}
\newcommand{\homp}{\mathrm{hom}_{\ast}}

\newcommand{\Imm}{\mathrm{Im}\,}
\newcommand{\Map}{\mathrm{Map}}
\newcommand{\Mapp}{\mathrm{Map}_{\ast}}

\newcommand{\Resp}{\mathrm{Res}^{\bullet}}
\newcommand{\sk}{\mathrm{sk}}
\newcommand{\Sch}{\mathrm{Sch}}

\newcommand{\Sm}{\mathrm{Sm}}
\newcommand{\Thom}{\mathrm{Th}}

\newcommand{\thomet}{\mathrm{th}_{\et}}
\newcommand{\Tor}{\mathrm{Tor}}
\newcommand{\Tot}{\mathrm{Tot}}
\newcommand{\Ch}{{\mathcal C}}
\newcommand{\Chk}{{\mathcal C}/\hEt k}
\newcommand{\Dh}{{\mathcal D}}
\newcommand{\Eh}{{\mathcal E}}
\newcommand{\hEh}{\hat{\mathcal E}}
\newcommand{\Fh}{{\mathcal F}}

\newcommand{\Hh}{{\mathcal H}}
\newcommand{\Hhp}{{\mathcal H}_{\ast}}
\newcommand{\hHh}{\hat{{\mathcal H}}}
\newcommand{\hHhp}{\hat{{\mathcal H}}_{\ast}}

\newcommand{\LCM}{L_{\Ch}\mathcal{M}}

\newcommand{\LeU}{L_{\et}U(k)}
\newcommand{\LNU}{L_{\mathrm{Nis}}U(k)}
\newcommand{\LU}{LU(k)}
\newcommand{\Lh}{{\mathcal L}}
\newcommand{\Mh}{{\mathcal M}}
\newcommand{\MVh}{\mathcal{MV}_k}

\newcommand{\Oh}{{\mathcal O}}

\newcommand{\Qh}{\mathcal{Q}}

\newcommand{\Sh}{{\mathcal S}}
\newcommand{\Shp}{{\mathcal S}_{\ast}}
\newcommand{\SHh}{{\mathcal{SH}}}
\newcommand{\SHhS}{{\mathcal{SH}^{S^1}(k)}}
\newcommand{\SHhP}{{\mathcal{SH}^{\Pro^1}(k)}}
\newcommand{\hSHh}{{\hat{\mathcal{SH}}}}

\newcommand{\hSHhk}{{\hat{\mathcal{SH}}}/\hEt k}
\newcommand{\hSHhtk}{{\hat{\mathcal{SH}}}_2/\hEt k}
\newcommand{\hSh}{\hat{\mathcal S}}

\newcommand{\hShp}{\hat{\mathcal S}_{\ast}}
\newcommand{\hShpk}{\hat{\mathcal S}_{\ast}/\hEt k}
\newcommand{\hSp}{\mathrm{Sp}(\hShp)}
\newcommand{\hSpk}{\mathrm{Sp}(\hShp/\hEt k)}
\newcommand{\Sp}{\mathrm{Sp}(\Sh)}
\newcommand{\SpC}{\mathrm{Sp}(\mathcal{C},T)}
\newcommand{\SpsC}{\mathrm{Sp}^{\Sigma}(\mathcal{C},K)}
\newcommand{\SpshSh}{\mathrm{Sp}^{\Sigma}(\hSh,S^1)}

\newcommand{\SymK}{\mathrm{Sym}(K)}
\newcommand{\SymKb}{\overline{\mathrm{Sym}(K)}}

\newcommand{\SpMVP}{\mathrm{Sp}(\MVh,\Pro^1 \wedge \cdot)}

\newcommand{\SpLUP}{\mathrm{Sp}(\LU,\Pro^1 \wedge \cdot)}

\newcommand{\Xh}{\mathcal{X}}

\newcommand{\hKU}{\hat{KU}}

\newcommand{\hMU}{\hat{MU}}

\newcommand{\tP}{\tilde{P}}
\newcommand{\tQ}{\tilde{Q}}

\newtheorem{theorem}{Theorem}[section]
\newtheorem{lemma}[theorem]{Lemma}
\newtheorem{prop}[theorem]{Proposition}
\newtheorem{defn}[theorem]{Definition}
\newtheorem{defnprop}[theorem]{Definition and Proposition}
\newtheorem{cor}[theorem]{Corollary}
\newtheorem{convention}[theorem]{Convention}
\newtheorem{example}[theorem]{Example}
\newtheorem{remark}[theorem]{Remark}
\newtheorem{assump}[theorem]{Assumption}
\newtheorem{conjecture}[theorem]{Conjecture}

\newenvironment{proof}{\noindent {\bf Proof}}{\mbox{}\hspace*{\fill}$\Box$}
\begin{document}
\title{Profinite Etale Cobordism}
\author{Gereon Quick}
\date{}
\maketitle
\begin{abstract}
In the present paper we construct a new cohomology theory for smooth schemes of finite type over a field, called profinite \'etale cobordism. The motivation for this theory is to develop an \'etale topological cobordism theory for the theories of algebraic cobordism of Voevodsky and Levine-Morel, that is easier to compute and admits interesting maps from algebraic cobordism to this new theory. Furthermore we hope to simplify the application of algebraic cobordism to questions in arithmetics. An important feature of \'etale cobordism is the existence of a convergent Atiyah-Hirzebruch spectral sequence starting from \'etale cohomology.\\
For the construction of the theory, we prove that there is a stable model structure on the category of simplicial profinite spectra. Furthermore we construct an \'etale realization on the stable motivic category of $\Pro^1$-spectra.\\
We prove that, over a separably closed field or a finite field, \'etale cobordism is an oriented cohomology theory which leads to canonical natural transformations from algebraic cobordism of Levine/Morel to \'etale cobordism. Over a separably closed field there are natural transformations from the algebraic cobordism of Voevodsky to the \'etale cobordism as well. Finally, the Atiyah-Hirzebruch spectral sequence gives rise to the conjecture that, over a separably closed field, algebraic and profinite \'etale cobordism with finite coefficients are isomorphic after inverting a Bott element.
\end{abstract}
\newpage

\setcounter{page}{1}

\tableofcontents

\newpage

\section{Introduction}

\subsection{Motivation and main results}

The aim of this paper is the construction of a new cohomology theory for smooth schemes over a field, called profinite \'etale cobordism.\\ 
Since the 1990s two approaches to the theory of algebraic cobordism for smooth schemes have been made.  On the one hand, Morel and Voevodsky have developed $\A^1$-homotopy theory for schemes. Voevodsky used it for the construction of cohomology theories on schemes. In particular, he showed the existence of an algebraic cobordism theory $MGL^{\ast,\ast}$ in analogy to Thom's homotopical definition of complex cobordism in topology.\\
On the other hand, Levine and Morel used Quillen's insight for a geometric construction of complex cobordism for proving that there is another possible definition of algebraic cobordism via a purely geometric construction. They proved that this $\Omega^{\ast}$ is the universal object for oriented cohomology theories on smooth schemes over a field. It is conjectured that $\Omega^{\ast}$ is in fact a geometrical description of the $MGL^{2\ast,\ast}$-part of Voevodsky's theory. Hopkins and Morel recently proved that the canonical map $\Omega^{\ast} \to MGL^{2\ast,\ast}$ is surjective over a field of characteristic zero.\\
Both approaches have been used to prove famous results. Voevodsky used his construction for the proof of the Milnor Conjecture. Levine and Morel proved Rost's Degree Formula. Nevertheless, both theories are still hard to compute and most of their features have only been proved over fields of characteristic zero.\\ 
The purpose of this paper is the construction of an \'etale topological version of cobordism for smooth schemes that is easier to compute since it is closely related to the \'etale cohomology. The idea is due to Eric Friedlander who constructed in \cite{etaleK} a first version of an \'etale topological K-theory for schemes that turned out to be a powerful tool for the study of algebraic K-theory with finite coefficients. In particular, Thomason proved in his famous paper \cite{thomason} that algebraic K-theory with finite coefficients agrees with \'etale K-theory after inverting a Bott element. The aim of this paper is also to construct a candidate for an analogous statement for algebraic cobordism, see Conjecture \ref{Bottconjecture2}.\\ 
At the end of the 1970s, in \cite{snaith} Victor P. Snaith has already constructed a $p$-adic cobordism theory for schemes. His approach is close to the definition algebraic K-theory by Quillen. He defines for every scheme $V$ over $\F_q$, $q=p^n$, a topological cobordism spectrum $\underline{A\F}_{q,V}$. The homotopy groups of this spectrum are the $p$-adic cobordism groups of $V$. He has calculated these groups for projective bundles, Severi-Brauer schemes and other examples. \\
In this paper, we do not follow Snaith's construction, but we provide Friedlander's idea with a general setting. We consider the profinite completion $\hEt$ of Friedlander's \'etale topological type functor of \cite{fried} with values in the category of simplicial profinite sets $\hSh$. We construct a stable homotopy category $\hSHh$ for $\hSh$ and define general cohomology theories for profinite spaces. We apply these cohomology theories to $\hEt X$ for a scheme $X$ over a field $k$. This yields a general foundation for different \'etale topological cohomology theories.\\
When we apply this construction to the cohomology theory represented by the profinitely completed $\hMU$-spectrum of complex cobordism, we get an \'etale topological cobordism theory $\hMU_{\et}^{\ast}$ for schemes of finite type over a field, which we call profinite \'etale cobordism. Note that this theory depends on a fixed prime number $\ell$, which must be different from the characteristic of the base field $k$.\\ 
The main feature of $\hMU_{\et}^{\ast}$ is the existence of an Atiyah-Hirzebruch spectral sequence, see Proposition \ref{coeffofhMUet} and Theorem \ref{MUahss}:
\begin{theorem}
1. Let $k$ be a separably closed field. There are isomorphisms
$$\hMU_{\et}^{\ast}(k)\cong MU^{\ast}\otimes_{\Z} \Z_{\ell}~\mathrm{and}~
\hMU_{\et}^{\ast}(k;\Zln)\cong MU^{\ast}\otimes_{\Z} \Zln.$$
2. Let $k$ be a field. For every smooth scheme $X$ over k, there is a convergent spectral sequence $\{E_r^{p,q}\}$ with 
$$E_2^{p,q}= H_{\et}^p(X;\Zln \otimes MU^q) \Longrightarrow \hMU_{\et}^{p+q}(X;\Zln).$$
\end{theorem}
These properties gives rise to several applications.\\
As a first corollary, we deduce from comparison results for \'etale cohomology that over $k=\C$ profinite \'etale and complex cobordism with finite coefficients are naturally isomorphic, see Theorem \ref{complexvar}.\\
On the one hand, based on the results of Panin \cite{panin}, they enable us to prove, see Theorem \ref{transfer}:
\begin{theorem}
Let $k$ be a separably closed field. Profinite \'etale cobordism $\hMU_{\et}^{\ast}(-)$, resp. $\hMU_{\et}^{\ast}(-;\Zln)$, is an oriented cohomology theory on $\Sm/k$
\end{theorem}
This implies that we have a unique natural morphism $\theta_{\hMU}:\Omega^{\ast}(X) \to \hMU_{\et}^{2\ast}(X)$ for every $X$ in $\Sm/k$. 
On the other hand, via a stable \'etale realization functor of motivic spectra, we construct a natural map $\phi:MGL^{\ast,\ast}(X) \to \hMU_{\et}^{\ast}(X)$, see Section 8.2, Theorem \ref{triangle} and Corollary \ref{MGLsurjection}: 
\begin{theorem}
This natural map $\phi$ fits into a commutative triangle for every $X$ in $\Sm/k$
$$\begin{array}{ccc}
\Omega^{\ast}(X) & \stackrel{\theta_{MGL}}{\longrightarrow} &  MGL^{2\ast,\ast}(X)\\
\theta_{\hMU} \searrow &  & \swarrow \phi \\
 & \hMU_{\et}^{2\ast}(X). & 
\end{array}$$
There is a similar triangle for $\Zln$-coefficients. \\
\end{theorem}
For a field $k$, we start the study of the absolute Galois group $G_k=\Gal(k_s/k)$ on \'etale cobordism, where $k_s$ denotes a separable closure of $k$. We deduce from the previous results that the action of $G_k$ on $\hMU^{\ast}_{\et}(k_s)$, resp. $\hMU^{\ast}_{\et}(k_s;\Zln)$, is trivial, see Theorem \ref{Galoisaction}. Together with the Atiyah-Hirzebruch spectral sequence applied to Galois cohomology, this enables us to determine the \'etale cobordism of a finite field $k=\F_q$, $\ell \not | ~q$, see Theorem \ref{hMUofFq}.\\
Finally, the Atiyah-Hirzebruch spectral sequence together with the results of Levine \cite{bott} on motivic cohomology with inverted Bott element yield good reasons for the following conjecture, see Section 8.4 and Conjecture \ref{Bottconjecture}:
\begin{conjecture}\label{Bottconjecture2}
Let $X$ be a smooth scheme of finite type over a separably closed field $k$ of characteristic different from $\ell$. Suppose that $\ell$ is odd or that $\ell^{\nu}\geq 4$. Let $\beta \in MGL^{0,1}(k;\Zln)$ be the Bott element. The induced morphism
$$\phi: MGL^{\ast,\ast}(X;\Zln)[\beta ^{-1}] \to \hMU_{\et}^{\ast}(X;\Zln)$$
is an isomorphism.
\end{conjecture}
At the end of Section 8, we will explain a strategy to prove this conjecture, which is close to Levine's new proof of Thomason's K-theory theorem in \cite{levineKtheory}. In particular, we would like to use the Atiyah-Hirzebruch spectral sequence of Hopkins and Morel for algebraic cobordism in \cite{homo}.
\subsection{Survey of the paper}
The outline of the paper is as follows. Towards the construction of profinite \'etale cobordism, there are several technical problems to solve. The first main result of this paper is that it is possible to establish a stable model structure on profinite spectra. It is based on Morel's $\Zl$-cohomological model structure on $\hSh$ of \cite{ensprofin}, see Theorem \ref{thmprostable}:
\begin{theorem}
There is a stable model structure on profinite spectra $\hSp$ such that the suspension functor $S^1\wedge-$ is an equivalence on the corresponding homotopy category $\hSHh$.
\end{theorem}
This result cannot be deduced in the same way as the stable structure on simplicial spectra, since $\hSh$ is not proper. We use a generalized theorem on left Bousfield localization of model categories similar to a result of Hirschhorn \cite{hirsch} and a generalized result on stable model structure on spectra similar to the work of Hovey in \cite{hovey}. We have discussed these theorems in the appendix.\\
The next step is the study of generalized cohomology theories on profinite spaces. Our knowledge on profinite completion of spectra and the results of Bousfield-Kan on $\ell$-completion of spaces enables us to determine the coefficients of $\hMU$. We also construct an Atiyah-Hirzebruch spectral sequence for profinite spaces. The proof of the existence of this spectral sequence is mainly due to Dehon in \cite{dehon}. Profinitely completed cobordism has already been used for other purposes. In particular, Francois-Xavier Dehon uses the completed cobordism theory $\hMU^{\ast}(-)$ for the study of Lannes' $T$-functor in \cite{dehon}. Many of his ideas were fruitful for this paper.\\  
With the stable category on $\hSh$ we give this approach a general setting. Hence the discussions of the first sections may already be interesting on their own. \\
In order to clarify the argument and the relevance of the discussions of Sections 2-4 for the reader, we mention the following point.\\  
Since our definitions are all compatible with \'etale cohomology in a certain sense, see Remark \ref{remarketalecohom}, we can deduce from these facts about profinite cohomology theories the existence of the mentioned Atiyah-Hirzebruch spectral sequence starting from \'etale cohomology and converging to profinite \'etale cobordism. Furthermore, we can determine the coefficients $\hMU_{\et}^{\ast}(k)$ for a separably closed field $k$.\\
In Sections 5 and 6, we dedicate our attention to the study of the completed \'etale topological type functor $\hEt$ on schemes of finite type over a field with values in $\hSh$. This functor is due to Artin-Mazur and Friedlander. It has been extended to the $\A^1$-homotopy category of schemes  independently by Daniel Isaksen \cite{a1real} and Alexander Schmidt \cite{schmidt}.\\ 
We use this functor to construct a stable \'etale realization of motivic spectra, see Section 6 and Theorem \ref{P1real}: 
\begin{theorem}
The functor $\hEt$ yields an \'etale realization of the stable motivic homotopy category of $\Pro^1$-spectra: $$\LhEt: \SHhP \to \hSHh.$$
\end{theorem}
This is another main result of the paper and may be interesting on its own. Although we do not know if $\hEt MGL$ is isomorphic to $\hMU$ in $\hSHh$ over an arbitrary base field, this realization is the key ingredient in the construction of the map $\phi:MGL^{\ast,\ast}(X) \to \hMU_{\et}^{\ast}(X)$ from algebraic to profinite \'etale cobordism.\\
In Sections 7 and 8 we start the study of \'etale cohomology theories and the comparison with algebraic cobordism. These two sections contain the study of profinite \'etale cobordism and give the proofs for the main results that we mentioned above. \\
The reader could even start with Section 7 if one accepts first the existence of a stable category of profinite spectra, second the existence of the functor $\hEt$ on schemes and motivic spectra and third the existence of an Atiyah-Hirzebruch spectral sequence for profinite generalized cohomology theories and the isomorphism $\hMU^{\ast} \cong MU^{\ast} \otimes_{\Z} \Z_{\ell}$. For the last fact, we even give another proof in Section 7.\\
The purpose of the appendix is to prove a general stabilization theorem for left proper fibrantly generated model categories, see Theorem \ref{stablestructure}: 
\begin{theorem}
Let $\Ch$ be a left proper fibrantly generated simplicial model category with all small limits and let $T$ be a left adjoint endofunctor on $\Ch$. There is a model structure on the category of spectra $\SpC$ such that $T$ becomes a Quillen equivalence on $\SpC$.\\ 
This model structure is again a left proper fibrantly generated simplicial model structure. 
\end{theorem}
We apply this theorem in Section 3 to the left proper fibrantly generated model category of profinite spaces. This stabilization is the only purpose of the appendix. But since we gave the proofs in a very general setting, we have postponed them to the end of the paper.\\
In Appendix 1, we prove first a theorem on left Bousfield localization of fibrantly generated model categories, based on \cite{hirsch}.\\
In Appendix 2, this theorem is applied to an intermediate model structure on spectra in order to get a stable category, based on \cite{hovey}. The ideas for the proofs are due to the papers \cite{hovey} of Hovey and \cite{hirsch} of Hirschhorn plus an idea of Christensen-Isaksen in \cite{prospectra}. Nevertheless, our results are slightly more general and might be of interest for other situations. 

I would like to use this opportunity to thank first of all Prof. Dr. Christopher Deninger for the suggestion of the topic of this thesis which satisfied all my interests and gave me the opportunity to learn a lot. He helped me with important suggestions and gave me the motivating impression to work on something that captured his interest. I would like to thank Dr. Christian Serp\'e for many very helpful discussions and suggestions for further studies. I am grateful to Prof. Marc Levine, Prof. Fabien Morel and Prof. Dr. Alexander Schmidt, who answered all my questions and explained to me both the great picture and the little details. I thank Francois-Xavier Dehon for an answer to a specific problem on K\"unneth isomorphisms for cobordism with coefficients.

\section{Profinite spaces}

We introduce the category of (pointed) profinite spaces $\hSh$ (resp. $\hShp$) and its model structure defined by Morel \cite{ensprofin} for a fixed prime number $\ell$. An important construction is the Bousfield-Kan-$\ell$-completion of profinite spaces, which is an explicit fibrant replacement functor in this model structure.\\
Via this fibrant replacement functor, we introduce homotopy groups of profinite spaces and study their behavior under the profinite completion functor $\Sh \to \hSh$. In particular, we show that these homotopy groups are always pro-$\ell$-groups.\\
Then we define a completion functor $\pro-\Sh \to \hSh$ from the category of pro-simplicial sets to profinite spaces and show that the respective $\Zl$-homotopy categories are equivalent. This compares our approach to the one of \cite{compofpro}.\\
We finish this section with a technical statement. The model structure of $\hSh$ is fibrantly generated. This is necessary in order to apply the localization and stabilization theorems of the appendix to $\hSh$. At this point, the category $\hSh$ is less well behaved than $\Sh$ and $\pro-\Sh$, but it has properties that are good enough for our purposes.  
    
\subsection{The profinite $\ell$-completion of Bousfield-Kan}

Let $\Eh$ denote the category of sets and let $\Fh$ be the full subcategory of $\Eh$
whose objects are finite sets. Let $\hEh$ be the category of compact and totally disconnected
topological spaces. We may identify $\Fh$ with a full subcategory of $\hEh$ in the obvious way. The limit functor $\lim :\pro-\Fh \to \hEh$, which sends a pro-object $X$ of $\Fh$ to the limit in $\hEh$ of the diagram corresponding to $X$, is an equivalence of categories.\\
We denote by $\Sh$ (resp. $s\Fh$, resp. $\hSh$) the category of simplicial sets (resp. simplicial finite sets, resp. simplicial profinite sets). The objects of $\Sh$ (resp. $\hSh$) will be called {\em spaces} (resp. {\em profinite spaces}). There is a subtle distinction between the two categories pro-$s\Fh$ and $\hSh$. The obvious functor between them is not an equivalence. See \cite{calclim} for more details and a counter example.\\
For a profinite space $X$ we define the ordered set $\Qh(X)$ of simplicial open equivalence relations on $X$. For every element $Q$ of $\Qh(X)$ the quotient $X/Q$ is a simplicial finite set and the map $X \to X/Q$ is a map of profinite spaces. In fact, when we consider the limit $\lim_{Q\in \Qh(X)} X/Q$ in $\hSh$, the map $X \to \lim_{Q\in \Qh(X)} X/Q$ is an isomorphism, cf. \cite{quillen2} , Lemma 2.3.\\
The forgetful functor $\hEh \to \Eh$ admits a left adjoint $\compl:\Eh \to \hEh$. It induces dimensionwisely a functor $\compl:\Sh \to \hSh$, which we call profinite completion. It is left adjoint to the forgetful functor $|\cdot|:\hSh \to \Sh$ which sends a profinite space $X$ to its underlying simplicial set $|X|$. By adjunction, the profinite completion of a simplicial set $Z$ may be identified with the filtered colimit in $\hSh$ of the simplicial finite subsets of $Z$. \\
Let $X$ be a profinite space. The continuous cohomology $H^{\ast}(X;\pi)$ of $X$ with coefficients in the profinite abelian group $\pi$ is defined as the cohomology of the complex $C^{\ast}(X;\pi)$ of continuous cochains of $X$ with values in $\pi$, i.e. $C^n(X;\pi)$ denotes the set $\Hom_{\hat{\Eh}}(X_n,\pi)$ of continuous maps $X_n \to \pi$.\\ 
\begin{remark}\label{remarkcohom}{\rm (\cite{ensprofin}, \S 1.2.)}\\
1. Let $X$ be a simplicial profinite set. If we denote by $C_{\ast}(X)$ the chain complex defined in degree $n$ to be the free profinite abelian group on the profinite set $X_n$, then the universal property of the free profinite group yields for every profinite abelian group $\pi$ an isomorphism $C^n(X)=\Hom_{\hEh}(X_n, \pi) \cong \Hom(C_n(X), \pi)$ for each $n$ where the last $\Hom$  denotes the morphisms of profinite abelian groups. This shows the relation to the classical definition of the cohomology of simplicial sets, cf. \cite{may} p. 5.\\
2. Let $\pi$ be a finite abelian group. There is a natural isomorphism between the continuous cohomology with coefficients in $\pi$ of the profinite completion $\hat{Y}$  of a simplicial set $Y$ and its ordinary cohomology with coefficients in $\pi$, i.e. 
$$H^{\ast}(Y;\pi) \cong H^{\ast}(\hat{Y};\pi).$$ This isomorphism exists already on the level of complexes, since there is a natural bijection $\Hom_{\hEh}(\hat{X}_n,\pi)\cong \Hom_{\Eh}(X_n,\pi)$. We will explain later how to generalize this result to an arbitrary abelian pro-$\ell$-group $\pi$.\\
3. Let $\pi$ be a finite abelian group and let $\{X_s\}$ be a cofiltered diagram of profinite spaces. The canonical map 
$$\colim_s H^{\ast}(X_s;\pi) \stackrel{\cong}{\longrightarrow} H^{\ast}(\lim_s X_s;\pi)$$ 
is an isomorphism, cf. \cite{dehon}, Lemme 1.1.1.\\ 
4. Let $\pi$ be a finite abelian group, the complex $C^{\ast}(X;\pi)$ is naturally identified with the filtered colimit of complexes  $\colim_Q C^{\ast}(X/Q;\pi)$, with $Q $ running through $\Qh(X)$, and the cohomology  $H^{\ast}(X;\pi)$ is naturally isomorphic to the colimit of the cohomologies $H^{\ast}(X/Q;\pi)$.\\
5. There is a K\"unneth formula, i.e. for two profinite spaces $X$ and $Y$ the map
$H^{\ast}(X;\Zl) \otimes H^{\ast}(Y;\Zl) \to H^{\ast}(X\times Y;\Zl)$ is an isomorphism, cf. \cite{ensprofin}.
\end{remark} 
Let $\ell$ be a fixed prime number. Fabien Morel has shown in \cite{ensprofin} that the category $\hSh$ can be given the structure of a closed model category in the sense of \cite{homalg}. The weak equivalences are the maps inducing isomorphisms in continuous cohomology with coefficients $\Zl$; the cofibrations are the degreewise monomorphisms and the fibrations are the maps that have the right lifting property with respect to the cofibrations that are also weak equivalences. One should note that there are functorial factorizations for any map in $\hSh$ into a trivial cofibration followed by a fibration, respectively into a cofibration followed by a trivial fibration. For example, we get a functorial factorization by fixing the construction of the proof of Proposition 1 on page 355 of \cite{ensprofin} on the existence of factorizations. All the constructions done there are functorial. Furthermore, all small limits and products exist in $\hSh$.\\
Morel has also given an explicit construction of fibrant replacements in $\hSh$, cf. \cite{ensprofin}, 2.1.
It is based on the $\Zl$-completion functor of \cite{bouskan}. Let $Y$ be a simplicial set. We denote by $\Resp Y$ its cosimplicial $\Zl$-resolution, by $\Tot_s (\Resp Y)$ its $s$-th total space and by 
$P^t \Tot_s \Resp Y$ its $t$-th Postnikov decomposition, cf. \cite{bouskan}, part I. If $Y$ is a finite set in each degree, then  $\Tot_s \Resp Y$ is also finite in each degree for all $s\geq 0$. Hence the total space $\Tot \Resp Y$, i.e. the $\ell$-completion of Bousfield-Kan of $Y$, which is the limit of the tower
$\Tot_s \Resp Y$, has a natural structure of a profinite space, which we denote by $\hat{Y}^{\ell}$.  
Now let $X$ be a profinite space. We denote by $\hat{X}^{\ell}$ the limit in $\hSh$ of the 
$\widehat{X/Q}^{\ell}$, $Q$ running through $\Qh(X)$. We call this space the {\em $\ell$-completion of Bousfield-Kan} of the profinite space $X$. For a simplicial set $Y$ we denote by $\hat{Y}^{\ell} \in \hSh$ its {\em profinite} $\ell$-completion defined as the limit $\lim_{Q,s,t} P^t\Tot_s \Resp (X/Q)$ in $\hSh$, where $X/Q$ are the simplicial finite quotients of $X$. Its underlying simplicial set is isomorphic in $\Ho(\Sh)$ to the $\ell$-completion in $\Sh$ of Sullivan \cite{sullivan}, which is the limit of the pro-Artin-Mazur-$\ell$-completion, cf. \cite{cohomodp}.\\   
Morel shows that the natural map $\theta_X: X \to \hat{X}^{\ell}$ defines a functorial fibrant replacement in $\hSh$, \cite{ensprofin}, 2.1, Prop. 2. The spaces $P^t \Tot_s (\Resp X/Q)$ are in fact $\ell$-finite spaces, i.e. fibrant simplicial finite sets whose homotopy groups for every choice of basepoint are finite $\ell$-groups, trivial except for a finite number of them. The profinite space $\hat{X}^{\ell}$ may be identified with the filtered limit  in $Q$, $s$ and $t$ of the $\ell$-spaces $P^t \Tot_s (\Resp X/Q)$. Hence $\hat{X}^{\ell}$ is a pro-$\ell$-space in the terminology of \cite{ensprofin}. \\
For the study of generalized cohomology theories on profinite spaces it will be crucial to check the compatibility of the $\ell$-completion with other constructions in $\hSh$, resp. $\hShp$, such as products, smash products and function spaces.

\subsection{Pointed profinite spaces and the simplicial structure}

The category $\hSh$ has a pointed analogue $\hShp$ whose objects are maps $\ast \to X$ in $\hSh$, where $\ast$ denotes the constant simplicial set equal to a point. Its morphisms are maps in $\hSh$ that respect the basepoints. The forgetful functor $\hShp \to \hSh$ has a left adjoint, which consists in adding a disjoint basepoint $X \mapsto X_+$. The $\ell$-completion of a pointed profinite space is naturally pointed and $X \mapsto X_+$ is compatible with $\ell$-completion. $\hShp$ has the obvious induced closed model category structure.\\
The product for two pointed profinite spaces $X$ and $Y$ in $\hShp$ is the smash product 
$X \wedge Y \in \hShp$, defined in the usual way as the quotient $(X \times Y)/(X \vee Y)$ in $\hSh$. For a pointed profinite space $X$ we define its reduced cohomology with coefficients in the profinite abelian group $\pi$, denoted $\tilde{H}^{\ast}(X;\pi)$, to be the kernel of the induced morphism
$H^{\ast}(X;\pi) \to H^{\ast}(pt ;\pi)$. It is clear that  $\tilde{H}^{\ast}(X_+;\pi) \cong H^{\ast}(X;\pi)$. We also have a K\"unneth formula.\\
\begin{lemma}\label{smashKunneth}
For two pointed profinite spaces $X$ and $Y$, the canonical map
$\tilde{H}^{\ast}(X;\Zl) \otimes \tilde{H}^{\ast}(Y;\Zl) \to \tilde{H}^{\ast}(X \wedge Y;\Zl)$ is an isomorphism.
\end{lemma} 
\begin{proof}
The canonical map is an isomorphism in the case that $X$ and $Y$ are simplicial finite sets. Hence it is an isomorphism in general since the tensor product commutes with filtered colimits.
\end{proof}

\begin{prop}\label{productl}
Let $X$ and $Y$ be simplicial sets whose $\Zl$-cohomology is finite in each degree. Then the map 
$\widehat{X\times Y}^{\ell} \to \hat{X}^{\ell} \times \hat{Y}^{\ell}$ is a weak equivalence in $\hSh$.
\end{prop}
\begin{proof}
We have $\widehat{X\times Y}^{\ell} = \lim_Q \widehat{(X\times Y)/Q}^{\ell} =
\lim_{Q_1,Q_2} \widehat{(X/Q_1 \times Y/Q_2)}^{\ell}$, which is due to the fact that the finite quotients 
$(X\times Y)/Q$ are in 1-1-correspondence with finite quotients $X/Q_1\times Y/Q_2$;
and we have $ \hat{X}^{\ell} \times \hat{Y}^{\ell} = 
\lim_{Q_1} \widehat{X/Q_1}^{\ell} \times \lim_{Q_1} \widehat{Y/Q_2}^{\ell}$. The hypothesis implies that $X$ and $Y$ are nilpotent spaces in the sense of \cite{bouskan}. Thus the same is true for the finite quotients $X/Q_1$ and $Y/Q_2$. Then \cite{bouskan} VI, \S 6 Proposition 6.5, says that 
$$H^{\ast}(\widehat{(X/Q_1 \times Y/Q_2)}^{\ell};\Zl) \to H^{\ast}(\widehat{X/Q_1}^{\ell} \times \widehat{Y/Q_2}^{\ell};\Zl)$$ is an isomorphism. We conclude by taking colimits and by Remark \ref{remarkcohom}.
\end{proof}

\begin{prop}\label{smashl}
Let $X$ and $Y$ be pointed simplicial sets whose $\Zl$-cohomology is finite in each degree. Then the map $\widehat{X\wedge Y}^{\ell} \to \hat{X}^{\ell} \wedge \hat{Y}^{\ell}$ is a weak equivalence in $\hShp$.
\end{prop}
\begin{proof}
This follows as the previous proposition from \cite{bouskan} VI, \S 6, Proposition 6.6.
\end{proof}

The categories $\hSh$ and $\hShp$ have natural simplicial structures. There is a simplicial set 
$\Map_{\hShp}(X,Y)=\Mapp(X,Y)$ for all profinite spaces $X$ and $Y$. It is naturally pointed by the map $X \to \ast \to Y$. It is characterized by the bijection
$$\Hom_{\Sh_{\ast}}(Z, \Mapp(X,Y)) \cong \Hom_{\hShp}(X\wedge Z ,Y)$$
for every pointed simplicial finite set $Z$. The space $\Mapp(X,Y)$ is given in the $n$-th dimension by the set $\Hom_{\hShp}(X \wedge \Delta[n]_+ ,Y)$, where $\Delta[n]_+$ is the well-known pointed simplicial finite set. For a pointed simplicial finite set $Z$ the space $\Mapp(Z,Y)$ has a natural structure of a profinite space considered as the filtered limit in $\hShp$ of the simplicial finite sets $\Mapp(\sk_nZ,Y/Q)$. We denote this profinite space by $\homp(Z,Y) \in \hShp$. For an arbitrary simplicial set $Z$ we define the profinite space $\homp(Z,Y)$ by taking the filtered limit in $\hShp$ of the profinite spaces $\homp(Z_s,Y)$ over all simplicial finite subsets $Z_s$ of $Z$. This defines a right adjoint to the functor $X \mapsto Z \wedge X$, $\hShp \to \hShp$. We define the functor $X\otimes \cdot := X\wedge \cdot : \Sh_{\ast} \to \hShp$ for a profinite space $X$ by taking the filtered limit in $\hShp$ of the profinite spaces $X/Q \wedge Z_s$ for all simplicial finite subsets $Z_s$ of the simplicial set $Z$ and all finite quotients $X/Q$ of $X$. These functors define a simplicial structure on $\hShp$ and similarly on $\hSh$.\\
Let $Z \leftarrow X \rightarrow Y$ be a diagram of profinite spaces and let $X\to Y$ be a cofibration. We denote its colimit by $Z \cup_X Y$. On the level of continuous cochains we get that $C^{\ast}(X,\Zl)$ is the colimit of the induced diagram $C^{\ast}(Z,\Zl) \leftarrow C^{\ast}(Z \cup_X Y,\Zl) \to C^{\ast}(Y,\Zl)$. This yields a Mayer-Vietoris long exact sequence in cohomology
$$\ldots \to H^n(Z \cup_X Y;\Zl) \to H^nY \oplus H^nZ \to H^nX \to H^{n+1}(Z \cup_X Y;\Zl) \to \ldots$$
Using this sequence Dehon shows that $\hSh$ and $\hShp$ carry in fact a {\em simplicial} model structure, cf. \cite{dehon} \S 1.4. \\
As an application, we define for $X \in \hShp$ the profinite space $\hat{\Omega}X \in \hShp$. Let $S^1$ be the simplicial finite set $\Delta[1]/\partial \Delta[1]$. We set $\hat{\Omega}X:=\homp(S^1,X)$. There is a natural isomorphism $\Hom_{\hShp}(S^1 \wedge X, Y) \cong \Hom_{\hShp}(X,\hat{\Omega}Y))$. In addition, $\hat{\Omega}X$ is fibrant if $X$ is fibrant.\\
This adjunction may be extended to the homotopy category $\Ho(\hShp)$ of $\hShp$. For every profinite space $X$ and every fibrant profinite space $Y$, we have a natural bijection
$\Hom_{\Ho(\hShp)}(S^1 \wedge X, Y) \cong \Hom_{\Ho(\hShp)}(X,\hat{\Omega}Y))$. Since $S^1\wedge X$ is a cogroup object and $\hOmega X$ a group object in $\Ho(\hShp)$, we conclude that the previous bijection is in fact an isomorphism of groups.\\
The following proposition has been proved by Francois-Xavier Dehon in \cite{dehon} using its Eilenberg-Moore spectral sequence for profinite spaces. We give an alternative proof using the results of Bousfield and Kan \cite {bouskan}.\\

\begin{prop}\label{omegal}
Let $X$ be a simply connected simplicial set whose $\Zl$-cohomology is finite in every degree. Then the map $\widehat{\Omega X}^{\ell} \to \hat{\Omega} (\hat{X}^{\ell})$ is a weak equivalence in $\hShp$. 
This map is in fact a homotopy equivalence.
\end{prop}     
\begin{proof}
Since $S^1$ is a simplicial finite set, the finite quotients $\hat{X}/Q$ of $\hat{X}$ and $\homp(S^1,X)/Q$ of $\homp(S^1,X)$ are in 1-1 correspondence. Hence we may write 
$$\hOmega (\hat{X}^{\ell}) = \homp(S^1,\lim_Q \widehat{X/Q}^{\ell}) = 
\lim_Q \homp(S^1, \widehat{X/Q}^{\ell}) = \lim_Q \hOmega (\widehat{X/Q}^{\ell});$$ 
and $$\widehat{\Omega X}^{\ell} = \lim_Q \widehat{(\homp(S^1,X)/Q)}^{\ell} =
\lim_Q \widehat{(\hOmega(X/Q))}.$$ 
By our hypothesis $X$ is a nilpotent space in the sense of \cite{bouskan} and hence $X/Q$ is nilpotent, too; by \cite{bouskan} VI, \S 6, Prop. 6.5, we know that
$$H^{\ast}(\widehat{(\hOmega(X/Q))}^{\ell}; \Zl) \cong H^{\ast}(\hOmega (\widehat{X/Q}^{\ell});\Zl).$$ 
Hence we get 
$$\begin{array}{rcl}
H^{\ast}(\widehat{(\hOmega(X))}^{\ell}; \Zl)&  \cong & \colim_Q H^{\ast}(\widehat{(\hOmega(X/Q))}^{\ell});\Zl)\\
  & \cong & \colim_Q H^{\ast}(\hOmega (\widehat{X/Q}^{\ell});\Zl) \cong H^{\ast}(\hOmega (\widehat{X}^{\ell});\Zl),
\end{array}$$ 
which is what we had to prove.\\  
The last statement follows from the fact that the spaces in question are cofibrant and fibrant and the fact that weak equivalences between such objects are homotopy equivalences. 
\end{proof}

\subsection{Homotopy groups of profinite spaces}

\begin{defn}
Let $X \in \hSh$ be a pointed profinite space and let $\hat{X}^{\ell}$ be its $\ell$-completion. 
We denote by $\pi_0X$ the coequalizer in $\hEh$ of the diagram $d_0,d_1:X_1 \stackrel{\to}{\to} X_0$. 
We define $\pi_nX$ for $n\geq 1$ to be the group $\pi_n|\hat{X}^{\ell}|$ where $|\hat{X}^{\ell}|$ is the underlying fibrant simplicial set of $\hat{X}^{\ell}$. We call $\pi_nX$ the {\em $n$-th homotopy group of $X$}.   
\end{defn}

\begin{prop}\label{piequiv}
A map $f:X\to Y$ in $\hShp$ is a weak equivalence if and only if $\pi_{\ast}(f)$ is an isomorphism.
\end{prop}
\begin{proof}
The maps $X\to \hat{X}^{\ell}$ and $Y\to \hat{Y}^{\ell}$ are weak equivalences. Hence $f$ is a weak equivalence if and only if $\hat{f}^{\ell}$ is a weak equivalence. Since all objects in $\hShp$ are cofibrant, the map $\hat{f}^{\ell}$ is a weak equivalence if and only if $\hat{f}^{\ell}$ is a homotopy equivalence. This finishes the proof. 
\end{proof}

We can already conclude from the definition 
$\pi_nX = \Hom_{\Ho(\Sh_{\ast})}(S^n,|\widehat{X}^{\ell}|)$ that $\pi_nX$ has a natural structure of a profinite group. The following proposition shows that our definition coincides with the one given by Dehon in \cite{dehon}.

\begin{prop}\label{piomegahat}
1. For every profinite space $X$, the natural map $$\Hom_{\Ho(\hSh)}(\ast , X) \to \pi_0 X$$ is an isomorphism. \\
2. For every pointed profinite space $X$ and every $n\geq 0$ we have
\begin{equation}\label{eqpiomegahat}
\pi_n(\hat{\Omega}X)\cong \pi_{n+1}(X).
\end{equation}
In particular, we have $\pi_n X \cong \pi_0(\hOmega^nX)$.
\end{prop}
\begin{proof}
The first statement is Proposition 1.3.4 a) of \cite{dehon}.\\
The second statement follows from the fact that $S^n$ is a simplicial finite set for all $n\geq 1$. For $n=0$, we use the first statement to conclude by adjunction
$$\begin{array}{rcl}
\pi_0(\hOmega X) & \cong & \Hom_{\Ho(\hShp)}(\ast , \hOmega X) \cong 
\Hom_{\Ho(\hShp)}(\ast , \hOmega \hat{X}^{\ell}) \\ 
  & \cong & \Hom_{\Ho(\hShp)}(S^1 , \hat{X}^{\ell}) \cong \Hom_{\Ho(\Sh_{\ast})}(S^1, |\hat{X}^{\ell}|) \\
 & =  & \pi_1 |\hat{X}^{\ell}| = \pi_1 X.
\end{array}$$
For $n \geq 1$, we have the following sequence of isomorphisms using Proposition \ref{omegal} and adjunction
$$\begin{array}{rcl}
\pi_n(\hat{\Omega}X) & = & \pi_n(|\widehat{\hat{\Omega}X}^{\ell}|) =  \Hom_{\Ho(\Sh_{\ast})}(S^n,|\widehat{\hat{\Omega}X}^{\ell}|)
\cong \Hom_{\Ho(\hShp)}(\widehat{S^n},\widehat{\hat{\Omega}X}^{\ell})\\
 & =  & \Hom_{\Ho(\hShp)}(S^n,\widehat{\hat{\Omega}X}^{\ell}) \cong \Hom_{\Ho(\hShp)}(S^n,\hat{\Omega}(\hat{X}^{\ell}))\\
 & \cong & \Hom_{\Ho(\hShp)}(S^{n+1},\hat{X}^{\ell}) \cong \Hom_{\Ho(\hShp)}(S^{n+1},\widehat{X}^{\ell})\\
  &  \cong & \Hom_{\Ho(\Sh_{\ast})}(S^{n+1},|\widehat{X}^{\ell}|)\\ 
 & = & \pi_{n+1}(X).
\end{array}$$
\end{proof}

As a corollary, we see that $\pi_n X$ has the natural structure of a pro-$\ell$-group for every $n\geq 1$, see also \cite{cohomodp}, Cor. 1.4.1.3. 

\begin{prop}\label{pro-l-group}
For every pointed profinite space $X$ and every $n\geq 0$, we have an isomorphism
$\pi_n X \cong \lim_{Q,s,t} \pi_n (P^t \Tot_s(\Resp X/Q))$.\\
In particular, for every $n \geq 1$ the group $\pi_n X$ has the structure of a pro-$\ell$-group.  
\end{prop}
\begin{proof}
We know that there is an isomorphism $\hat{X}^{\ell} \cong \lim_{Q,s,t} P^t \Tot_s \Resp (X/Q)$. Since $\pi_0$ commutes with limits and sends weak equivalences to homeomorphisms in $\hEh$, we conclude that we have an isomorphism 
$$\pi_n X \cong \lim_{Q,s,t} \pi_n (P^t \Tot_s \Resp (X/Q))$$ for every $n \geq 0$. The description of the homotopy groups of the simplicial finite sets 
$P^t \Tot_s\Resp (X/Q)$ finishes the proof. 
\end{proof} 

We conclude this subsection with a collection of results on the homotopy groups of the $\ell$-completion of a simplicial set, see also \cite{bouskan} VI, \S 5. 

\begin{prop}\label{piofZl-completion}
Let $X$ be a pointed connected simplicial set.\\ 
1. We have an isomorphism 
$$\pi_1(\widehat{X}^{\ell}) \cong \widehat{\pi_1(X)}^{\ell}.$$
In particular, if $\pi_1X$ is a finitely generated abelian group, then $\pi_1(\widehat{X}^{\ell}) \cong  \Z_{\ell} \otimes_{\Z} \pi_1X$.\\
2. If $X$ is in addition simply connected, and its higher homotopy groups are finitely generated, then we have for all $n \geq 2$
$$\pi_n \hat{X}^{\ell} \cong \widehat{\pi_n X}^{\ell} \cong \Z_{\ell} \otimes_{\Z} \pi_nX,$$
where $\Z_{\ell}$ denotes the $\ell$-adic integers. 
\end{prop}
\begin{proof} 
We deduce this proposition from the methods of \cite{bouskan}.
The hypothesis on $X$ implies that $X$ is a nilpotent space in the sense of \cite{bouskan}, i.e. $\pi_1X$ acts nilpotently on all $\pi_iX$. We denote by $\Zl_{\infty}X$ the Bousfield-Kan 
$\Zl$-completion in $\Sh$ of $X$. The assumptions imply that $X$ is a $\Zl$-good space, i.e. $H^{\ast}(X;\Zl)  \cong H^{\ast}(\Zl_{\infty}X;\Zl)$, cf. \cite{bouskan} VI, \S 5, Proposition 5.3. Considering the isomorphism $H^{\ast}(X;\Zl)  \cong H^{\ast}(\hat{X};\Zl) \cong H^{\ast}(|\hat{X}^{\ell}|;\Zl)$, this implies that $|\hat{X}^{\ell}|$ is $\Zl$-complete by \cite{bouskan} V, \S 3, Proposition 3.3. Then Proposition 5.4 of \cite{bouskan}, VI, \S 5, says that the two $\ell$-completions $|\hat{X}^{\ell}|$ and $\Zl_{\infty}X$ are weakly equivalent. Then \cite{bouskan} VI, \S 5, Proposition 5.2 yields the result.\\
Alternatively, one could deduce the result from  \cite{sullivan}, Theorem 3.1, and the fact that $|\widehat{X}^{\ell}|$ is isomorphic in $\Ho(\Sh)$ to the $\ell$-completion of Sullivan.  
\end{proof}

\begin{example}\label{pi1s1}
In the special case of the simplicial circle $S^1$, viewed as a profinite space, one has 
$$\pi_1(S^1)= \pi_1(\widehat{S^1}^{\ell})\cong \widehat{\Z}^{\ell}=\Z_{\ell}.$$
\end{example}
The last proposition is crucial for our studies of cohomology theories on $\hSh$, since it will allow us to compute the coefficients for some cohomology theories, e.g. for profinite cobordism $\widehat{MU}$.

\subsection{Comparison with the category of pro-simplicial sets}

Let $\pro-\Sh$ be the category of pro-objects in $\Sh$. Its objects are cofiltered diagrams $X(-):I \to \Sh$ and its morphisms are defined by 
$$\Hom_{\pro-\Sh}(X(-),Y(-)):=\lim_{t\in J}\colim_{s\in I} \Hom_{\Sh}(X(s),Y(t)),$$
see \cite{artinmazur} for more details. The cohomology with $\Zl$-coefficients is defined to be
$$H^{\ast}(X(-),\Zl):=\colim_s H^{\ast}(X(s),\Zl).$$ 
Isaksen has constructed several model structures on $\pro-\Sh$, cf. \cite{modstruc} and \cite{compofpro}. We are interested in the $\Zl$-cohomological model structure of \cite{compofpro} in which the weak equivalences are morphisms inducing isomorphisms in $\Zl$-cohomology and the cofibrations are levelwise monomorphisms. \\
We want to compare $\pro-\Sh$ with $\hSh$. These two categories are not equivalent, but we will see that completion induces an equivalence on the level of homotopy categories. We define a completion functor $\compl:\pro-\Sh \to \hSh$ as the composite of two functors. First we apply $\compl :\Sh \to \hSh$ levelwise to get a functor $\pro-\Sh \to \pro-\hSh$. Then we take the limit in $\hSh$ of the underlying diagram. The next lemma shows that the functor $\compl$ has good properties with respect to the model structures on $\pro-\Sh$ and $\hSh$. 

\begin{lemma}\label{Z/ell-comparison}
1. Let $X \in \pro-\Sh$ be a pro-simplicial set. Then we have a natural isomorphism of cohomology groups
$$H^{\ast}(X;\Zl) \stackrel{\cong}{\longrightarrow} H^{\ast}(\hat{X};\Zl).$$
2. A morphism $f:X\to Y$ of pro-simplicial sets induces an isomorphism in $\Zl$-cohomology if and 
only if the morphism $\hat{f}:\hat{X}\to \hat{Y}$ in $\hSh$  induces an isomorphism in continuous $\Zl$-cohomology. \\
3. The functor $\compl :\pro-\Sh \to \hSh$ preserves monomorphisms.\\
4. The functor $\compl :\pro-\Sh \to \hSh$ preserves finite limits and finite colimits.
\end{lemma}
\begin{proof}
1. This follows from the definition of $\compl$ and Remark \ref{remarkcohom}.\\
2. Suppose $f:X\to Y$ in $\pro-\Sh$  induces an isomorphism in $\Zl$-cohomology. 
We have the natural sequence of commutative squares
$$\begin{array}{ccc}
H^{\ast}(X;\Zl) & \stackrel{f^{\ast}}{\longrightarrow} & H^{\ast}(Y;\Zl)\\
=\downarrow &  & \downarrow =\\
\colim_s H^{\ast}(X_s;\Zl) & \longrightarrow & \colim_t H^{\ast}(Y_t;\Zl)\\
\cong \downarrow &   & \downarrow \cong \\
\colim_s H^{\ast}(\widehat{X_s};\Zl) & \longrightarrow & \colim_t H^{\ast}(\widehat{Y_t};\Zl)\\
\cong \downarrow &  & \downarrow \cong \\
H^{\ast}(\hat{X};\Zl) & \stackrel{\hat{f}^{\ast}}{\longrightarrow} & H^{\ast}(\hat{Y};\Zl)
\end{array}$$
which proves the second assertion.\\
3. This is clear.\\
4. Since  $\compl :\Sh \to \hSh$ is a left adjoint functor it preserves all colimits
and it commutes with finite limits. The functor $\lim:\pro-\hSh \to \hSh$ commutes with all limits and finite 
colimits. Hence the last statement holds, too. 
\end{proof} 

\begin{cor}\label{pro-comparison}
The functor $\compl :\pro-\Sh \to \hSh$ induces an equivalence of homotopy categories
$$\Ho(\pro-\Sh) \stackrel{\sim}{\longrightarrow} \Ho(\hSh),$$
where the left hand side is the $\Zl$-cohomological homotopy category of \cite{compofpro}.
\end{cor}
Isaksen has considered a different functor $\pro-\Sh \to \hSh$ in \cite{compofpro} and has shown that it induces an equivalence on the level of homotopy categories, too. Our approach fits better in the later picture for the comparison of generalized cohomology theories.

\subsection{The model structure on $\hSh$ is fibrantly generated}

For our studies of generalized cohomology theories and a stable structure on profinite spaces, we have to show the technical point that the model structure on $\hSh$, hence also on $\hShp$, is fibrantly generated. The necessity of this point becomes clear from the localization and stabilization results of the appendix, cf. Theorem \ref{leftloc}.\\
We recall some notations and constructions from \cite{ensprofin}. Let $n \geq 0$ be a non-negative integer and $S$ be a profinite set. The functor $\hSh^{\mathrm{op}}\to \Eh$, $X \mapsto \Hom_{\hEh}(X_n,S)$ is representable by a simplicial profinite set, which we denote by $L(S,n)$. It is given by the formula
$$L(S,n):\Delta^{\mathrm{op}} \to \hat{\Eh} , ~[k] \mapsto S^{\Hom_{\Delta}([n],[k])}.$$
If $M$ is a profinite abelian group, then $L(M,n)$ has a natural structure of a simplicial profinite abelian group and the abelian group $\Hom_{\hSh}(X,L(M,n))$ can be identified with the group $C^n(X;M)$ of continuous $n$-cochains of $X$ with values in $M$. Furthermore, for every $k$ $L(M,\ast)([k])$ may be considered in the usual way as an abelian cochain complex. 
For a profinite abelian group $M$, let $Z^n(X;M)$ denote the abelian group of $n$-cocycles of the complex $C^{\ast}(X;M)$. The functor $\hSh^{\mathrm{op}} \to \Eh$, $X \mapsto Z^n(X;M)$ is also representable by a simplicial profinite abelian group, which we denote by $K(M,n)$, called the profinite Eilenberg-MacLane space of type $(M,n)$. We may define $K(M,n)([k])$ to be the subgroup of all cocycles of $L(M,n)([k])$; see also \cite{may} \S 23, for these constructions.
\\
The natural homomorphism $C^n(X;M)\to Z^{n+1}(X;M)$ given by the differential defines a natural morphism of simplicial profinite abelian groups $L(M,n) \to K(M,n+1)$. \\
Consider the two sets of morphisms
$$P:=\{L(M,n)\to K(M,n+1), ~ K(M,n)\to \ast | M ~ \mathrm{abelian~ pro-\ell-group}, n\geq 0 \}$$ 
and
$$Q:=\{L(M,n)\to \ast |  M~\mathrm{abelian~ pro-\ell-group}, n\geq 0\}$$
of Lemme 2 of \cite{ensprofin}.

\begin{theorem}\label{fibgen}
The simplicial model structure on $\hSh$, in which the weak equivalences are the 
$\Z/\ell$-cohomological isomorphisms and the cofibrations are the dimensionwise monomorphisms, is left proper and fibrantly generated with $P$ as the set of generating fibrations and $Q$ as the set of generating trivial fibrations.
\end{theorem}
\begin{proof}
The left properness follows from the fact that all objects in $\hSh$ are cofibrant, cf. Corollary 13.1.3 of \cite{hirsch}.\\
We write Fib and Cof for the classes of fibrations and cofibrations of the model structure on $\hSh$ of \cite{ensprofin}. Note that the subcategory $Q$-cocell of relative $Q$-cocell complexes consists here of limits of pullbacks of elements of $Q$ and look at Definition 2.1.7 of \cite{hoveybook} for the notations $P-\proj$ and $Q-\fib$ used below. In order to prove the second statement we check the six conditions of the dual of Theorem 2.1.19 of \cite{hoveybook}. We only assume that the category has all small limits but only finite colimits. But since we use cosmall instead of small objects, the theorem holds for this kind of fibrantly generated model categories as well. We check now the six conditions of \cite{hoveybook}, Theorem 2.1.19.\\
1. It is clear that the weak equivalences satisfy the two-out-of-three axioms.\\
2. and 3. We have to show that the codomains of the maps in $P$ and $Q$ are cosmall relative to $P$-cocell and $Q$-cocell, respectively. Note that the terminal object $\ast$ is cosmall relative to all maps. Hence it remains to show that the objects $K(M,n)$ are cosmall relative to $Q$-cocell. It suffices to show that they are cosmall relative to a filtered sequence of maps $\ldots \to Y_1 \to Y_0$ of maps in $\hSh$. Consider the canonical map $f:\colim_{\alpha} \Hom_{\hSh}(Y_{\alpha},K(M,n)) \to \Hom_{\hSh}(\lim_{\alpha}Y_{\alpha},K(M,n)).$ We have to show that it is an isomorphism. By the definition of the spaces $K(M,n)$ this map is equal to the map
$$\colim_{\alpha} Z^n(Y_{\alpha},M) \to Z^n(\lim_{\alpha}Y_{\alpha},M).$$
But this map is already an isomorphism on the level of complexes
$\colim_{\alpha} C^n(Y_{\alpha},M) \cong C^n(\lim_{\alpha}Y_{\alpha},M)$ as in\cite{serre}, Proposition 8; hence $f$ is an isomorphism.\\
4. We know that $L(M,n)\to \ast$ is a trivial fibration. Since trivial fibrations are
preserved under pullbacks and limits, we get
$Q-\mathrm{cocell} \subseteq W \cap \mathrm{Fib} \subseteq W \cap P-\fib$ which is what we had to show.\\
5. Given a diagram
\begin{equation}\label{diagram1}
\begin{array}{ccc}
A & \to & L(M,n)\\
\downarrow &   & \downarrow\\
B & \to & \ast
\end{array}
\end{equation}
we have to show that there is a lift if the map $f:A\to B$ is in $P-\proj$, i.e. has the left lifting property with respect to $P$. But if $F$ is in $P-\proj$, then we get a lifting $B\to K(M,n+1)$ for any map $A\to K(M,n+1)$. Hence the diagram (\ref{diagram1})
yields a diagram
\begin{equation}\label{diagram2}
\begin{array}{ccc}
A & \to & L(M,n)\\
\downarrow &   & \downarrow\\
B & \to & K(M,n+1)\\
\parallel &  & \downarrow\\
B & \to & \ast
\end{array}
\end{equation}
and we know that there is a lifting $B\to L(M,n)$ in the upper rectangle which is also the lift of the diagram (\ref{diagram1}) above. Hence $f \in Q-\proj$.\\
It remains to show that $P-\proj \subseteq W$. Let $f:A\to B$ be a map in
$P-\proj$. By definition of the spaces $K(\Z/\ell,n)$ and the definition of $P-\proj$ we get that
$f^{\ast}:Z^n(B,\Z/\ell)\to Z^n(A,\Z/\ell)$ is surjective for all $n\geq 0$.
Hence it is enough to show that $f^{\ast}(\Imm(C^nA\to Z^{n+1}A))\subseteq \Imm(C^nB\to Z^{n+1}B).$
The other lifting property of maps in $P-\proj$ is equivalent to the surjectivity
of the map $C^nB\to C^nA\times_{Z^{n+1}A}Z^{n+1}B$ from which we get the
desired result that $f^{\ast}:H^n(B,\Z/\ell)\to H^n(A,\Z/\ell)$  is an isomorphism for all $n\geq0$.
\\
6. We show that $W\cap Q-\proj \subseteq P-\proj$. Let $f:A\to B$ be a map that belongs to $W$ and $Q-\proj$. Since $f^{\ast}$ is an isomorphism on $\Z/\ell$-cohomology, it can be shown that $f$ induces an isomorphism on cohomology with coefficients in any abelian pro-$\ell$-group $M$. This implies already the surjectivity of the map $f^{\ast}:Z^n(B,M) \to Z^n(A,M)$ for all $n\geq 0$.
The lifting property of maps in $Q-\proj$ implies the surjectivity of the induced map
$f^{\ast}:C^n(B,M) \to C^n(A,M)$ for all $n\geq 0$. Using the isomorphism on cohomology
it follows that $C^nB \to C^nA\times_{Z^{n+1}A}Z^{n+1}B$ is surjective which is equivalent
to the other desired lifting property of maps in $P-\proj$.
Hence $f\in P-\proj$.
\\
Now we have proved using theorem 2.1.19 of \cite{hoveybook} that the quadruple
$(\hSh, W, Q-\proj, P-\fib)$ is a fibrantly generated model category.
It remains to show that it coincides
with the given structure ($\hSh$, $W$, Cof, Fib).
Since $L(M,n)\to \ast$ is a trivial fibration, we know already Cof $\subseteq Q-\proj$.
It remains to show that every map in $Q-\proj$ is a monomorphism in each dimension. Then the cofibrations and weak equivalences of the two structures coincide and hence the fibrations coincide as well.
\\
Let $f:X\to Y$ be a map in $Q-\proj$. Hence every map $X\to L(M,n)$ can be lifted to a map $Y\to L(M,n)$. This means by using the definition of the spaces $L(M,n)$ that the map $C^n(Y,M)\to C^n(X,M)$ induced by $f$ is surjective for all $n\geq 0$. This is equivalent to the surjectivity of the maps $\Hom_{\hEh}(Y_n,M)\to \Hom_{\hEh}(X_n,M)$ for all $n\geq 0$ and all abelian pro-$\ell$-groups $M$.
If we choose $M=F(X_n)$ to be the free abelian pro-$\ell$-group on the set $X_n$ defined
in \cite{serre}, then we see that the map $X_n \to F(X_n)$, sending $x$ to its class in $F(X_n)$, is in the image of $\Hom_{\hEh}(Y_n,M)\to \Hom_{\hEh}(X_n,M)$ only if $X_n \hookrightarrow Y_n$ is
a monomorphism for each $n$. Hence $f$ is also a cofibration.
\end{proof}

\section{Profinite spectra}

We introduce the usual notion of spectra on $\hSh$. The results of the appendix show that there is a stable model structure on $\hSp$. This is the main technical result for the construction of generalized cohomology theories on schemes via the profinite \'etale topological type functor.\\
We study the behavior of the profinite completion functor $\compl:\Sp \to \hSp$ and show that it factorizes through stable equivalences and is in fact a left Quillen functor with right adjoint the forgetful functor.\\
The following subsection is dedicated to stable homotopy groups and their relation to stable equivalences. We show that these homotopy groups behave well under completion of spaces. In particular, we study them for profinite Eilenberg-MacLane spectra and the completed complex cobordism spectrum $\hMU$. In this last case, we get $\pi_{\ast}(\hMU) \cong \Z_{\ell} \otimes \pi_{\ast}(MU)$. \\
Furthermore, we give an explicit construction of a fibrant replacement functor for spectra satisfying some conditions among which is $MU$. These results are similar to those of \cite{dehon}.\\  
At the end of this section we deduce the existence of Postnikov-decompositions for profinite spectra from general model category theory.

\subsection{The stable structure of profinite spectra}

In order to stabilize the category $\hSh$ of profinite spaces we begin with the usual
definition of spectra.
\begin{defn}\label{defspectrum}
A {\rm spectrum} $X$ of simplicial profinite sets consists of a sequence $X_n \in \hShp$ of pointed
profinite spaces for $n\geq0$ and maps $\sigma_n:S^1 \wedge X_n \to X_{n+1}$ in $\hShp$, where
$S^1 = \Delta^1/\partial \Delta^1 \in \hShp$.\\
A {\rm map} $f:X \to Y$ of spectra consists of maps $f_n:X_n \to Y_n$ in $\hShp$ for $n\geq0$
such that $\sigma_n(1\wedge f_n)=f_{n+1}\sigma_n$.\\
We denote by $\hSp$ the corresponding category and call it the category of {\rm profinite spectra}.
\end{defn}

We apply the results of Appendix B to the category $\hShp$. By Theorem \ref{fibgen} its model structure is simplicial, left proper and fibrantly generated. Note that $S^1\wedge \cdot$ is a left Quillen endofunctor since it takes monomorphisms to monomorphisms and preserves cohomological equivalences by Lemma \ref{smashKunneth}.

\begin{theorem}\label{thmprostable}
There is a {\rm stable model structure on $\hSp$} for which the prolongation $S^1 \wedge \cdot:\hSp \to \hSp$ is a Quillen equivalence.\\
In particular, the {\rm stable equivalences} are the maps that induce an isomorphism on
all generalized cohomology theories, represented by profinite $\hat{\Omega}$-spectra;
the {\rm stable cofibrations} are the maps $i:A \to B$ such that $i_0$ and the induced maps $j_n: A_n \amalg_{S^1\wedge A_{n-1}}S^1 \wedge B_{n-1} \to B_n$ are monomorphisms;
the {\rm stable fibrations} are the maps with the right lifting property
with respect to all maps that are both stable equivalences and stable cofibrations.
\end{theorem}
\begin{proof}
By Theorem \ref{fibgen}, we know that the model structure on $\hShp$ is fibrantly generated, left proper and simplicial. Hence we can apply Theorem \ref{stablestructure} to $\hShp$, $S^1 \wedge \cdot$. The fact that $S^1 \wedge \cdot$ is a Quillen equivalence for this model structure is implied by Theorem \ref{Tstable} for $T=S^1 \wedge \cdot$.
\end{proof}

\begin{remark}
The stable fibrant objects are exactly the $\hat{\Omega}$-spectra, i.e. spectra $E$ such that each
$E_n$ is a fibrant object and the adjoint structure maps $E_n \to \hat{\Omega}E_{n+1}$ are weak equivalences in $\hShp$ for all $n\geq 0$.
\end{remark}

Define the functor $\hat{\Omega}:\hSp \to \hSp, E \mapsto \hat{\Omega} (E)$, where $\hat{\Omega}(E)_n :=\hat{\Omega}(E_n)$.
As usual one can prove the following lemma.
\begin{lemma}\label{additivity}
1. The functor $$\hSHh \to \hSHh, ~ E \mapsto S^1 \wedge E$$ is an equivalence of categories.\\
2. The homotopy category $\hSHh$ of profinite spectra is an additive category.
\end{lemma}
\begin{proof}
Choose a functorial stable replacement
$$E \mapsto E^f.$$
Then $E \mapsto \hat{\Omega}(E^f)$ is an inverse to $E\mapsto S^1 \wedge E$ using the two stable equivalences
$$E \to \hat{\Omega}((S^1 \wedge E)^f)$$
and
$$S^1 \wedge \hat{\Omega}(E^f) \to E^f$$
which are induced by the adjunction morphisms.\\
As $S^1$ is a cogroup object in $\mathrm{Ho}(\hShp)$, its codiagonal $\phi:S^1 \to S^1 \vee S^1$ induces, via first assertion of the lemma and $\phi \wedge \id_E:S^1 \wedge E \to S^1 \wedge E \vee S^1 \wedge E$, for any spectrum $E$ a morphism $E\to E\vee E$. This yields a natural structure of a cogroup object on $E$. Being natural in $E$, this structure has to be abelian. Thus the category $\hSHh$ is additive.
\end{proof}

\begin{remark}\label{triangular}
Since $\hSp$ is a pointed model category we may conclude by the general theory of \cite{homalg} I, \S 2 and \S 3, that $\hSp$ has fiber and cofiber sequences, which satisfy the axioms of a triangulation. 
Since the above theorem tells us that $\Sigma$ is an equivalence on $\hSHh$, we may call $\hSp$ 
a {\em stable} category. 
\end{remark} 

\subsection{Profinite completion of spectra}

\begin{lemma}\label{smash}
The functor $\compl:\Shp \to \hShp$ commutes with smash products.\\
Furthermore, the two functors $\Omega$ and $\hat{\Omega}$ agree after applying the forgetful functor
$|\cdot|:\hShp \to \Shp$, i.e. $\Omega(|Z|)\cong |\hat{\Omega}(Z)|$ for every pointed profinite space
$Z\in \hShp$.
\end{lemma}
\begin{proof}
Since $\compl$ is left adjoint to the forgetful functor it preserves colimits. Hence it remains to check that it commutes with finite products. But this follows from the fact that the finite quotients of a product $X\times Y$ of two spaces are in bijective correspondence with pairs of finite quotients of $X$ and $Y$ respectively. Since $\compl$ is defined by applying a limit we see immediately that it commutes with finite products.
\\
For every $X\in \Shp$ we have natural bijections coming from adjunction
$$\begin{array}{rcccl}
\Hom_{\Sh}(S^1 \wedge X, |\hat{\Omega}Z|) & \cong & \Hom_{\hSh}(\hat{X},\hat{\Omega}Z) &
\cong & \Hom_{\hSh}(S^1 \wedge \hat{X},Z) \\
  & \cong & \Hom_{\hSh}(\widehat{S^1\wedge X},Z) & \cong & \Hom_{\Sh}(S^1\wedge X,|Z|) \\
  & \cong & \Hom_{\Sh}(X,\Omega |Z|).
\end{array}$$
Hence the two Hom-functors defined by $\Omega(|Z|)$ and $|\hat{\Omega}(Z)|$ agree
on $\Sh$ and hence by Yoneda the two spaces $\Omega(|Z|)$ and $|\hat{\Omega}(Z)|$ are
isomorphic for every profinite space $Z$.
\end{proof}

Let $\Sp$ be the stable model structure of simplicial spectra defined in \cite{bousfried}. We may extend the profinite completion to spectra. Let $\compl :\Sp \to \hSp$ be the profinite completion  applied levelwise that takes the spectrum $X$ to the profinite spectrum $\hat{X}$ whose structure maps are given, using Lemma \ref{smash}, by
$$S^1 \wedge \hat{X}_n \cong \widehat{S^1 \wedge X_n} \stackrel{\hat{\sigma}}{\to} \hat{X}_{n+1}.$$
Since the two functors $\compl$ and $|\cdot|$ on $\hShp$ commute with smash products we get
that commutative diagrams
$$\begin{array}{ccc}
X_{n+1} & \to & |Y_{n+1}| \\
\uparrow &   & \uparrow \\
S^1 \wedge X_n & \to & S^1 \wedge |Y_n|
\end{array}$$
are in bijective correspondence with commutative diagrams
$$\begin{array}{ccc}
\hat{X}_{n+1} & \to & Y_{n+1} \\
\uparrow &   & \uparrow \\
S^1 \wedge \hat{X}_n & \to & S^1 \wedge Y_n
\end{array}$$
for every $X\in \Sp$ and $Y\in \hSp$. Let $|\cdot|: \hSp \to \Sp$ be the levelwise applied forgetful functor. Hence $\compl$ and $| \cdot |$ form an adjoint pair of functors. 

\begin{prop}\label{quillenpair}
The functor $\compl :\Sp \to \hSp$ preserves weak equivalences and cofibrations.\\
The functor $|\cdot|: \hSp \to \Sp$ preserves fibrations and weak equivalences between
fibrant objects.\\
In particular, $\compl$ induces a functor on the homotopy categories and the adjoint pair $(\compl , |\cdot|)$ is a Quillen pair of adjoint functors.
\end{prop}

\begin{proof}
Let $i:A \to B$ be a cofibration in $\Sp$. Since $\compl:\Sh_{\ast} \to \hShp$ preserves cofibrations and pushouts as a left Quillen functor, the maps $\hat{i_0}$ and $\hat{j_n}$ are cofibrations in $\hShp$. Hence $\hat{i}$ is a cofibration in $\hSp$. 
\\
Now let $f:X \to Y$ be a weak equivalence in $\Sp$. We want to show that $\hat{f}$ is a weak equivalence in $\hSp$. By \cite{hirsch}, Theorem 9.7.4, this is equivalent to the statement that $\Map(\tilde{\hat{f}},E):\Map(\tilde{\hat{Y}},E) \to \Map(\tilde{\hat{X}},E)$ is a weak equivalence of simplicial sets for every fibrant object $E$  and some cofibrant approximation $\tilde{\hat{f}}$ of $\hat{f}$ in $\hSp$. Since $\compl$ and $X\wedge \cdot$ are left adjoints, we have natural isomorphisms extending the adjunction of $\compl$ and $|\cdot|$ for every $U\in \Sp$ and $V\in \hSp$
\begin{equation}\label{Mapadjointness}\Map_{\Sp}(U,|V|) \cong \Map_{\hSp}(\hat{U},V).\end{equation}
Hence $\Map_{\hSp}(\hat{f},E)$ is a weak equivalence of simplicial sets if and only if
$\Map_{\Sp}(f,|E|)$ is a weak equivalence. Since $\compl$ preserves cofibrations, the completion of
$\tilde{f}$ is a cofibrant approximation of $\hat{f}$.
We conclude the first part by \cite{hirsch}, Theorem 9.7.4, and the following lemma,
reminding the fact that the fibrant objects of $\Sp$ are exactly the
$\Omega$-spectra.

\begin{lemma}\label{omegaspectra}
If $E$ is an $\hat{\Omega}$-spectrum, then $|E|$ is an $\Omega$-spectrum.
\end{lemma}
This is the analogue of the fact that the underlying simplicial set $|X|$ of a fibrant profinite space is a Kan simplicial set, cf. \cite{ensprofin}, \S 2.1, Proposition 1.

\begin{proof}
Let $E$ be an $\hat{\Omega}$-spectrum. This implies that each $E_n$ is fibrant and $E_n \to \hat{\Omega}E_{n+1}$
is a weak equivalence for each $n\geq 0$. By \cite{ensprofin}, $\S$ 2, Proposition 1, this implies
that $|E_n|$ is fibrant and $|E_n| \to \Omega |E_{n+1}|$ is a weak equivalence for each $n$.
Hence $|E|$ is an $\Omega$-spectrum.
\end{proof}

We continue the proof of the second statement of the proposition. The fact that $|\cdot|$ preserves fibrations follows now from adjunction since $\compl$ preserves trivial cofibrations. Now let $f:E \to F$ be a weak equivalence between fibrant objects in $\hSp$. Again by \cite{hirsch}, Theorem 9.7.4, we have to show that $\Map_{\Sp}(W,|f|)$ is a weak equivalence for every cofibrant object $W$ of $\Sp$. By the isomorphism (\ref{Mapadjointness}) this is equivalent to that $\Map_{\hSp}(\hat{W},f)$ is a weak equivalence for every such $W$. But since $\hat{W}$ is also cofibrant in $\hSp$ and since $f$ is a weak equivalence in $\hSp$, we get that $\Map_{\hSp}(\hat{W},f)$ is a weak equivalence of simplicial sets.\\
The last statement follows from general model category theory.
\end{proof}

\subsection{Stable homotopy groups and fibrant replacements}

For a level fibrant replacement functor of a profinite spectrum $E$ we may use the explicit construction of the $\ell$-completion functor $X\mapsto \widehat{X}^{\ell}$ in $\hShp$. Since all construction in \cite{bouskan} and \cite{ensprofin}
are compatible with products there are natural maps
$$X\times \widehat{Y}^{\ell} \to \widehat{X}^{\ell} \times \widehat{Y}^{\ell} \stackrel{\sim}{\to} \widehat{X\times Y}^{\ell}$$
and similarly maps
$$X\wedge \widehat{Y}^{\ell} \to \widehat{X}^{\ell} \wedge \widehat{Y}^{\ell} \stackrel{\sim}{\to} \widehat{X\wedge Y}^{\ell}.$$
For a simplicial spectrum $E$ this yields structure maps
$$S^1 \wedge \widehat{E_n}^{\ell} \to (\widehat{S^1 \wedge E_n})^{\ell} \stackrel{\widehat{\sigma_n}^{\ell}}{\to} \widehat{E_{n+1}}^{\ell},$$
which ensure that we may associate to every profinite spectrum $E$ a level fibrant spectrum $E^{lf}$
such that $E\to E^{lf}$ is a level equivalence.

We may also define homotopy groups of profinite spectra. For $n \in \Z$ we set
\begin{equation}\label{defofpi}
\pi_n(E):= \pi_n (E^{lf})=\colim_k \pi_{n+k}(E^{lf}_k).
\end{equation}
For morphisms of spectra we set $\pi_n (g):=\pi_n(g^{lf})$.

\begin{remark}\label{profinitestablepi}
In order to calculate the stable homotopy groups of a spectrum $E$, we may replace it by a fibrant spectrum $E^f$. Then its homotopy groups are the homotopy groups of the infinite loop space $E^f_0$. By Proposition \ref{pro-l-group} we know that the homotopy groups of a profinite space have a natural structure of pro-$\ell$-groups. Hence we conclude that the stable homotopy groups of any profinite spectrum have the structure of pro-$\ell$-groups. 
\end{remark}

We may describe stable equivalences via stable homotopy groups.
\begin{prop}\label{stablepiequiv}
A map of profinite spectra $g:E\to F$ is a stable equivalence if and only if
$\pi_n(g)$ is an isomorphism for all $n \in \Z$.
\end{prop}
\begin{proof}
Since we are interested in the homotopy type of the map, we may consider a stable fibrant replacement $g^f: E^f \to F^f$ of $g$. The map $g$ is a stable equivalence if and only if $g^f$ is a stable equivalence. 
But a map between $\hOmega$-spectra is a stable equivalence if and only if it is a level equivalence, e.g. see \cite{hovey}. But this means that all maps $\pi_n(g^f_k)$ are isomorphisms for all $n$ and $k$, which finishes the proof by Proposition \ref{piequiv}.
\end{proof}

\begin{cor}\label{reflectswe}
Let $f:X \to Y$ be a map in $\hSp$ such that $|f|:|X| \to |Y|$ is a stable equivalence
of $\Sp$. Then $f$ is a stable equivalence of $\hSp$.
\end{cor}

\begin{proof}
This follows immediately from the proposition above.
\end{proof}

From model category theory we know that fibrant replacements with respect to the stable model structure exist for every $E \in \hSp$. We would like to construct explicit stable fibrant replacements in $\hSp$. We will employ the Bousfield-Kan-$\ell$-completion functor defined on profinite spaces. Since the $\ell$-completion behaves well only under certain conditions, cf. \cite{bouskan}, we give a construction only for profinite spectra satisfying the following hypotheses.\\  
The idea for the following construction is inspired by Dehon \cite{dehon}. Let $E$ be a $(-1)$-connected profinite spectrum, i.e. each profinite space $E_n$ is $(n-1)$-connected for all $n \geq 0$ and let $|E| \in \Sp$ be its underlying simplicial spectrum. We may consider an $\Omega$-spectrum $R(|E|)$ stable equivalent to $|E|$, i.e. a stable fibrant replacement of $|E|$ in $\Sp$. After the application of the levelwise $\ell$-completion we get maps 
$\widehat{(R|E|)_n}^{\ell} \to \widehat{\Omega (R|E|)_{n+1}}^{\ell}$ that are weak equivalences in $\hShp$ for all $n$. Then Proposition \ref{omegal} implies that all the composite maps $$\widehat{(R|E|)_n}^{\ell} \to \widehat{(\Omega R|E|_{n+1})}^{\ell} \to  \hOmega \widehat{(R|E|)_{n+1}}^{\ell}$$ are weak equivalences in $\hShp$ for all $n \geq 1$. It is clear that the resulting profinite spectrum is also level fibrant. From Corollary \ref{reflectswe}, we deduce that it is also stable equivalent to $E$. Hence we have constructed an explicit stable fibrant replacement of $E$ in $\hSp$.

\begin{defnprop}\label{stablefibrantreplacement}
Let $E$ be a $(-1)$-connected profinite spectrum. We suppose that the $\Zl$-cohomology of each $E_n$, $n \geq 1$, is finite dimensional in each degree. We consider the profinite  spectrum $\hat{E}^{\ell} \in \hSp$ whose spaces are defined by $\hat{E}^{\ell}_0 :=\hOmega \widehat{(R|E|)_1}^{\ell}$ for $n=0$ and by $\widehat {(R|E|)_n}^{\ell}$ for all $n\geq 1$ with structure maps adjoint to the above defined maps $\widehat{(R|E|)_n}^{\ell} \to \hOmega \widehat{(R|E|)_{n+1}}^{\ell}$ for $n \geq 1$ and the obvious map for $n=0$. By the preceding discussion, the profinite spectrum $\hat{E}^{\ell}$ is a {\rm stable fibrant replacement} of $E$ in $\hSp$ and we call it the {\rm $\ell$-completion of $E$}. 
\end{defnprop}
The reason why we use this definition for $\hat{E}^{\ell}_0$ is that if $E_0$ has infinitely many connected components, the natural structure map will not be a weak equivalence any more.\\
This discussion also leads to an $\ell$-completion of spectra $E \in\Sp$ which is inspired by \cite{dehon}, \S 1.6.
\begin{defn}\label{lcompletionofspectra}
Let $E$ be a $(-1)$-connected spectrum in $\Sp$. We suppose that the $\Zl$-cohomology of each $E_n$, $n \geq 1$, is finite dimensional in each degree. Let $RE$ be an equivalent $\Omega$-spectrum. Then we define the {\rm profinite $\ell$-completion} $\hat{E}^{\ell}$ to be the profinite spectrum whose spaces are defined by $\hat{E}^{\ell}_0:=\hOmega \widehat{(R(E)_1)}^{\ell}$ for all $n=0$ and by $\widehat {(R(E)_n)}^{\ell}$ for all $n\geq 1$ with structure maps adjoint to the above defined maps $\widehat{R(E)_n}^{\ell} \to \hOmega \widehat{R(E)_n}^{\ell}$ for $n \geq 1$ and the obvious map for $n=0$. 
\end{defn}

\begin{remark}
As above, for $E\in \Sp$, it is clear that $\hat{E}^{\ell}$ is stable equivalent to $\hat{E}$ in $\hSp$. From the homotopy theoretic point of view we may consider either $\hat{E}^{\ell}$ or $\hat{E}$. The point is that we want to apply the comparison results of Proposition \ref{piofZl-completion} to the situation of spectra. The previous definition gives us the form of the $\hat{E}$ that fits well with the results in Proposition \ref{piofZl-completion}.
\end{remark}

Similar versions of the following facts have also been proved by Dehon \cite{dehon}.
\begin{prop}\label{stablepiofcompletion}
Let $E$ be a $(-1)$-connected spectrum and suppose that the $\Zl$-cohomology of each $E_n$, $n \geq 1$, is finite dimensional in each degree. Then the stable homotopy groups of the profinite $\ell$-completion $\hat{E}^{\ell}$ are given by the following isomorphism for all $n$
$$\pi_n \hat{E}^{\ell} \cong \widehat{\pi_n E}^{\ell} \cong \Z_{\ell} \otimes_{\Z} \pi_n E.$$
\end{prop}

\begin{proof}
Since $\hat{E}^{\ell}$ is an $\hOmega$-spectrum, the stable homotopy group $\pi_n \hat{E}^{\ell}$ is equal to the homotopy group $\pi_{n+1} ((\hat{E}^{\ell})_{1})$ for all $n$. The hypothesis on $E$ and the definition of $\hat{E}^{\ell}$ imply that the spaces $\hat{E}^{\ell}_n$ satisfy the conditions of Proposition \ref{piofZl-completion}. This completes the proof.
\end{proof}

\begin{cor}\label{piofMU}
Let $MU$ be the simplicial spectrum representing complex cobordism. For the profinite $\ell$-completion $\hMU^{\ell}$ of $MU$ there is an isomorphism
$$\pi_{\ast}\hMU =\pi_{\ast}\hMU^{\ell} \cong  \Z_{\ell} \otimes_{\Z} \pi_{\ast} MU \cong \Z_{\ell} \otimes_{\Z} \Lee_{\ast}$$
where $\Lee{\ast}$ denotes the Lazard ring, cf. \cite{adams}. 
\end{cor}
\begin{proof}
Since $MU$ is $(-1)$-connected and its cohomology groups are finitely generated in each degree, cf. \cite{adams}, $MU$ satisfies the hypothesis of Proposition \ref{stablepiofcompletion} above. 
\end{proof}

We finally consider the profinite completion $\hKU$ of the $\Omega$-spectrum representing
complex $K$-theory, see e.g \cite{may2}: 
$\hKU_{2i}= \widehat{BU\times  \Z}$ and $\hKU_{2i+1}= \hat{U}$
for all $i\geq 0$. By Lemma \ref{piomegahat} and Proposition \ref{piofZl-completion}
we get 
\begin{cor}\label{hatKU}
$\pi_{2i}(\hKU)\cong \Z_{\ell} ~ \mathrm{and} ~ \pi_{2i+1}(\hKU)= 0 ~ \mathrm{for ~ all} ~ i.$
\end{cor} 

\section{Generalized cohomology theories on profinite spaces}

We define generalized cohomology theories on $\hSh$ via the stable category of profinite spectra in the classical way. Such profinitely completed cohomology theories have already been studied by Dehon in \cite{dehon}. \\
It is the main advantage of the stable category of profinite spectra that it gives a canonical and general setting for cohomology theories on profinite spaces and the profinite completion of cohomology theories on simplicial sets. In fact, our main objective are the cohomology theories represented by the completion of simplicial spectra such as $\hMU$. The results of the previous section enable us to calculate the coefficients of this theory. \\
An important tool for the calculation of cohomology groups is the profinite Atiyah-Hirzebruch spectral sequence. The construction of this spectral sequence is the analog of the topological Atiyah-Hirzebruch spectral sequence. It has already been constructed by Dehon in a more restricted setting. We slightly generalize his proof.\\    
At the end we will prove a K\"unneth formula for $\hMU$ with $\Zln$-coefficients, similar to and inspired by the results of \cite{dehon}.   

\subsection{Generalized cohomology theories}

In analogy to the stable homotopy theory of simplicial spectra the above construction enables us to
consider generalized cohomology theories for profinite spaces. We define them to be the functors represented by profinite spectra. Let $E$ be a spectrum in $\hSp$, we set
\begin{equation}
E^n(X):= \Hom_{\hSHh}(X,E[n]),
\end{equation} 
where $\Hom_{\hSHh}(X,E[n])$ denotes the set of maps that lower the dimension by $n$, and call this the $n$-th cohomology group of $X$ with values in $E$.
We set $E^{\ast}(X):=\bigoplus_n E^n(X)$.
For a pointed profinite space $X$ we define its $n$-th cohomology groups for the spectrum $E$ by
\begin{equation}\label{defgencoho}
E^n(X):= \Hom_{\hSHh}(\Sigma ^\infty (X),E\wedge S^n)
\end{equation}
For a profinite space $X$ without a chosen basepoint, let $X_+$ be the profinite space $X\coprod pt$ with additional basepoint. We define the cohomology of $X$ to be the one of $X_+$.\\
For a pair $(X,A)$ of profinite spaces we define the relative cohomology by
\begin{equation}\label{defrelcoho}
E^n(X,A):= E^n(X/A).
\end{equation} 
We have the identity $E^{\ast}(X,pt)=E^{\ast}(X)$.
Just as for simplicial spectra \cite{switzer}, 8.21, one can prove the following isomorphism
\begin{equation}\label{colimrel}
\Hom_{\hSHh}(\Sigma ^\infty (X),E\wedge S^n) \cong
\colim_k \Hom_{Ho(\hShp)}(\Sigma^k (X), E_{n+k}).
\end{equation}
In particular, for an $\hat{\Omega}$-spectrum $E$, we get:
\begin{equation}\label{omegacohom}
E^n(X)\cong  \Hom_{Ho(\hShp)}(X,E_n).
\end{equation}

\subsection{Continuous cohomology theories}

\begin{defn}
Let $E$ be a spectrum in $\Sp$. We define the {\em continuous cohomology represented by $E$} to be the theory $\hat{E}^{\ast}(\cdot)$ on $\hShp$.
\end{defn}

\begin{remark}
1. Note that since profinite completion commutes with smash-products by Lemma \ref{smash}, the completion of a multiplicative cohomology theory is multiplicative, too.\\
2. Let $E$ be a $(-1)$-connected spectrum. We suppose that the $\Zl$-cohomology of each $E_n$, $n \geq 1$, is finite dimensional in each degree. Since the spectra $\hat{E}$ and $\hat{E}^{\ell}$ are equivalent in $\hSp$, the {\em continuous $\ell$-completed cohomology theory} on $\hShp$ represented by $\hat{E}^{\ell}$ is isomorphic to $\hat{E}$, i.e. for all $X \in \hShp$ we have 
$\hat{E}^{\ell \ast}(X)=\hat{E}^{\ast}(X)$.
\end{remark}

\begin{example}
By abuse of notations we write $MU$ for the simplicial spectrum  $\mathrm{Sing}(MU)$ and we will write $KU$ for $\mathrm{Sing}(KU)$. There are the obvious examples for profinite spectra representing generalized cohomology theories:\\
1. The Eilenberg-MacLane spectra $H\pi$ given by $H\pi_n =K(\pi,n)$ for an abelian pro-$\ell$-group $\pi$.\\
2. The profinite completion $\hKU$ of the simplicial spectrum representing complex K-theory.\\
3. The profinite completion $\hMU$ of the simplicial spectrum representing complex cobordism.
\end{example}

\begin{prop}\label{coeffofcompletion}
Let $E$ be a multiplicative $(-1)$-connected spectrum and suppose that the $\Zl$-cohomology of each $E_n$, $n \geq 1$, is finite dimensional in each degree. Then the coefficient ring 
$\hat{E}^{\ast}:= \hat{E}^{\ast}(pt)$ of the profinite completion satisfies the following isomorphism 
$$\hat{E}^{\ast}(pt) \cong \widehat{E^{\ast}(pt)}^{\ell} \cong \Z_{\ell} \otimes_{\Z} E^{\ast}(pt).$$
For the homotopy groups of $\hat{E}$ we get
$$\pi_{\ast}(\hat{E}) \cong \widehat{\pi_{\ast}(E)}^{\ell} \cong \Z_{\ell} \otimes_{\Z} \pi_{\ast}(E).$$
\end{prop}
\begin{proof}
This is a corollary of Proposition \ref{stablepiofcompletion}, since $\hat{E}$ and $\hat{E}^{\ell}$ are equivalent objects in $\hSp$.
\end{proof} 

Since the spectrum $MU$ satisfies the hypotheses of the previous proposition, we get the following result.
\begin{cor}\label{coeffofcompletedMU}
For the homotopy groups of $\hMU$ there is an isomorphism
$$\pi_{\ast}(\hMU) \cong \widehat{\pi_{\ast}(MU)}^{\ell} \cong \Z_{\ell} \otimes_{\Z} \pi_{\ast}(MU).$$
For the coefficient ring of the profinite cobordism $\hMU^{\ast}$ we have the following isomorphism
$$\hMU^{\ast} = \hMU^{\ell \ast} \cong \Z_{\ell} \otimes_{\Z} MU^{\ast}.$$
\end{cor}

Let $G$ be a finite group. There is a simplicial Moore spectrum $MG$. For a simplicial spectrum $E$ we define $EG:=E\wedge MG$ to be the spectrum with coefficients in $G$ corresponding to $E$, cf. \cite{adams}, Chapter III. 
\begin{defn}\label{cohomwithcoeff}
Let $E$ be a simplicial spectrum and $G$ a finite abelian group. We define the continuous cohomology theory with $G$-coefficients to be 
$$\hat{E}^{\ast}(\cdot ; G):=\widehat{EG}^{\ast}(\cdot).$$ 
\end{defn}
By \cite{adams} III, Prop. 6.6, there are exact sequences
\begin{equation}\label{torsequencepi}
0 \to \pi_n(E)\otimes G \to \pi_n(EG) \to \Tor_1^{\Z}(\pi_{n-1}(E),G) \to 0
\end{equation}
\begin{equation}\label{torsequencecohom}
0 \to E^n(X)\otimes G \to (EG)^n(X) \to \Tor_1^{\Z}(E^{n+1}(X),G) \to 0
\end{equation}
for all spectra $E$ and all spaces $X$.
 
\begin{cor}\label{hMUZl^k}
For the homotopy groups of the profinite cobordism with $\Zln$-coefficients $\hMU\Zln$ there is an isomorphism
$$\pi_{\ast}(\hMU\Zln) \cong \widehat{\pi_{\ast}(MU\Zln)}^{\ell} \cong \Zln \otimes_{\Z} \pi_{\ast}(MU).$$
For the coefficient ring of the profinite cobordism $\hMU^{\ast}$ we have the following isomorphism
$$(\hMU\Zln)^{\ast} = (\widehat{MU\Zln}^{\ell})^{\ast} \cong \Zln \otimes_{\Z} MU^{\ast}.$$
\end{cor}
\begin{proof}
Since $\pi_{\ast}(MU)$ and $MU^{\ast}$ have no torsion, the assertions follow from the previous exact sequences and Proposition \ref{coeffofcompletion}.
\end{proof}

For the Eilenberg-MacLane spectra recall the construction of the spaces $K(\pi,n)$ in the previous section. We have the following results. 
\begin{prop}\label{Kpin}
Let $\pi$ be an abelian pro-$\ell$-group.\\ 
1. The Eilenberg-MacLane spaces $K(\pi,n)$ represent continuous cohomology in $\hHh$, i.e. for every profinite space $X$  and every $n\geq 0$ there is an isomorphism $$H^n(X;\pi) \cong \Hom_{\hHh}(X, K(\pi,n)).$$
2. The Eilenberg-MacLane spectrum $H\pi=H(\pi)$ represents continuous cohomology with coefficients in $\pi$, i.e. $[\Sigma ^{\infty}(X),H\pi[n]]\cong H^n(X;\pi)$ for every $n$ and every $X\in \hSh$. \\
3. Let $Y$ be a simplicial set. We denote by $|\pi|$ the underlying abstract group. There is a natural isomorphism $$H^{\ast}(Y;|\pi|) \cong H^{\ast}(\hat{Y};\pi).$$ 
In particular, if $G$ is a finitely generated abelian group, there are isomorphisms
\begin{equation}\label{ZlG} 
H^{\ast}(Y;\Z_{\ell}\otimes_{\Z} G)\cong H^{\ast}(\hat{Y}; \hat{G}^{\ell})
\end{equation} and 
\begin{equation}\label{ZlkG} 
H^{\ast}(Y;\Zln\otimes_{\Z} G)\cong H^{\ast}(\hat{Y};\Zln\otimes_{\Z} G)
\end{equation} for every $\nu$. 
\end{prop}

\begin{proof}
1. Since every space in $\hSh$ is cofibrant and since the spaces $K(\pi,n)$ are fibrant if $\pi$ is an abelian pro-$\ell$-group, cf. \cite{ensprofin}, \S 1.4, Lemme 2, we have $$\Hom_{\Ho(\hSh)}(X, K(\pi,n)) \cong \Hom_{\hSh}(X, K(\pi,n))/ \sim$$ where $\sim$ denotes simplicial homotopy. Then the proof of the analogue result \cite{may}, Theorem 24.4, works here as well. \\  
2. We know that $K(\pi,n) \to \hat{\Omega}K(\pi,n+1)$ is a homotopy equivalence, hence it is a $\Z/\ell$-equivalence and $H\pi$ is an $\hat{\Omega}$-spectrum. From (\ref{omegacohom}) and $H^n(X;\pi) \cong \Hom_{\hHh}(X,K(\pi,n))$ the assertion follows.\\
3. Since $\pi$ is an abelian pro-$\ell$-group, $K(\pi,n)$ is a fibrant profinite space for each $n$. Hence by adjunction we have natural isomorphisms $\Hom_{\Hh}(Y, |K(\pi,n)|) \cong \Hom_{\hHh}(\hat{Y}, K(\pi,n))$. But by the construction of Eilenberg-MacLane spaces we have $|K(\pi,n)| \cong K(|\pi|,n)$ in $\Sh$ and we deduce the desired result from the first point above.  
\end{proof}

\begin{remark}
1. The third statement (\ref{ZlG}) of the proposition is a useful generalization of Remark \ref{remarkcohom}. The case (\ref{ZlkG}) is in fact a special case of (\ref{ZlG}), but also of Remark \ref{remarkcohom} with $\pi = \Zln\otimes_{\Z} G$, since this is a finite abelian group if $G$ is finitely generated abelian.\\
2. Since the model structures on $\hSh$ and $\hSp$ depend on the chosen prime number $\ell$, it is clear that one cannot  represent continuous cohomology with arbitrary coeffiecients in $\hSHh$. The homotopy groups of profinite spaces are pro-$\ell$-groups and hence only cohomology with pro-$\ell$-coefficients is representable in $\hSHh$.\\
3. It follows from the proposition that the mod $\ell$-Steenrod algebra of continuous cohomology operations is identical to the usual one. 
\end{remark}


\subsection{Postnikov decomposition}

We show the existence of a Postnikov decomposition for every profinite spectrum. Although some books use the existence of arbitrary colimits in the usual category of spectra for this purpose, there is another way to do so without referring to colimits. This enables us to use this construction in our profinite setting.\\
First we note that by general nonsense on pointed model categories we get the existence of fiber and cofiber sequences in $\hSp$. We can construct for every profinite spectrum $E$ a connective covering $E_{\geq 0} \to E$ which induces an isomorphism $\pi_k(E_{\geq n})=\pi_k(E)$ for all $k\geq 0$ and with $\pi_k(E_{\geq 0})=0$ for all $k<0$. We may construct it by considering the diagram

\begin{equation}
\begin{array}{ccc}
  &  & P(H\pi_0(E)[1])\\
  &  & \downarrow\\
E & \longrightarrow & H\pi_0(E)[1]
\end{array}
\end{equation}
where $H\pi$ denotes the Eilenberg-MacLane spectrum for a profinite abelian group $\pi$ and $P(H(\pi_0(E)[1])$ is the path object of $ H\pi_0(E)[1]$. We define $E_{\leq 0}$ to be the fiber product of this diagram.  Then we obtain inductively spectra $E_{\leq n}$ and maps $E \to E_{\leq n}$ for every $n$ and we call $E_{\leq n}$ the {\em Postnikov $n$-truncation of $E$}. It has the property  $\pi_k(E_{\leq n})=\pi_k(E)$ for $k\leq n$ and $\pi_k(E_{\leq n})=0$ for $k> n$. We get in fact a {\em Postnikov tower} $ \ldots \to E_{\leq n} \to E_{\leq n-1} \to \ldots$ for every $E$. The cofibre of each morphism $E_{\leq n} \to E_{\leq n-1}$ is identified with the Eilenberg-MacLane spectrum $H(\pi_n(E))[n+1]$. \\
Again by general nonsense \cite{homalg} I,\S 3, we know that the morphism $E \to E_{\leq n-1}$ fits into a fiber sequence $E_{\geq n} \to E \to E_{\leq n-1}$ and we call the spectrum $E_{\geq n}$ the {\em $n$-connective covering of $E$} which has the property $\pi_k(E_{\geq n})=\pi_k(E)$ for $k\geq n$ and $\pi_k(E_{\geq n})=0$ for $k< n$. 

\subsection{Atiyah-Hirzebruch spectral sequence}

Let $X$ be a profinite space. For every integer $p$ we denote by $\sk_pX$ the profinite subspace of $X$ which is generated by the simplices of degree less or equal $p$. We call $\sk_pX$ the $p$-skeleton of $X$. For every $k$ the set of $k$-simplices of $\sk_pX$ is closed in $X_k$ such that $\sk_pX$ is simplicial profinite subset of $X$, cf. \cite{dehon}. In addition, $\sk_pX$ is equal to the colimit of the diagram \begin{equation}\label{skeleton}
\begin{array}{ccc}
X_p \times \sk_{n-1}\Delta[p] & \longrightarrow & \sk_{p-1}X\\
\downarrow &   & \downarrow\\
X_p \times \Delta[p] & \longrightarrow & \sk_pX
\end{array}
\end{equation}
which is induced by the obvious map $X_p \times \Delta[p] \to X$, where we consider $X_p$ as a constant simplicial profinite set. 
This description implies that we get morphisms $\sk_{p-1}X \to \sk_pX$, which are cofibrations; and it implies that $X$ is equal to the colimit in $\hSh$ of the sequence $\ldots \sk_{p-1}X \to \sk_pX \to \ldots$  
One should note that if $X/-$ denotes the diagram of the simplicial finite quotients of $X$, then the limit of the diagram $\sk_p(X/-)$ is equal to $\sk_pX$. With this knowledge in hands we may continue to construct in the classical way the Atiyah-Hirzebruch spectral sequence for profinite cohomology theories. In the following we will write $X^p$ for $\sk_pX$.

Let $\pi$ be an abelian profinite group. We define the cochain complex $$D^p(X,\pi):=\tilde{H}^p(X^p/X^{p-1};\pi)$$ with differential $d:D^p(X;\pi) \to D^{p+1}(X;\pi)$ defined as the composite
$$\tilde{H}^p(X^p/X^{p-1};\pi) \to H^{p+1}(X^{p+1};\pi) \to \tilde{H}^{p+1}(X^{p+1}/X^p;\pi).$$ One verifies easily that this defines a cochain complex. For the construction of the profinite Atiyah-Hirzebruch spectral sequence the following result will be important in order to identify the $E_2^{p,q}$-term.

\begin{prop}\label{cellcomplex}
For every profinite space $X$, the cohomology of $X$ is given by the above cochain complex, i.e.
$$H(D^{\ast}(X;\pi),d) \cong H^{\ast}(X;\pi).$$
\end{prop}

\begin{proof}
The proof is essentially the same as for CW-complexes. We consider the exact couple arising from the long exact sequence for the pair $(X^p,X^{p-1})$:
\begin{equation}
\begin{array}{ccc}
H^{p+q}(X^{p-1};\pi) & \longleftarrow & H^{p+q}(X^p;\pi)\\
\partial \searrow &   & \nearrow  \\
 &  \tilde{H}^{p+q}(X^p/X^{p-1};\pi). &  
\end{array}
\end{equation}
Via Diagram (\ref{skeleton}) we see that $X^p/X^{p-1}$ is isomorphic to a limit of joints of $p$-spheres and that hence $\tilde{H}^{p+q}(X^p/X^{p-1};\pi)=0$ for $q>0$. Since $d_1:E_1^{p,0} \to E_1^{p+1,0}$ is equal to $d:D^p(X;\pi) \to D^{p+1}(X;\pi)$ we have
$$E_2^{\ast,0}=E_{\infty}^{p,0}=H(D^{\ast}(X;\pi),d).$$ 
From the definition of continuous cohomology for profinite spaces and from $(\sk_{p+r}X)_p=X_{p+r}$ for $r\geq 1$, it is clear that $$H^p(X^{p+r};\pi)=H^p(X;\pi)~\mathrm{for}~r \geq 1.$$ 
By induction on $p$, we show in addition that $$H^{p+r}(X^p;\pi)=0~\mathrm{for}~r \geq 1.$$ This is clear for $p=0$. Using (\ref{skeleton}) we get the exact sequence
$$0=\tilde{H}^{p+r}(X^p/X^{p-1};\pi) \to \tilde{H}^{p+r}(X^p;\pi) \to \tilde{H}^{p+r}(X^{p-1};\pi)
\to \tilde{H}^{p+r+1}(X^p/X^{p-1};\pi)=0.$$ By induction, $H^{p+r}(X^{p-1};\pi)=0$ for $r \geq 1$ and, by exactness, $H^{p+r}(X^p;\pi)=0$ for $r \geq 1$, too.
One continues just as in topology to conclude the desired isomorphism $E_2^{\ast,0}\cong H^{\ast}(X;\pi)$, see e.g. \cite{mccleary}, Theorem 4.11.
\end{proof}

Let $E$ be a profinite spectrum.  We want to construct an Atiyah-Hirzebruch spectral sequence for $E$.
We consider the skeletal filtration $\sk_0X \subset \ldots \subset \sk_pX \subset \sk_{p+1}X \subset \ldots \subset X$. This filtration yields a filtration on $E^{\ast}X$ defined by 
$F^pE^{\ast}X:=\mathrm{Ker}(E^{\ast}X \to E^{\ast}X^{p-1})$.

The proof of the following theorem is essentially the one given by Adams in \cite{adams}. We use the special consideration of the profinite setting by Dehon in \cite{dehon}, Proposition 2.1.9. But we are able to prove a slightly more general statement than Dehon for any profinite spectrum.\\
By Remark \ref{profinitestablepi} we know that the homotopy groups $\pi_qE$ have a natural profinite structure and hence also the coefficient groups $E^q$. Thus we may consider continuous cohomology with coefficients in $E^q$.

\begin{theorem}\label{ahss}
For any profinite spectrum $E$ and for any profinite space $X$ there is a spectral sequence $\{E_r^{p,q}\}$ with $E_2^{p,q}\cong H^p(X;E^q) \Longrightarrow E^{p+q}(X)$.
The spectral sequence converges strongly, in the sense of \cite{adams} III, \S 8.2, to the graded term $e_{\infty}^{p,q}=F^pE^{\ast}X/F^{p+1}E^{\ast}X$ of the filtration on $E^{\ast}X$ if ${\lim_r}^1E_r^{p,q}=0$. In particular, the spectral sequence converges if $H^p(X;E^q)=0$ for $p$ sufficiently large or if $H^p(X;E^q)$ is finite for all $p$.  
We call it the {\em profinite Atiyah-Hirzebruch spectral sequence}. 
\end{theorem} 

\begin{proof} 
The construction of the spectral sequence is a standard argument. The cofibre sequence $X^{p-1}_+ \to X^p_+ \to X^p_+/X^{p-1}_+$ induces an exact couple in cohomology
\begin{equation}\label{excouple}
\begin{array}{ccc}
\bigoplus\limits_{p} E^{\ast} (X^{p-1}) & \longleftarrow &
\bigoplus\limits_{p} E^{\ast} (X^p)\\
\partial \searrow &   & \nearrow  \\
 &  \bigoplus\limits_{p}  E^{\ast} (X^p/X^{p-1}) &  
\end{array}
\end{equation}
which leads to a spectral sequence with $E_1^{p,q}= E^{p+q}(X^p/X^{p-1}) \Longrightarrow E^{p+q}(X)$, the differentials of which are of bidegree $(r,1-r)$. The convergence properties are described in and follow from \cite{adams} III, \S 8.2 or \cite{etaleK}, Proposition A.2. In particular, if $H^p(X;E^q)=0$ for $p\geq n+1$, the spectral sequence degenerates at the $E_{n+1}$-term and $\oplus_pE_{n+1}^{p,\ast}$ is isomorphic to the quotient $F^pE^{\ast}X/F^{p+1}E^{\ast}X$. Hence it remains to identify the $E_2$-term.

{\em Claim}: We have a natural isomorphism 
$$E_1^{p,q}= E^{p+q}(X^p/X^{p-1}) \cong \tilde{H}^p(X^p/X^{p-1};\pi_pE_{p+q}).$$
If we can prove this claim, we get, together with Proposition \ref{cellcomplex} and the fact $\pi_pE_{p+q} =E^q$, that $E_2^{p,q}\cong H^p(X;E^q)$ and we are done. 

Hence it remains to prove the claim. In the following we suppose that $E$ is a fibrant profinite spectrum, i.e. an $\hOmega$-spectrum. This identification is due to \cite{dehon}. We follow his argumentation. By the characterizing diagram (\ref{skeleton}) of the $p$th-skeleton $X^p$ of $X$ we see that $X^p/X^{p-1} \cong \lim_{F \in \Qh(X_p)} (\vee_F S^p)$ is a limit in $\hShp$ of joints of $p$-spheres. The limit comes from the profinite structure of $X_p$; $F$ runs through the finite sets equal to quotients of $X_p$. Since $E$ is fibrant we can identify $E_{p+q}$ with the limit $\lim_t E_{p+q}(t)$ where $E_{p+q}(t)$ is a finite $\ell$-space for every $t$, i.e. $\pi_iE_{p+q}(t)$ is a finite $\ell$-group for every $i$ and $t$ and  $\pi_i E_{p+q}(t)\neq 0$ only for a finite number of $i$. 
Since $\vee_F S^p$ is a simplicial finite set and since $\pi_i E_{p+q}(t)$ is finite, we conclude that we have an isomorphism $$\lim_t \colim_{F \in \Qh(X_p)} \tilde{H}^p(\vee_F S^p;\pi_pE_{p+q}(t)) \cong \tilde{H}^p(X^p/X^{p-1};\pi_p E_{p+q}).$$ 
By \cite{switzer}, Proposition 6.39, we have 
$$\begin{array}{rcl}\tilde{H}^p(\vee_F S^p;\pi_p E_{p+q}(t)) & \cong & \Hom_{\hHhp}(\vee_F S^p, K(\pi_pE_{p+q}(t),p))\\
 & \cong & \Hom_{\Ab}(\pi_p(\vee_F S^p), \pi_p E_{p+q}(t)).\end{array}$$ 
Hence we get an isomorphism
$$\tilde{H}^p(X^p/X^{p-1};\pi_p E_{p+q}) \cong \Hom_{\pro-\Ab}(\pi_p(\vee_- S^p), \pi_p E_{p+q}(-)).$$ 
On the other hand we have the natural map
$$\begin{array}{rcl}E^{p+q}(X^p/X^{p-1}) & = & \Hom_{\hHhp}(X^p/X^{p-1},E_{p+q}) \\
 & \to & \Hom_{\pro-\Hhp}(X^p/X^{p-1},E_{p+q}(-)) \\
  & \to & \Hom_{\pro-\Ab}(\pi_p(\vee_- S^p), \pi_pE_{p+q}(-)).\end{array}$$
Hence in order to prove the claim, it suffices to show the following lemma.

\begin{lemma}\label{lemma2.1.10} {\rm Cf. Lemma 2.1.10 of \cite{dehon}.}
Let $S(-)$ be a filtered diagram of finite joints of $p$-spheres with limit $S$ in $\hShp$ and let $E$ be a fibrant profinite spectrum. Then the composite
$$\Hom_{\hHhp}(S,E_{p+q}) \to \Hom_{\pro-\Hhp}(S(-),E_{p+q}(-)) \to 
\Hom_{\pro-\Ab}(\pi_p(\vee_- S^p), \pi_pE_{p+q}(-))$$
is an isomorphism of abelian groups. 
\end{lemma}
\begin{proof}
The proof is the same as for Lemma 2.1.10 of \cite{dehon}. For the sake of completeness, we reformulate it here in more detail. For every pair of indices $F$ and $t$ we have a sequence of isomorphisms using the facts that $\pi_pS^p\cong \Z$, $S^p$ considered in $\Shp$, and that $E_{p+q}(t)$ is fibrant:
$$\begin{array}{rcl}
\Hom_{\Hhp}(\vee_F S^p, E_{p+q}(t)) & \cong & \prod_F  \Hom_{\Hhp}(S^p, E_{p+q}(t))\\
 & \cong & \prod_F \pi_p(E_{p+q}(t)) \cong \prod_F  \Hom_{\Ab}(\pi_p(S^p), \pi_p(E_{p+q}(t)))\\
 & \cong &  \Hom_{\Ab}(\pi_p(\vee_F S^p), \pi_p(E_{p+q}(t)))\end{array}$$
from which we deduce that the map 
$$\Hom_{\pro-\Hhp}(\vee_- S^p,E_{p+q}(-)) \to \Hom_{\pro-\Ab}(\pi_p(\vee_- S^p), \pi_p(E_{p+q}(-)))$$
is an isomorphism. Hence it remains to show that the left hand side is isomorphic to 
$\Hom_{\hHhp}(\lim_{F \in \Qh(X_p)} (\vee_F S^p), E_{p+q})$.
Since $E_{p+q}(t)$ is a finite $\ell$-space for every $t$ and since $\vee_F S^p$ is a simplicial finite set, Proposition 1.3.2 of \cite{dehon} implies that 
$$\colim_{F \in \Qh(X_p)} \Hom_{\Hhp}(\vee_F S^p, E_{p+q}(t)) \stackrel{\cong}{\longrightarrow}
\Hom_{\hHhp}(\lim_{F \in \Qh(X_p)} \vee_F S^p, E_{p+q}(t))$$ is an isomorphism for every $t$. Furthermore, Corollary 1.4.3 of \cite{dehon} implies that 
$$\Hom_{\hHhp}(\lim_{F \in \Qh(X_p)} \vee_F S^p, \lim_t E_{p+q}(t)) \cong \lim_t \Hom_{\hHhp}(\lim_{F\in \Qh(X_p)} \vee_F S^p, E_{p+q}(t))$$ is an isomorphism if ${\lim_t}^1 \pi_1(\homp(\lim_{F \in \Qh(X_p)} \vee_F S^p, E_{p+q}(t)))$ vanishes. If we could show that ${\lim_t}^1 \pi_1(\homp(\lim_{F \in \Qh(X_p)} \vee_F S^p, E_{p+q}(t)))=0$ we would get
$$\lim_t \colim_{F \in \Qh(X_p)} \Hom_{\Hhp}(\vee_F S^p, E_{p+q}(t)) \cong 
\Hom_{\hHhp}(\lim_{F \in \Qh(X_p)} \vee_F S^p, E_{p+q}(t)).$$ Hence we have to show that ${\lim_t}^1 \pi_1(\homp(\lim_{F \in \Qh(X_p)} \vee_F S^p, E_{p+q}(t)))$ vanishes. 
Again since $\vee_F S^p$ is a simplicial finite set and since $E_{p+q}(t)$ is fibrant we get for every $t$
$$\begin{array}{rcl}\pi_1(\homp(\lim_{F \in \Qh(X_p)} \vee_F S^p, E_{p+q}(t))) & \cong & 
\colim_{F \in \Qh(X_p)} \pi_1(\homp(\vee_F S^p, E_{p+q}(t))\\
 & \cong & \colim_{F \in \Qh(X_p)} \pi_1(\vee_F \Omega^p(E_{p+q}(t)))\\
 & \cong & \colim_{F \in \Qh(X_p)} \oplus_F \pi_{p+1}(E_{p+q}(t)).
\end{array}$$
Since $\pi_{p+1}E_{p+q}(t)$ is a finite group, the pro-abelian group 
$\{ \colim_{F\in \Qh(X_p)} \oplus_F \pi_{p+1}(E_{p+q}(t))\}_t$ is pro-isomorphic to a tower of surjections. Hence the ${\lim}^1$-term of the tower $\{\pi_1(\homp(\lim_{F\in \Qh(X_p)} \vee_F S^p, E_{p+q}(t))\}_t$ vanishes.  
\end{proof}
\end{proof}

\subsection{A K\"unneth formula for profinite cobordism}

We finish this section with a discussion of a K\"unneth isomorphism in $\hMU\Zln$ for the product $X\times Y$ in $\hSh$ for profinite spaces $X$ and $Y$ satisfying certain conditions. The whole argumentation is inspired by the work of Dehon \cite{dehon}, where he proves the corresponding result for $\hMU$. We say that $Y \in \hSh$ is {\em without $\ell$-torsion} if the reduction map $\Zl^n \to \Zl$ induces a surjection $H^{\ast}(Y;\Zl^n)\to H^{\ast}(Y;\Zl)$ for every positive $n$, cf. \cite{dehon}, Definition 2.1.5. This is equivalent to the fact that $H^{\ast}(Y;\Z_{\ell})$ has no torsion. 
For the application to \'etale cohomology theories in later sections, we only need the following version of a K\"unneth isomorphism for finite dimensional profinite spaces.
 
\begin{theorem}\label{KunnethforhMU}
Let $X$ and $Y$ be finite dimensional profinite spaces and let $X$ be without $\ell$-torsion. Then we have for every $\nu$ K\"unneth isomorphisms
$$\hMU^{\ast}(X) \otimes_{\hMU^{\ast}} \hMU^{\ast}(Y)
\stackrel{\cong}{\longrightarrow}\hMU^{\ast}(X \times Y)$$ and
$$\hMU^{\ast}(X;\Zln) \otimes_{(\hMU\Zln)^{\ast}} \hMU^{\ast}(Y;\Zln)
\stackrel{\cong}{\longrightarrow}\hMU^{\ast}(X \times Y;\Zln).$$
\end{theorem}
\begin{proof}
The first assertion is a restricted version of \cite{dehon}, Proposition 4.3.2.\\
The part for $\Zln$-coefficients follows in a similar way. I am very grateful to Francois-Xavier Dehon for an explanation of this point. We will prove the assertion following his argument.\\
We suppose first that $X$ and $Y$ are finite dimensional without $\ell$-torsion. The point is that the morphism of coefficients $(\hMU\Zln)^{\ast} \to (H\Z \wedge MU\Zln)^{\wedge \ast}$ is injective and, since $X$ has no $\ell$-torsion, the map $E_2^{\ast,\ast}(X;MU\Zln) \to E_2^{\ast,\ast}(X;H\Z \wedge MU\Zln)$ is also injective. Since the spectral sequence for $H\Z \wedge MU\Zln$ degenerates at the $E_2$-stage, the same is true for the spectral sequence for $MU\Zln$. One can deduce from this fact that we have an isomorphism $\hMU^{\ast}(X;\Zln)/f^1 \cong H^{\ast}(X;\Zl)$ as described in the proof of \cite{dehon}, Proposition 2.1.8, where $f^{\ast}$ is an $\ell$-adic filtration on $MU^{\ast}$-modules defined in \cite{dehon}, \S 2.1.\\
Since $Y$ and $X \times Y$ have no $\ell$-torsion, the same argument is valid for them as well. By Corollaire 2.1.4 of \cite{dehon}, it suffices to check that the map
$$\hMU^{\ast}(X;\Zln)/f^1 \otimes_{(\hMU\Zln)^{\ast}}\hMU^{\ast}(Y;\Zln)/f^1 \to \hMU^{\ast}(X \times Y;\Zln)/f^1$$
is an isomorphism. But this follows from the previous discussion and the K\"unneth isomorphism for $\Zl$-cohomology. This finishes the case that both $X$ and $Y$ have no $\ell$-torsion.\\ 
For the case that $Y$ has $\ell$-torsion, we use induction on the skeletal filtration of $Y$. Let $Y^n$ be the $n$-th skeleton of $Y$. If $\hMU^{\ast}(X;\Zln)\otimes_{(\hMU\Zln)^{\ast}}\hMU^{\ast}(Z;\Zln) \to \hMU^{\ast}(X \times Z;\Zln)$ is an isomorphism for $Z=Y^n$ and for $Z=Y^{n+1}/Y^n$, then it is also an isomorphism for $Z=Y^{n+1}$. This follows from the long exact sequence associated to the cofibre sequence $Y^n \to Y^{n+1} \to Y^n/Y^{n+1}$ and the exactness of the tensor product with  $\hMU^{\ast}(X;\Zln)$. The exactness of $\hMU^{\ast}(X;\Zln) \otimes -$ is due to the fact that $X$ has no $\ell$-torsion. We deduce the assertion of the theorem by induction on $n$ from the proof of Proposition 4.3.2 of \cite{dehon} using the point that we consider only the $n$-th skeleton of $Y$ instead  of the whole space. 
\end{proof}

\begin{remark}
There is a version of the previous theorem with the assumption that $Y$ is an arbitrary profinite space and $X$ is a profinite space without $\ell$-torsion whose $\Zl$-cohomology is finite in each degree.\\
One proves the $\Zln$-coefficient case by taking the limit over the skeletal filtration. 
\end{remark}

\subsection{Comparison with generalized cohomology theories of pro-spectra}

Isaksen constructs in \cite{gencohomofpro} a stable model structure on the category of pro-spectra and defines generalized cohomology theories of pro-spectra. If $E$ is a pro-spectrum then the $r$-th cohomology $E_{\pro}^r(X)$ of a pro-spectrum $X$ with coefficients in $E$ is the set $[X,E]^{-r}_{\pro}$ maps that lower the degree by $r$ in the stable homotopy category of pro-spectra. In addition, for two pro-spectra $X$ and $E$ he shows the existence of an Atiyah-Hirzebruch spectral sequence $E_2^{p,q}=H^{-p}(X;\pi_{-q}E)\Rightarrow [X,E]^{p+q}_{\pro}$ which is in particular convergent if $E$ is a constant pro-spectrum and $X$ is the suspension spectrum of a finite dimensional pro-space. We show that this definition of cohomology of pro-spectra gives in the cases of interest the same groups as our definition via profinite spectra.\\
First we construct a functor $\compl :  \pro-\Sp \to \hSp$ that induces morphisms on cohomology theories. As for pro-spaces we define $\compl$ to be the composition of the profinite completion on spectra followed by taking the homotopy limit in $\hSp$ of the corresponding diagram in $\hSp$. Note that the homotopy limit exists by \cite{hirsch}, Chapter 18, keeping in mind that all limits exist, see  Convention \ref{convention} on limits. Since the completion commutes with suspension this functor obviously agrees with our previous definition for the case of the suspension spectrum of a space. A map $f:E \to F$ in $\pro-\Sp$ is a weak equivalence if $f$ is a levelwise $n$-equivalence for all $n$. This implies that $f$ induces an isomorphism $\pi_k f$ on pro-homotopy groups, cf. \cite{gencohomofpro}, Theorem 8.4. We want to show that $\compl$ sends a weak equivalence $f$ to a stable equivalence in $\hSp$. For this we may suppose that $f$ is a level map $\{f_s:E_s \to F_s\}$ and that $\pi_k f_s$ is an isomorphism of stable homotopy groups for every $k$ and every $s$. We know that the completion $\compl : \Sp \to \hSp$ sends stable equivalences to stable equivalences. Since the homotopy limit is by construction well behaved with respect to levelwise weak equivalences $\holim_s \hat{f_s}$ is also a stable equivalence. Hence the functor induces maps $\Hom_{\Ho(\pro-\Sp)}(X,E) \to \Hom_{\hSHh}(\hat{X},\hat{E})$ and hence also maps $E_{\pro}^{\ast}(X) \to \hat{E}^{\ast}(\hat{X})$. These maps yield morphisms of Atiyah-Hirzebruch spectral sequences. 

\begin{theorem}\label{proMUcomparison}
For any $k$, consider the constant pro-spectrum $MU\Zln$. Let $\{X_s\}$ be a finite dimensional pro-space. Then we have an isomorphism 
$$MU^{\ast}_{\pro}(\{X_s\};\Zln)\cong \hMU^{\ast}(\hat{X};\Zln).$$
\end{theorem} 
\begin{proof}
As pointed out above the completion functor induces a map of spectral sequences. But for a constant pro-spectrum the completion is just the completion defined in Section 3 on profinite spectra.\\ 
Hence our previous results on the coefficients of $MU\Zln$ and $\hMU\Zln$ imply that it induces an isomorphism of coefficients $(MU\Zln)^{\ast}\cong (\hMU\Zln)^{\ast}$ and hence an isomorphism of $E_2$-terms. Since both spectral sequences converge for the given $X$, the abutments are also isomorphic. 
\end{proof}

\section{Profinite \'etale realization on the unstable $\A^1$-homotopy category}

The construction of the \'etale topological type functor $\Et$ from locally noetherian schemes to pro-simplicial sets is due to Artin-Mazur and Friedlander. The construction of the $\A^1$-homotopy category of schemes makes naturally arise the question if one could enlarge this functor to the category of spaces and if this functor behaves well with respect to the model structure. These questions have been answered by Isaksen and Schmidt. The first step into this direction was the construction of a model structure on $\pro-\Sh$. In view of $\Et$, there are in fact at least two interesting structures. It turns out that the $\Zl$-cohomological model structure of $\pro-\Sh$ fits better with the $\A^1$-localized homotopy category. This is due to the fact that over a field of positive characteristic $p>0$, the \'etale fundamental group of $\A^1$ is highly non-trivial. To avoid this obstruction one has to complete away from the characteristic of the ground field using \'etale cohomology. The projection $X\times \A^1 \to X$ induces an isomorphism in \'etale cohomology $H_{\et}^{\ast}(X;\Zl) \to H_{\et}^{\ast}(X\times \A^1;\Zl)$ for every prime $\ell \neq p$.\\
We will adapt the constructions of Isaksen \cite{a1real} and Dugger \cite{dugger} to $\hSh$ and check that we still get an induced left derived functor on the homotopy categories. We will calculate $\hEt$ for some examples.      

\subsection{The functor $\hEt$}

\begin{defn} {\rm Definition 4.4 of \cite{fried}.}\\
Let $X_{\cdot}$ be a locally noetherian simplicial scheme. The {\rm \'etale topological type} of $X_{\cdot}$ is defined to be the following pro-simplicial set
$$\Et X = \pi \circ \Delta :HRR(X_{\cdot}) \to \Sh$$
sending a hypercovering $U_{\cdot,\cdot}$ of $X_{\cdot}$ to the simplicial set of connected components of the diagonal of $U_{\cdot,\cdot}$. If $f: X_{\cdot} \to Y_{\cdot}$ is a map of locally noetherian simplicial schemes, then the {\rm strict \'etale topological type} of $f$ is the strict map $\Et f:\Et X_{\cdot} \to \Et Y_{\cdot}$ given by the functor $f^{\ast}:HRR(Y_{\cdot}) \to HRR(X_{\cdot})$ and the natural transformation $\Et X_{\cdot} \circ \Et f \to \Et Y_{\cdot}$.
\end{defn}
We refer the reader to \cite{fried} and \cite{a1real} for a detailed discussion of the category $HRR(X_{\cdot})$ of rigid hypercoverings of $X_{\cdot}$ and rigid pullbacks.\\ 
Isaksen now uses the insight of Dugger \cite{dugger} that one can construct the $\A^1$-homotopy category in a universal way. Starting from an almost arbitrary category $\Ch$, Dugger constructs an enlargement of $\Ch$ that carries a model structure and is universal for this property. He also shows how to enlarge functors $\Ch \to \Mh$ from $\Ch$ to a model category $\Mh$. In \cite{a1real}, Isaksen takes Dugger's model and extends the functor $\Et$ via this general method to the $\A^1$-homotopy category. The idea is that $\Et X$ should be the above $\Et X$ on a representable presheaf $X$ and should preserve colimits and the simplicial structure.\\
For our purpose, we would like to define $\hEt$ directly in the way Dugger suggests. But the problem is that $\compl:\pro-\Sh \to \hSh$ is not a left adjoint functor and does not preserve all small colimits.
Therefore, we define $\hEt$ to be the composition of $\Et$ followed by completion. To be precise we make the following definition:

\begin{defn}\label{defet}
If $X$ is a representable presheaf, then $\Et X$ is the \'etale topological type of $X$.
If $P$ is a discrete presheaf, i.e. each simplicial set $P(U)$ is $0$-dimensional, i.e. $P$ is
just a presheaf of sets, then $P$ can be written as a colimit $\colim_i X_i$ of representables and we define $\Et P:=\colim_i \Et X_i$. Finally, an arbitrary simplicial presheaf can be written as the coequalizer
of the diagram
$$\coprod_{[m]\to[n]} P_m \otimes \Delta^n \rightrightarrows \coprod_{[n]} P_n \otimes \Delta^n,$$
where each $P_n$ is discrete. Define $\Et P$ to be the coequalizer of the diagram
$$\coprod_{[m]\to[n]} \Et P_m \otimes \Delta^n \rightrightarrows \coprod_{[n]} \Et P_n \otimes \Delta^n.$$
We define the {\rm profinite \'etale topological type functor} $\hEt$ to be the composition of $\Et$ and the profinite completion functor $\pro-\Sh \to \hSh$:
\begin{equation}\label{hEt}
\hEt := \widehat{(\cdot)} \circ \Et : \Delta^{\mathrm{op}}\mathrm{PreShv}(\mathrm{Sm}/S) \to \hSh.
\end{equation}
\end{defn}

For computations the following remark is crucial. 
\begin{remark}\label{remarketalecohom}
1. By \cite{fried}, Proposition 5.9, we know that $H^{\ast}_{\et}(X;M) \cong H^{\ast}(\Et X,M)$ where $H^{\ast}(Z;M)$ denotes the cohomology of a pro-simplicial set $Z$ with coefficients in the local coefficient system $M$ corresponding to the sheaf $M$. For a finite abelian group $\pi$, we have in addition a natural isomorphism by Lemma \ref{Z/ell-comparison}: $H^{\ast}(Z;\pi) \cong H^{\ast}(\hat{Z};\pi)$ for every pro-simplicial set $Z$. Hence we get as well $$H^{\ast}_{\et}(X;\pi) \cong H^{\ast}(\hEt X,\pi)$$ for every locally noetherian scheme $X$ and every finite abelian group $\pi$.\\
2. For a morphism $g:Z \to X$ of schemes over $k$, there are relative cohomology groups $H^{\ast}(\hEt X,\hEt Z;\Zl)$ fitting in a natural long exact sequence 
$$\ldots \to H^{\ast}(\hEt X,\hEt Z;\Zl) \to H^{\ast}(\hEt X;\Zl) \to H^{\ast}(\hEt Z;g^{\ast}\Zl) \to \ldots$$
For a closed immersion $Z \hookrightarrow X$ with open complement $U$, in \cite{fried}, Corollary 14.5, Friedlander has shown that these relative cohomology groups coincide with the \'etale cohomology of $X$ with support in $Z$, i.e. we have a natural isomorphism $$H_{\et,Z}^{\ast}(X;\Zl)\cong H^{\ast}(\hEt X,\hEt U;\Zl).$$
\end{remark}

Although $\hEt$ does not commute with products we have a weaker compatibility.
\begin{prop}\label{hEtofproduct}
Let $X$ and $Y$ be smooth schemes of finite type over a separably closed field of characteristic $p \neq \ell$. Then the canonical map is a weak equivalence in $\hSh$
$$\hEt(X\times Y) \stackrel{\simeq}{\longrightarrow}\hEt X \times \hEt Y.$$ 
\end{prop}
\begin{proof}
This follows from the K\"unneth formula for smooth schemes over a separably closed field proved in \cite{sga412} [Th. finitude], Corollaire 1.11. It provides an isomorphism 
$$H_{\et}^{\ast}(X\times Y;\Zl) \cong H_{\et}^{\ast}(X;\Zl) \otimes H_{\et}^{\ast}(Y;\Zl).$$ 
On the other hand, we have a K\"unneth formula for the cohomology of profinite spaces which yields an isomorphism $$H^{\ast}(\hEt X\times \hEt Y;\Zl) \cong H^{\ast}(\hEt X;\Zl) \otimes H^{\ast}(\hEt Y;\Zl).$$ This implies that $\hEt(X\times Y) \longrightarrow \hEt X \times \hEt Y$is a weak equivalence in $\hSh$.
\end{proof}

In the following example, we consider objects in the category
$\Delta^{\mathrm{op}}\mathrm{PreShv}(\mathrm{Sm}/k)$ of simplicial presheaves over $\Sm/k$.  

\begin{example}\label{SpecR}
Let $R$ be a strict local henselian ring, i.e. a local henselian ring with separably closed residue field.
Then $\Spec R$ has no nontrivial \'etale covers and the \'etale topological type of $\Spec R$ is a contractible space, i.e. $\hEt (\Spec R)\cong \ast$.
\end{example}

\begin{example}\label{hEta}
Let $k$ be a separably closed field with $\mathrm{char}(k)\neq \ell$. Let $c:\Sh \to \Delta^{\mathrm{op}}\mathrm{PreShv}(\Sm/k)$ be the functor that sends a simplicial set $Z$ to the constant presheaf defined by $Z$. For every $n$, $cZ([n])$ is isomorphic to a disjoint union of copies of $\Spec k$. Hence, since $\Et$ commutes with coproducts and since $\Et (\Spec k) = \ast$ by Example \ref{SpecR} above, its \'etale topological type $\Et cZ$ is just the simplicial set $Z$ itself.\\ 
This shows that the pro-simplicial set $\Et Z$ is in fact just a simplicial set. The completion $\widehat{(\cdot)}:\Sh \to \hSh$ preserves all colimits. Hence the same argument holds for $\hEt$, $\hEt(cZ) = \widehat{\Et cZ}$ and $\hEt (cZ) = Z$ for every $Z\in \hSh$.\\
In particular, if $S^1$ denotes the simplicial circle, then $S^1([n])=\coprod_{S^1([n])}\Spec k$ and
$(S^1 \wedge \mathcal{X})([n]) = \coprod_{S^1([n])} \mathcal{X}([n])$.
This yields an isomorphism in $\hShp$
\begin{equation}\label{s1smash}
S^1 \wedge \hEt (\Xh) \stackrel{\cong}{\longrightarrow} \hEt (S^1) \wedge \hEt (\Xh)
\stackrel{\cong}{\longrightarrow} \hEt (S^1 \wedge \Xh).
\end{equation}
\end{example}

Furthermore, Friedlander proved in \cite{fried} that the \'etale fundamental group $\pi_1^{\et}(X)$ of a scheme is isomorphic as pro-group to the pro-fundamental group of $\Et X$. From Proposition \ref{piofZl-completion} we deduce the following
\begin{prop}\label{pi1et}
Let $X$ be a connected locally noetherian scheme. The fundamental group of $\hEt X$ is isomorphic to the $\ell$-completion of the \'etale fundamental group of $X$, i.e.
$$\pi_1^{\ell}(\hEt X) \cong \widehat{\pi_1^{\et}(X)}^{\ell}.$$     
\end{prop}

\begin{example}\label{pia1}
Let $k$ be a separably closed field with $\mathrm{char}(k)\neq \ell$. 
Since $H^{\ast}_{\et}(\A^1_k;\Z/\ell)\cong H^{\ast}_{\et}(k;\Z/\ell)$ we know that 
$\hEt \A^1_k \simeq \ast$, i.e. that $\hEt \A^1_k$ is contractible in $\hSh$ and
$$\pi^{\ell}_1(\hEt \A^1_k) = \{1\}.$$
\end{example}

\begin{example}\label{gm}
Let $k$ be a separably closed field with $\mathrm{char}(k)\neq \ell$. The space $\Ge_m$ is connected and its \'etale fundamental group is $\hat{\Z}$. Its $\ell$-completion is hence equal to $\Z_{\ell}$ and we get $\hEt \Ge_m \cong K(\Z_{\ell},1)$. Hence $\hEt \Ge_m$ is isomorphic to $S^1$ in $\hShp$ by (\ref{pi1s1}).
\end{example}

\begin{example}\label{p1}
Let $k$ be a separably closed field with $\mathrm{char}(k)\neq \ell$. Let $\Pro_k^1$ be the projective line over $k$. Since $\Pro_k^1$ is connected and $\pi^{\ell}_{1, \et} (\Pro_k^1)=0$ is trivial, cf.  e. g. \cite{milne} I, Example 5.2 f), $\hEt \Pro_k^1$ is a simply connected space. Apart from $H^0$ its only nonzero \'etale cohomology group is $H^2_{\et}(\Pro^1_k ,\Z/\ell) \cong \Zl$, see \cite{milne} VI, Example 5.6 with a chosen isomorphism $\Zl \to \Zl(1)$. Hence $\hEt \Pro_k^1$ is isomorphic in $\hShp$ to the simplicial finite set $S^2$.
\end{example}

\begin{example}\label{p1overk}
Let $k$ be a field of characteristic $p \neq \ell$. The \'etale realization of the projective line $\Pro^1_k$ over $k$ is given by the isomorphism $\hEt \Pro^1_k \stackrel{\cong}{\to} S^2 \times \hEt k$ in $\hSh$. For, the projective bundle formula in \'etale cohomology 
$H_{\et}^{\ast}(\Pro^1_k; \Zl) \cong H_{\et}^{\ast}(k;\Zl) \oplus H_{\et}^{\ast}(k;\Zl)$ and the K\"unneth formula in $\hSh$ yield the identification.
\end{example}

\begin{example}\label{Fq}
Let $k=\F_q$ be a finite field with $\mathrm{char}(k)=p\neq \ell$. The \'etale topological type of $k$ is isomorphic to $S^1$ in $\hSh$. For, $\hEt k$ is a connected space by \cite{fried} Proposition 5.2. Its $\ell$-completed fundamental group is the $\ell$-completion of the absolute Galois group of $k$ by Proposition \ref{pi1et}, i.e. $\pi^{\ell}_{1,\et}(k)=\Z_{\ell}$. Since all cohomology groups $H^i(k;\Zl)$ vanish for $i>1$ by \cite{serre}, $\Et k$ is a space of dimension one and weakly equivalent to $S^1$, hence it is isomorphic to $S^1$ in $\hSh$. 
\end{example}

\subsection{Unstable profinite \'etale realization}

We want to show that $\hEt$ induces a functor on the homotopy category. As indicated above we choose the universal model of \cite{dugger} for the unstable $\A^1$-homotopy category of smooth schemes over a base field $k$.\\
Dugger showed that after some localization the projective model structure on $U(\Sm/k):=\Delta^{\mathrm{op}}\mathrm{PreShv}(\Sm/k)$ is a model for the $\A^1$-homotopy category $\Hh _{\A^1}(k)$ of \cite{mv}. In this projective model structure the weak equivalences (fibrations) are objectwise weak equivalences (fibrations) of simplicial sets. The cofibrations are the maps having the left lifting property with respect to all trivial fibrations.
Then one takes the left Bousfield localization of this model structure at the set $S$ of maps:
\begin{enumerate}
\item\label{rel1} for every finite collection $\{X^a\}$ of schemes with disjoint union $X$, the
map $\coprod X^a \to X$ from the coproduct of the presheaves represented by each $X^a$
to the presheaf represented by $X$;
\item\label{rel2} every \'etale (Nisnevich) hypercover $U \to X$;
\item\label{rel3} $X\times \A^1 \to X$ for every scheme $X$.
\end{enumerate}
We call this the {\em \'etale (Nisnevich) $\A^1$-local projective model structure}
according to \cite{a1real}, and denote it by $\LeU=L_{S, \et} U(\Sm /k)$ (resp. $\LNU=L_{S, \mathrm{Nis}} U(\Sm /k)$).\\
Proposition 8.1 of \cite{dugger} states that $\LNU$ is Quillen equivalent to the Nisnevich $\A^1$-localized model category $\MVh$ of \cite{mv} and the analogue holds for the \'etale case.\\

The following lemma is clear from the theory of Bousfield localization.

\begin{lemma}\label{lflocalized}
Let $F:\Ch \to \Dh$ be a functor between model categories. Let $S$ be a set of maps in $\Ch$ and let $QC$ denote a fibrant replacement 
for objects $C$ of $\Ch$. Suppose that the total left derived functor $LF$ of $F$ and that the left Bousfield localizations $\Ch/S$ and $\Dh/FQ(S)$ exist.
Then $F$ induces a functor 
$$F/S:\Ch/S \to \Dh/FQ(S)$$  
and if $F$ sends the maps in $S$ into weak equivalences in $\Dh$, $F$ induces a total left derived functor on the localized category 
$\mathrm{Ho}(\Ch/S) \to \mathrm{Ho}(\Dh)$.
\end{lemma}

Since $\Et$ is a left Quillen functor on $U(\mathrm{Sm} /k)$ by \cite{dugger}, Proposition 2.3, in order to show that $\Et$ induces a left Quillen functor on $L_T U(\mathrm{Sm} /k)$, it suffices to show that $\Et$ takes the relations defined above into weak equivalences in $\pro-\Sh$. 

\begin{theorem}\label{LhEt}
Let $\ell$ be a prime different from the characteristic of $k$. With respect to the \'etale (Nisnevich) $\A^1$-local projective model structure on simplicial presheaves on $\Sm/k$, the functor $\hEt$ induces a functor $\LhEt$ from the \'etale (Nisnevich) $\A^1$-homotopy category to the $\Zl$-cohomological homotopy category of $\hSh$.
In particular, $\LhEt X$ is just $\hEt X$ for every scheme in $\Sm/k$, and hence $\hEt$ preserves $\A^1$-weak equivalences between smooth schemes over $k$.
\end{theorem}
\begin{proof}
Since every Nisnevich hypercover is also an \'etale hypercover it suffices to check this for the \'etale case. The first part of the proof is the one of \cite{a1real}, Theorem 2.6.
The functor $\Et$ above is exactly the functor
$\mathrm{Re}:U(\Sm/k)=\Delta^{\mathrm{op}}\mathrm{Pre}(\Sm/k) \to \pro -\Sh$
of \cite{dugger}, Proposition 2.3, hence it is a left Quillen functor. In order to see that $\Et$ is in fact a left Quillen functor on $L_S(U(\Sm/k))$, by Theorem 3.1.6 in \cite{hirsch} it remains to show
that it takes cofibrant approximations of maps in $S$ to weak equivalences of $\pro -\Sh$. This is done in \cite{a1real}.\\
Note that the condition on $\ell$ is needed to ensure that the projection $X \times \A^1 \to X$ induces an isomorphism on \'etale cohomology $H_{\et}^{\ast}(X;\Z/\ell) \stackrel{\cong}{\to} H_{\et}^{\ast}(X\times \A^1;\Z/\ell)$. The functor $\Et$ takes weak equivalences between cofibrant objects to weak equivalences.\\
By Lemma \ref{Z/ell-comparison} the completion functor $\widehat{(\cdot)}$ preserves weak equivalences and cofibrations. Hence the composition $\hEt$ sends weak equivalences between
cofibrant objects into weak equivalences and the total left derived functor
$$\LhEt: \Hh _{\A^1}^{\et}(k) \to \hHh$$ exists.\\
The last statement of the theorem follows from the definition of the total left derived functor
and the fact that all representable presheaves are cofibrant
in $L_S(U(\Sm/k))$, cf. \cite{dugger}.
\end{proof}

\begin{remark}\label{factorization}
Since every Nisnevich hypercover is also an \'etale hypercover, and since $\Et$ and hence also
$\hEt$ send all \'etale hypercovers to weak equivalences, the functor from the
Nisnevich $\A^1$-homotopy category factors through the \'etale $\A^1$-homotopy category:
\begin{equation}
\LhEt : \Hh_{\A^1}^{\mathrm{Nis}}(k) \to \Hh_{\A^1}^{\et}(k) \stackrel{\LhEt}{\to} \hHh
\end{equation}
where the first functor corresponds to sending a presheaf to its associated sheaf in the
\'etale topology.
\end{remark}

We conclude this section with a result of \cite{a1real} on distinguished square. We will deduce from this result that we get a Mayer-Vietoris sequence for profinite \'etale cohomology theories.\\
We recall the definition of an {\em elementary distinguished square} of \cite{mv}. It is a diagram
\begin{equation}\label{dsquare}
\begin{array}{ccc}
U\times_X V & \longrightarrow & V \\
\downarrow &   & \downarrow p\\
U  & \stackrel{i}{\hookrightarrow} &  X
\end{array}
\end{equation}
of smooth schemes in which $i$ is an open inclusion, $p$ is \'etale and  $p:p^{-1}(X-U) \to X-U$ is an isomorphism, where $X-U$ and $p^{-1}(X-U)$ are given the reduced structure. In particular, the maps $i$ and $p$ form a Nisnevich cover of $X$.    As in \cite{a1real}, Theorem 2.10, we can prove the following excision theorem for elementary distinguished squares. A similar result has already been proved by Friedlander in \cite{fried}, Lemma 14.10.

\begin{theorem}\label{excision}
Given an elementary distinguished square of smooth schemes over $k$. Then the square
\begin{equation}
\begin{array}{ccc}
\hEt(U\times_X V) & \longrightarrow & \hEt V \\
\downarrow &   & \downarrow \\
\hEt U  & \longrightarrow & \hEt X
\end{array}
\end{equation}
is a homotopy pushout square of profinite spaces.
\end{theorem}

\begin{proof}
The proof is the exact analog of the proof of Theorem 2.10 of \cite{a1real} with $\mathrm{LEt}$ replaced by $\LhEt$.
\end{proof}

\section{Profinite \'etale realization on the stable $\A^1$-homotopy category}

We extend the results of the previous section to the stable $\A^1$-homotopy category. One of the reasons why we consider the model $\hSh$ of $\Ho(\hSh)$, instead of $\pro-\Sh$, is that it seems to be easier to extend the functor $\hEt$ to the category of profinite spectra rather than pro-spectra.\\ 
We consider $S^1$-spectra and motivic $\Pro^1$-spectra separately. The first point for both will be to choose the correct model for the stable $\A^1$-homotopy theory. Then we deduce the existence of an \'etale realization as a left derived functor on the stable homotopy category from \cite{hovey}.\\
The problem for the motivic version is that $\Et$ and hence also $\hEt$ do not commute with products in general. Fortunately, $\hEt$ commutes with the smash product by the simplicial circle $S^1$. But for $\Pro^1$ we have to construct an intermediate category and then we will show by a zig-zag of functors that we get a functor on the homotopy level.\\
This \'etale realization of the stable motivic category is the technical key point for the construction of a transformation from algebraic cobordism given by the $MGL$-spectrum to the profinite \'etale cobordism of the next section. 

\subsection{Etale realization of motivic spectra}

Let $k$ be the base field and let $\ok$ be its separable closure. For the category of motivic $\Pro^1_k$-spectra over $k$, we have to consider presheaves $\Xh$ pointed by a morphism $\Spec k \to \Xh$. In particular, $\Spec k$ is the initial and terminal object.  If we want to construct a stable \'etale realization, this forces us to consider the category $\hShpk$ of pointed profinite spaces relative over $\hEt k$. Its objects $(X,p,s)$ are pointed profinite spaces $X$ together with a projection morphism  $p:X \to \hEt k$ and a section morphism $s:\hEt k \to X$. The morphisms in this category are commutative diagrams in the obvious sense. The \'etale realization of a pointed presheaf $\Xh$ is naturally an object of $\hShpk$ via the images of the projection and section morphisms of $\Xh$. They are pointed by the composite $\ast = \hEt \ok \to \hEt k \to \Xh$.\\
Via the canonical maps $X \to \ast = \hEt \ok \to \hEt k$ and $\hEt k \to \ast \to X$ we may view every space $X$ in $\hShp$ as an object in $\hShpk$. In particular, the pointed space $S^2$ is naturally an object in $\hShpk$. \\
We consider the usual model structure on $\hShpk$ where weak equivalences (resp. cofibrations, fibrations) are those maps which are $\Zl$-weak equivalences (resp. cofibrations, fibrations) in $\hShp$ after forgetting the projection and section maps. This model structure on $\hShpk$ has the same properties as the one on $\hShp$. In particular, it is left proper and fibrantly generated and we construct a stable model structure on the category $\mathrm{Sp}(\hShpk, S^2 \wedge \cdot)$ of profinite $S^2$-spectra over $\hEt k$ exactly in the same way. Its homotopy category will be denoted by $\hSHhtk$.\\  

\begin{remark}
The smash product $X \wedge_{\hEt k} Z$ in $\hShpk$ of a space $X \in \hShp$ considered in $\hShpk$ and an object $Z \in \hShpk$ is canonically isomorphic to the smash product $X \wedge Z$ in $\hShp$. For, if maps $T \to X$ and $T \to Z$ commute over $\hEt k$, they commute over $\ast$, since the projection map of $X$ factors through $\ast =\hEt \ok \to \hEt k$. Hence both products are canonically isomorphic  via their universal property.
\end{remark}
 
When we want to extend $\hEt$ to $\Pro^1$-spectra, we have to take into account the problem that the \'etale topological type functor does not commute with products in general. However the projections to each factor induce a canonical map 
\begin{equation}\label{canonicalmap} 
\hEt(\Pro^1_k \wedge_k X) \to \hEt(\Pro^1_k) \wedge_{\hEt k} \hEt(X).\end{equation}
\begin{lemma}\label{Kunneth}
For every pointed presheaf $\Xh$ on $\Sm/k$ the sequence of canonical maps in $\hShpk$
\begin{equation}\label{}
S^2 \wedge \hEt \Xh \stackrel{\simeq}{\longrightarrow} \hEt (\Pro^1_k) \wedge_{\hEt k} \hEt(\Xh)
\stackrel{\simeq}{\longleftarrow} \hEt(\Pro^1_k \wedge_k \Xh)
\end{equation}
is a sequence of $\Zl$-weak equivalences.
\end{lemma}
\begin{proof}
The projective bundle formula for \'etale cohomology implies that the \'etale cohomology $H^{\ast}(\hEt \Pro^1_k;\Zl)$ is a free module of rank $2$ over $H^{\ast}(k;\Zl)$. We deduce from this fact on the one hand that we have a relative K\"unneth isomorphism for $\hEt \Pro^1_k$ in $\hShpk$, see e.g. \cite{smith} for a discussion of relative K\"unneth theorems: 
$$H^{\ast}(\hEt \Pro^1_k \wedge_{\hEt k} \hEt X;\Zl) \cong H^{\ast}(\hEt \Pro^1_k;\Zl) \otimes_{H^{\ast}(\hEt k;\Zl)} H^{\ast}(\hEt X;\Zl).$$ This implies that the canonical map (\ref{canonicalmap}) is a $\Zl$-weak equivalence in $\hShp$ for everey $X \in \Sm/k$.
On the other hand we deduce that the base extension map
$\hEt \Pro^1_{\ok} \wedge \hEt X \to  \hEt\Pro^1_k \wedge_{\hEt k} \hEt X$ is a $\Zl$-weak equivalence in $\hShp$. Moreover, we have an isomorphism $\hEt \Pro^1_{\ok} \cong S^2$ in $\hShp$ by Example \ref{p1}, which implies the assertion of the lemma for $X \in \Sm/k$.\\ 
If $\Xh$ denotes a presheaf on $\Sm/k$, we can replace $X$ by $\Xh$ and get the same results. For, $\Xh$ is isomorphic to the colimit of representable presheaves $\Xh =\colim_s X_s$. Since $\Et $ commutes with colimits, we get
$H^{\ast}(\Et \Xh;\Zl) \cong \lim_s H^{\ast}(\Et X_s;\Zl)$. Since each $X_s$ is a smooth scheme over $k$, the \'etale cohomology groups $H^i(\Et X_s;\Zl) = H_{\et}^i(X_s;\Zl)$ are finite $\Zl$-vector spaces in each degree. Hence the limit over all $s$ of $H^{\ast}(\Et X_s;\Zl)$ commutes with the functor$H^{\ast}(\Et \Pro^1_k;\Zl) \otimes_{H^{\ast}(\Et k;\Zl)}-$. 
\end{proof}

Hence if $\sigma_n:\Pro^1_k \wedge_k E_n \to E_{n+1}$ is the structure map of a motivic $\Pro^1_k$-spectrum, then $\hEt$ yields a sequence of maps 
\begin{equation}\label{structuremap}
S^2 \wedge \hEt E_n \stackrel{\simeq}{\longrightarrow} \hEt (\Pro^1_k) \wedge_{\hEt k} \hEt( E_n)
\stackrel{\simeq}{\longleftarrow} \hEt(\Pro^1_k \wedge_k E_n) \stackrel{\hEt \sigma_n}{\longrightarrow}
E_{n+1}
\end{equation}
where the first two maps are weak equivalences in $\hShpk$. Unfortunately, the map in the middle points to the wrong direction. This leads to the following construction.\\
Since there is no natural inverse map $\hEt (\Pro^1_k) \wedge_{\hEt k} \hEt(\Xh)
\to \hEt(\Pro^1_k \wedge_k \Xh)$ in $\hShpk$, we may only construct an \'etale realization on the level of the homotopy categories. Therefore, we consider an intermediate category $\Ch$ and deduce from a zig-zag of functors 
$$\mathrm{Sp}^{\Pro^1}(k) \stackrel{\tilde{\hEt}}{\to} \Ch \stackrel{i}{\hookleftarrow} \mathrm{Sp}(\hSh, S^2 \wedge \cdot)$$
the existence of a stable realization functor $\SHh^{\Pro^1}(k) \to \hSHh_2$.\\

In view of Lemma \ref{Kunneth}, it is natural to consider the following definition. The objects of the category $\Chk$ are sequences
$$\{ F_n, F'_n, F''_n; S^2 \wedge F_n \stackrel{\simeq p_n}{\longrightarrow} F'_n \stackrel{\simeq q_n}{\longleftarrow} F''_n \stackrel{r_n}{\longrightarrow} F_{n+1} \}_{n \in \Z}$$
where $F_n$, $F'_n$, $F''_n$ are pointed profinite spaces over $\hEt k$ and $p_n$, $q_n$ and $r_n$ are maps in $\hShpk$; furthermore the maps $p_n$ and $q_n$ are weak equivalences in $\hShp$.\\ 
The morphisms of $\Chk$ are levelwise morphisms of $\hShpk$ which make the obvious diagrams commute, where the map $S^2 \wedge E_n \to S^2 \wedge F_n$ is the map induced by $E_n \to F_n$. 
The functor $$\mathrm{Sp}(\hShpk, S^2 \wedge \cdot)\stackrel{i}{\hookrightarrow} \Chk$$ denotes the full embedding which sends $\{ F_n, S^2 \wedge F_n \stackrel{\sigma_n}{\to} F_{n+1} \}$ to  $\{F_n, S^2 \wedge F_n, S^2 \wedge F_n; S^2 \wedge F_n \stackrel{\mathrm{id}}{\to} S^2 \wedge
F_n \stackrel{\mathrm{id}}{\leftarrow} S^2 \wedge F_n \stackrel{\sigma_n}{\to} F_{n+1} \}$.\\ 
When we apply $\hEt$ levelwise, we get a functor $$\mathrm{Sp}^{\Pro^1}(k) \stackrel{\tilde{\hEt}}{\to} \Chk.$$
We define a class $W$ of maps in $\Chk$ as the image of the stable equivalences of $\mathrm{Sp}(\hShpk, S^2 \wedge \cdot)$ under the embedding $i$. 
Since the maps in $W$ are the images of weak equivalences in a model structure and since $i$ is a full embedding, it is clear that $W$ admits a calculus of fractions. This ensures that we may form the localized category $\Ho(\Chk):=\Ch[W^{-1}]$ in which the maps in $W$ become isomorphisms. We call the maps in $W$ weak equivalences or stable equivalences, by abuse of notation. We will call a map in $W$ a level equivalence if it is in the image of the level equivalences of $\mathrm{Sp}(\hShpk, S^2 \wedge \cdot)$ under $i$.\\ 
Furthermore, $i$ sends stable equivalences into weak equivalences by definition and hence it induces a functor on the homotopy categories, which we also denote by $i$. 

\begin{prop}\label{equivofcat}
The induced functor $i: \hSHhtk \to \Ho(\Chk)$ is an equivalence of categories. In particular, there is an inverse equivalence $j: \Ho(\Chk) \to \hSHhtk$.
\end{prop}

\begin{proof}
By \cite{MacLane}, IV, 4, Theorem 1, since $i$ is a full embedding, it suffices to show that for every $F \in \Ho(\Chk)$ there is a spectrum $E \in \hSHhtk$ such that $i(E)\cong F$ in $\Ho(\Chk)$. This implies already that $i$ is an equivalence and that there is an inverse functor $j$.\\
 The crucial point is to construct the structure map of a spectrum from the given data of $F$. For the convenience of notations, we will omit in the following proof the structure maps of spaces over $\hEt k$.\\ 
Let $R$ be a fixed fibrant replacement functor in $\hShpk$. We consider the category $\mathrm{Sp}(\hShpk, RS^2 \wedge \cdot)$ as a Quillen equivalent model for $\hSHhtk$. We apply $R$ on each level. Since $R$ commutes with products, we get the following sequence
\begin{equation}\label{sequence}
RS^2 \wedge RF_n \stackrel{\simeq Rp_n}{\longrightarrow} RF'_n \stackrel{\simeq Rq_n}{\longleftarrow} RF''_n \stackrel{Rr_n}{\longrightarrow} RF_{n+1}.
\end{equation}
In addition, the functor $R$ can be chosen such that the map $Rq_n$ is a trivial fibration between fibrant and cofibrant objects, see Proposition 8.1.23 of \cite{hirsch}. The sequence we get is still isomorphic in $\Ho(\Chk)$ to the initial one since they are level equivalent. By Proposition 9.6.4 of \cite{hirsch}, there is a right inverse $s_n$ of $Rq_n$ in $\hShpk$ such that $Rq_ns_n=\id_{RF'_n}$ and a homotopy $s_nRq_n \sim \id_{RF''_n}$. We denote by $E$ the resulting spectrum with structure maps $\sigma_n:=Rr_n\circ s_n\circ Rp_n$.\\ 
In order to check that $i(E)$ is isomorphic to $F$ in $\Ho(\Chk)$, we consider the following diagram
\begin{equation}
\begin{array}{ccc}
RS^2 \wedge RF_n & \stackrel{\id}{\longrightarrow} & RS^2 \wedge RF_n \\
\id \downarrow &   & \downarrow Rp_n \\
RS^2 \wedge RF_n & \stackrel{Rp_n}{\longrightarrow} & RF'_n \\
\id \uparrow &   & \uparrow Rq_n \\
RS^2 \wedge RF_n & \stackrel{s_n\circ Rp_n}{\longrightarrow} & RF''_n\\
\sigma_n \downarrow &   & \downarrow Rr_n \\
RF_{n+1} & \stackrel{\id}{\longrightarrow} & RF_{n+1}
\end{array}
\end{equation}
representing the $n$-th level of the canonical map $i(E) \to F$ in $\Ho(\Chk)$. The upper square commutes obviously. The lower square commutes by the definition of $\sigma_n:=Rr_n\circ s_n\circ Rp_n$. The middle square commutes by the construction of $s_n$ such that $Rq_n s_n =\id_{RF'_n}$. Hence this is in fact a morphism in $\Chk$. Since all horizontal maps are weak equivalences, the morphism is a level equivalence in $\Chk$ and hence it is an isomorphism in $\Ho(\Chk)$.
\end{proof}


We have to show that  $\thEtsp:\mathrm{Sp}^{\Pro^1}(k) \to \Chk$ has a total left derived functor. Therefor, we have to choose the right model for the stable motivic category $\SHhP$. By \cite{dugger}, we know that $\LU$ is a left proper cellular simplicial model category which allows us to apply the methods of \cite{hovey}.

\begin{prop}\label{p1model}
The canonical functors 
$$\mathrm{Sp}^{\Pro^1}(k) \to \SpMVP \to \SpLUP$$
are Quillen equivalences and we get
$$\SHhP \simeq \Ho^{\mathrm{stable}}(\SpMVP) \simeq \Ho^{\mathrm{stable}}(\SpLUP).$$ 
\end{prop}
\begin{proof}
The first equivalence follows from \cite{hovey}, Corollary 3.5, and the obvious fact that the cofibrations are mapped to cofibrations and the fibrant objects correspond to each other in both categories.\\
The second equivalence follows from \cite{hovey}, Theorem 5.7, taking into account that $\LU$ is Quillen equivalent to $\MVh$ and that all objects in $\MVh$ are cofibrant.
\end{proof}

We use $ \SpLUP$ as a model for $\SHhP$. 
\begin{theorem}\label{P1real}
The functor $\hEt$ induces an \'etale realization of the stable motivic homotopy category of $\Pro^1$-spectra: $$\LhEt: \SHhP \to \hSHhtk$$
defined to be the composite $\LhEt : \SHh(k) \stackrel{\LthEt}{\longrightarrow} \Ho(\Chk) 
\stackrel{j}{\longrightarrow} \hSHhtk$.
\end{theorem}
\begin{proof}
We have to show that stable equivalences in $\SpLUP$ are sent to isomorphisms in $\Ho(\Chk)$.\\
We know that $\thEtsp$ sends level weak equivalences between cofibrant objects in $\LU$ to weak equivalences in $\Chk$, since $\hEt$ sends weak equivalences between cofibrant objects to weak equivalences in $\hShpk$. Hence it induces a total left derived functor on the projective model structure of $\SpLUP$.\\
We use the notation $\Sigma^{\Pro^1}_n:\LU \to \SpLUP$ for the left adjoint to the $n$-th evaluation functor. It is given by $(\Sigma^{\Pro^1}_n \Xh)_m = (\Pro^1)^{m-n}\Xh$ if $m\geq n$ and $(\Sigma^{\Pro^1}_n \Xh)_m = \Spec k$ otherwise.\\ 
We denote by $F_n:\hShpk \to \Chk$ the composition of the corresponding functor $\widehat{\Sigma}_n:\hShpk \to \mathrm{Sp}(\hShpk,S^2\wedge \cdot)$ followed by the embedding 
$i:\mathrm{Sp}(\hShpk,S^2\wedge \cdot) \to \Chk$. \\ 
In order to show that there exists a derived functor on the stable structure it suffices to show that  $\hEtsp(\zeta^{\Xh}_n)$ is a stable equivalence for maps $\zeta^{\Xh}_n:\Sigma_{n+1}\Pro^1\wedge \Xh \to \Sigma_n \Xh$ in $ \SpLUP$ for all cofibrant presheaves $\Xh \in \LU$, since these are the maps at which we localize for the stabilization.\\ 
We consider the commutative diagram in $\Ho(\Chk)$
$$\begin{array}{cccc} 
\zeta^{\hEt \Xh}_n: & F_{n+1}(S^2\wedge \hEt \Xh) & \stackrel{\cong}{\longrightarrow} & F_n(\hEt \Xh)\\
  & \updownarrow \cong &  & \parallel \\
       & F_{n+1}(\hEt \Pro^1_k\wedge_{\hEt k} \hEt \Xh))  & 
        & F_n(\hEt \Xh)\\
  & \updownarrow \cong &  & \parallel \\
       & F_{n+1}(\hEt (\Pro^1 \wedge \Xh))  & 
        & F_n(\hEt \Xh)\\
 & \cong \updownarrow &  & \updownarrow \cong\\
\hEtsp\zeta^{\Xh}_n: & \hEtsp (\Sigma^{\Pro^1}_{n+1}(\Pro^1\wedge \Xh)) &
\longrightarrow & \hEtsp (\Sigma^{\Pro^1}_n \Xh).  
\end{array}$$
The upper and middle vertical isomorphism on the left hand side are given by the obvious level equivalences given by the canonical sequence of weak equivalences as in Lemma \ref{Kunneth}. The lower vertical isomorphisms are given by the following composition:\\
We discuss the isomorphism $F_n(\hEt \Xh) \cong \hEtsp(\Sigma^{\Pro^1}_n \Xh)$. The lower left hand isomorphism is constructed in the same way. We define an intermediate object $E^n \in \Chk$ given in degree $m\geq n$ by 
$$E^n_m= (\hEt \Pro^1_k)^{m-n}\wedge_{\hEt k} \hEt \Xh  , E_m^{n'}= (\hEt \Pro^1_k)^{m+1-n}\wedge_{\hEt k} \hEt \Xh, E_m^{n''}= E_m^{n'}$$ with the obvious strcuture maps induced by $S^2 \to \hEt \Pro^1_k$ respectively the identity; in degree $m<n$ it is defined by $E_m^n=E_m^{n'}=E_m^{n''}=\hEt k$ with identity maps and the map to the terminal object $\hEt k$ of $\hShpk$.\\
The object $E^n$ is defined such that there are canonical maps
$$F_n(\hEt \Xh) \stackrel{\alpha}{\longrightarrow} E^n \stackrel{\beta}{\longleftarrow} \hEtsp(\Sigma^{\Pro^1}_n \Xh)$$ induced in degree $m\geq n$ by the maps
$$(S^2)^{m-n}\wedge \hEt \Xh \longrightarrow (\hEt \Pro^1_k)^{m-n}\wedge_{\hEt k} \hEt \Xh 
\longleftarrow \hEt ((\Pro^1_k)^{m-n}\wedge_k \hEt \Xh).$$
One checks easily using the canonical weak equivalences of Lemma \ref{Kunneth} that the maps $\alpha$ and $\beta$ are both level equivalences in $\Chk$. Hence $\alpha$ and $\beta$ are isomorphisms in $\Ho(\Chk)$. Their composition is the isomorphism $F_n(\hEt \Xh) \cong \hEtsp(\Sigma^{\Pro^1}_n \Xh)$.\\  
We deduce from the diagram that $\zeta^{\hEt \Xh}_n$ and $\hEtsp\zeta^{\Xh}_n$ differ only by an isomorphism in $\Ho(\Chk)$ given by level equivalences. Since $\zeta^{\hEt \Xh}_n$ is a stable equivalence in $\Chk$ by definition, the map $\hEtsp\zeta^{\Xh}_n$ is in fact an isomorphism in $\Ho(\Chk)$.\\
The last point is to show that $\hEtsp$ factors through $\A^1$-weak equivalences. By Lemma \ref{lflocalized} and \cite{pi0} it suffices to show that  $\hEtsp$ sends the maps $\Sigma^{\infty}_{\Pro^1}((X\times\A^1)_+)[n] \to \Sigma^{\infty}_{\Pro^1}(X_+)[n]$ to isomorphisms for all $X\in \mathrm{Sm}/k$ and all $n \in \Z$.\\ 
We know that $\hEt(X\times \A^1) \to \hEt (X)$ is a weak equivalence in $\hSh$ for $X \in \Sm/k$. 
Hence $\hEtsp \Sigma^{\infty}_{\Pro^1}((X\times\A^1)_+)[n] \to \hEtsp \Sigma^{\infty}_{\Pro^1}(X_+)[n]$ is a level weak equivalence in $\Chk$. 
This shows that $\hEtsp$ induces a total derived functor on the category $\SHhP = \mathcal{SH}^{\Pro^1}_{\A^1}(k)$ of $\A^1$-localized $\Pro^1$-spectra. 
\end{proof}

As in Theorem \ref{LhEt}, since the suspension spectrum of a smooth scheme is a stable cofibrant object, we have $\LhEt = \hEt$ on such spectra and $\hEt$ sends stable $\A^1$-equivalences to the image in $\Chk$ of stable equivalences in $\hSHhtk$. We summarize this discussion in the following 


\subsection{Etale realization of $S^1$-spectra}

The \'etale realization of $S^1$-spectra is essentially simpler. The structure maps $\sigma_n:S^1_k \wedge_k E_n \to E_{n+1}$ induce by base extension canonical maps $S^1_{\ok} \wedge E_n \to S^1_k \wedge_k E_n \stackrel{\sigma_n}{\to} E_{n+1}$ in $\hShpk$. By Example \ref{hEta}, we have a canonical isomorphism $\hEt S^1_{\ok} = S^1$ in $\hShp$ and by Formula (\ref{s1smash}) we have an isomorphism $S^1 \wedge \hEt E_n \stackrel{\cong}{\to} \hEt (S^1_{\ok} \wedge E_n)$ in $\hShpk$. This implies that the functor $\hEt$ on presheaves induces a functor on $S^1$-spectra 
$$\hEt : \mathrm{Sp}^{S^1}(k) \to \hSpk$$ by sending a spectrum $E$ to the profinite spectrum $\hEt E$ given in degree $n$ by $(\hEt E)_n:=\hEt (E_n)$ with structure maps 
$$S^1 \wedge \hEt E_n = \hEt (S^1_{\ok} \wedge E_n) \stackrel{\simeq}{\longrightarrow} \hEt (S^1_k \wedge_k E_n) \stackrel{\hEt \sigma_n}{\longrightarrow} \hEt E_{n+1}$$ 
where $\sigma$ is the structure map of the given spectrum $E$ and the map in the middle is obviously a $\Zl$-weak equivalence in $\hShpk$.\\
The proof of the following theorem is essentially the same as for the previous one on $\Pro^1$-spectra with the obvious simplifications. 

\begin{theorem}
The functor $\hEtsp : \mathrm{Sp}^{S^1}(k) \to \hSpk$ admits a left derived
functor $\LhEtsp : \SHhS \to \hSHhk$ on the stable $\A^1$-homotopy category of
$S^1$-spectra. 
\end{theorem}

\subsection{Examples}

Let $MGL$ denote the motivic Thom spectrum defined in \cite{a1hom} representing algebraic cobordism.

\begin{theorem}\label{EtofMGL}
Let $k$ be a separably closed field. There is an isomorphism in $\hSHh_2$ 
$$\LhEt (MGL) \cong \hMU.$$ 
\end{theorem}
\begin{proof}
Let $G_k(n,N)$ be the Grassmannian over $k$ and $G_k(n,N)$ be the colimit over $N$. By the Thom isomorphisms in \'etale and singular cohomology, it suffices to show that there is a weak equivalence between $\hEt (G_k(n,N))$ and $\hat{G}_{\C}(n,N)$, the profinite completion of the corresponding simplicial set of the complex Grassmannian manifold.\\
If $\car k >0$ let $R$ be the ring of Witt vectors of $k$, otherwise $R=k$. By choosing a common embedding of $R$ and $\C$ into an algebarically closed field, Friedlander proved in \cite{etaleK}, 3.2.2, that there is a natural sequence of $\Zl$-weak equivalences in $\pro-\Sh_{\ast}$ 
$$\mathrm{Sing}(G_{\C}(n,N))  \to \Et(G_{\C}(n,N)) \to \Et(G_R(n,N)) \leftarrow \Et(G_k(n,N)),$$
where one should recall that weak equivalences in $\pro-\Hh_0$ in the sense of \cite{etaleK} correspond to $\Zl$-weak equivalences in $\pro-\Sh_{\ast}$ in the sense of \cite{compofpro}.\\
By taking colimits with respect to $N$ we get a sequence of $\Zl$-weak equivalences in $\pro-\Sh_{\ast}$ 
$$\mathrm{Sing}(G_{\C}(n)) \to \Et(G_{\C}(n))  \to \Et(G_R(n)) \leftarrow \Et(G_k(n)).$$ 
Since $\Zl$-weak equivalences are preserved under completion, this shows that there is an isomorphism $\hat{G}_{\C}(n) \cong \hEt (G_k(n))$ in $\hHhp$. Together with the Thom isomorphism, this proves that there is an isomorphism of spectra $\LhEt (MGL) \cong \hMU$ in $\hSHh_2$.   
\end{proof}

\section{Profinite \'etale cobordism}

The study of profinite \'etale cohomology theories is the main purpose of this paper. The idea is to apply generalized profinite cohomology theories represented by a profinite spectrum to the functor $\hEt$.\\
We show that every such profinite \'etale cohomology theory is in fact a cohomology theory on $\Sm/k$ in the sense of \cite{panin}. One should note that in order to get the $\A^1$-invariance, we have to complete away from characteristic of $k$ and cannot use a $\pi_{\ast}$-model structure, since $\pi_1^{\et}(\A^1_k)$ is non-trivial over a field of positive characteristic, see \cite{schmidt}.\\
We show that the naturally arising profinite \'etale K-theory represented by $\hKU$ for schemes over a separably closed field is isomorphic to the \'etale K-theory of Friedlander \cite{etaleK}.\\
The second theory we study is profinite \'etale cobordism represented by $\hMU$. The coefficients of this theory are determined by our knowledge about profinite cohomology theories. We deduce an Atiyah-Hirzebruch spectral sequence from Theorem \ref{ahss} for \'etale cobordism with finite coefficients $\hMU\Zln$ starting with \'etale cohomology.\\
We show in the last section that $\hMU_{\et}$ is in fact an oriented cohomology theory on $\Sm/k$ if $k$ is separably closed. In order to prove this, we have to use non-trivial results on $\hMU$ to be able to reduce the problems to questions in \'etale cohomology.\\
The ideas for the proof that $\hMU_{\et}$ is an oriented cohomology theory are not new. For our purpose we use the techniques of \cite{panin}. We define an orientation of $\hMU_{\et}$ with respect to smooth schemes over $k$. We use this orientation and a projective bundle formula to construct a Chern structure on $\hMU_{\et}$. The key ingredient for the proof of the projective bundle formula is the Atiyah-Hirzebruch spectral sequence of $\hMU$. From facts about \'etale cohomology and the results of \cite{panin}, we deduce that \'etale cobordism is an oriented cohomology theory on smooth schemes over a separably closed field $k$ in the sense of \cite{panin} and hence in the sense of \cite{lm}.\\
This result will enable us to compare profinite \'etale cobordism with algebraic cobordism in the next section.

For the convenience of the reader we recall some notations and definitions: we denote by $\hSh$, resp. $\hShp$, the category of simplicial profinite sets, resp. pointed simplicial profinite sets; $\hHh$, resp. $\hHhp$, is the corresponding homotopy category; $\hSp$ is the category of profinite spectra and $\hSHh$ denotes the stable homotopy category of profinite spectra; we write $\Sigma^{\infty}(Y)$ for the suspension spectrum of a profinite space $Y$, see Sections 2 and 3.\\ 
Every profinite spectrum $E$ determines a generalized cohomology theory $E^{\ast}(-)$ on profinite spaces, e.g. $\hMU$ and $\hMU\Zln$ yield profinite cobordism theories $\hMU^{\ast}(-)$ and $\hMU^{\ast}(-;\Zln)$, see Section 4.\\
The functor $\hEt: \Sch/k \to \hSh$ is the profinitely completed \'etale topological type functor, see Section 5; we denote by $\LhEt: \SHhP \to \hSHh$ the stable \'etale realization functor, see Section 6 and Theorem \ref{P1real}.
    
\subsection{Profinite \'etale cohomology theories}

Let $k$ be a fixed base field and let $\ell$ be a fixed prime different from the characteristic of $k$. We consider the category $\Sch /k$ of schemes of finite type over $\Spec k$. This implies in particular that all schemes are noetherian and we may apply the \'etale topological type functor. 
\begin{defn}\label{profinetalecohom}
Let $E\in \hSp$ be a profinite spectrum and let $X \in \Sch/k$. \\
1. We define the profinite \'etale cohomology of $X$ in $E$ to be the profinite cohomology theory represented by $E$ applied to the profinite space $\hEt X$, i.e. $$E^n_{\et}(X):=E^n(\hEt X)=\Hom_{\hSHh}(\Sigma^{\infty}(\hEt X), E[n]),$$
where we add a base-point if $X$ is not already pointed.\\ 
2. We define relative \'etale cohomology groups $E^n_{\et}(X,U)$ for an open subscheme $U\subset X$ by
$$E^n_{\et}(X,U):= E^n(\hEt (X)/\hEt (U)).$$
3. We define the \'etale cohomology group
$$E_{\et}^n(\Pro^{\infty}_k):= \lim_V E_{\et}^n(\Pro(V))$$ to be the limit over all finite dimensional vector spaces $V \hookrightarrow k^{\infty}$ of the just defined groups $E_{\et}^n(\Pro(V))$. 
\end{defn}

\begin{remark}\label{cohomofPinfty}
Since $\hEt$ is also a functor on spaces, we could have defined $E_{\et}^n(\Pro^{\infty}_k)$ just as for schemes to be the group $E^n(\hEt \Pro^{\infty}_k)$. But we will use results that are based on the first definition. Since $\hEt$ does not commute with arbitrary colimits, the two definitions differ in general. But there is a canonical map $E^{\ast}(\hEt \Pro^{\infty}) \to E^{\ast}_{\et}(\Pro^{\infty})$. 
\end{remark}

\begin{prop}\label{etaleEM}
If $\pi$ is a finite abelian group, then the \'etale cohomology $H^{\ast}(X;\pi)$ with constant coefficients $\pi$ is equal to the profinite \'etale cohomology theory represented by the Eilenberg-MacLane spectrum $H\pi$.
\end{prop} 
\begin{proof}
This follows directly from Remark \ref{remarketalecohom} and Proposition \ref{Kpin}.
\end{proof}

\begin{prop}\label{mayervietoris}
Let $U \stackrel{i}{\hookrightarrow} X \stackrel{p}{\leftarrow} V$ induce an elementary distinguished square, cf. (\ref{dsquare}). Let $E$ be a profinite spectrum. Then there is a Mayer-Vietoris long exact sequence of graded groups
$$\cdots E^n_{\et}(X)\to E^n_{\et}(U)\oplus E^n_{\et}(V)\to E^n_{\et}(U\times_X V) \to E^{n+1}_{\et}(X)\cdots.$$  
\end{prop}
\begin{proof}
This follows directly from Theorem \ref{excision}.
\end{proof} 

\begin{prop}\label{relativesequence}
Let $E$ be a profinite spectrum. Let $U\subset X$ be an open subscheme of $X$. We get a long exact sequence of cohomology groups 
$$\begin{array}{rcl}
E^{\ast}(\hEt X) \stackrel{j^{\ast}}{\to} & E^{\ast}(\hEt U) \stackrel{\partial}{\to}  & E^{\ast}(\hEt (X)/\hEt (U)) \stackrel{i^{\ast}}{\to}\\
 \stackrel{i^{\ast}}{\to} E^{\ast+1}(\hEt X) & \stackrel{j^{\ast}}{\to} E^{\ast+1}(\hEt U)
\end{array}$$
where $j:\hEt U \hookrightarrow \hEt X$ and $i:(\hEt X,\oslash) \hookrightarrow (\hEt X, \hEt U)$ denote the natural induced inclusions.
\end{prop}
\begin{proof}
Since $\hEt U \hookrightarrow \hEt X \to \hEt(X)/\hEt(U)$ is isomorphic to a cofiber sequence in $\hShp$, this is just the usual long exact sequence of $\Hom$-groups  induced by a cofiber sequence in a simplicial model category, see \cite{homalg}.
\end{proof}

On the category $\Sm/k$ of smooth schemes of finite type over $k$, every \'etale cohomology theory satisfies homotopy invariance. 

\begin{prop}\label{hominv}
Let $E$ be a profinite spectrum. Let $X \in \Sm/k$. The projection $p:X\times \A^1 \to X$ induces an isomorphism $p^{\ast}:E_{\et}^{\ast}(X) \stackrel{\cong}{\to} E_{\et}^{\ast}(X\times \A^1)$.
\end{prop}
\begin{proof}
This is clear since $p$ induces a weak equivalence in $\hSh$, see Theorem \ref{LhEt}, and hence induces isomorphisms on cohomology theories.
\end{proof}

More generally, this implies  
\begin{cor}\label{EH}
Let $E$ be a profinite spectrum. Let $V \to X$ be an $\A^n$-bundle over $X$ in $\Sm/k$. Then $p^{\ast}:E_{\et}^{\ast}(X) \stackrel{\cong}{\to} E_{\et}^{\ast}(V)$ is an isomorphism. 
\end{cor}
\begin{proof}
This follows immediately from the Mayer-Vietoris sequence of Proposition \ref{mayervietoris} and the previous proposition.
\end{proof}

In \cite{panin}, Definition 2.0.1, Panin gives the definition of a  cohomology theory on $\Sm/k$. The following theorem states that we have constructed a way to define such cohomology theories starting with a profinite spectrum and applying its associated cohomology theory to $\hEt$ on $\Sm/k$.
\begin{theorem}\label{Ecohom}
Let $E$ be a profinite spectrum. The \'etale cohomology theory $E_{\et}^{\ast}(-)$ represented by $E$ satisfies the axioms of a cohomology theory on $\Sm/k$ of \cite{panin}.
\end{theorem}
\begin{proof}
We check the axioms of a cohomology theory in the sense of \cite{panin}, Definition 2.0.1.\\ 
1. Localization: This is clear from Remark \ref{relativesequence}.\\
2. Excision: Let $e:(X',U') \to (X,U)$ be a morphism of pairs of schemes in $\Sm/k$ such that $e$ is \'etale and for $Z=X-U$, $Z'=X' -U'$ one has $e^{-1}(Z)=Z'$ and $e:Z' \to Z$ is an isomorphism. By \cite{milne} III, Proposition 1.27, we know that the morphism $e$ induces an isomorphism in \'etale cohomology $H^{\ast}(\hEt(X)/\hEt(U);\Zl) \cong H^{\ast}(\hEt(X')/\hEt(U');\Zl)$. Hence $\hEt(X')/\hEt(U') \to \hEt(X)/\hEt(U)$ is an isomorphism in $\hHhp$. Therefore, it induces the desired isomorphism 
$$E^{\ast}(\hEt(X)/\hEt(U)) \cong E^{\ast}(\hEt(X')/\hEt(U'))$$ for every \'etale cohomology theory.\\
3. Homotopy invariance: This is the content of Proposition \ref{hominv}.
\end{proof}

Later we will prove that profinite \'etale cobordism is in fact an {\em oriented} cohomology theory in the sense of \cite{panin} and of \cite{lm}.

\subsection{Profinite \'etale K-theory}

As a first application, we consider a comparison statement for profinite K-theory and the pro-K-theory of \cite{etaleK}, and hence for profinite \'etale K-theory and for Friedlander's \'etale K-theory. Let $\{X_s\} \in \pro-\Sh$ be a pro-object in $\Sh$ and let $BU$ be the simplicial set representing complex K-theory. Friedlander defines the K-theory of $\{X_s\}$ for $\epsilon=0,1$ and $k>0$ by
$$K^{\epsilon}(\{X_s\};\Zln)=\Hom_{\pro-\Hh}(\{\Sigma^{\epsilon}X_s\},\{P^nBU\wedge\Zln\})$$
where $\{P^nBU\}$ denotes the Postnikov tower of $BU$ considered as a pro-object in $\Sh$.
For a locally noetherian scheme $X$, Friedlander defines the \'etale K-theory by
$$K_{\et}^{\epsilon}(X;\Zln):=K^{\epsilon}(\Et X;\Zln).$$

Friedlander also establishes an Atiyah-Hirzebruch spectral sequence, \cite{etaleK}, Proposition 1.4,
$$E_2^p(\Zln)=H^n(\{X_s\};\Zln)\Longrightarrow K^{\ast}(\{X_s\};\Zln)$$
which is convergent provided that $E_2^p(\Zln)=0$ for $n$ sufficiently large or that $E_2^p(\Zln)$ is finite for each $p$. \\
Following Definition \ref{profinetalecohom} we set
\begin{defn}\label{profinetaleKU}
We define the {\rm profinite K-theory} of a profinite space $X$ to be the cohomology theory represented by the profinitely completed spectrum $\hKU$ and define {\rm profinite K-theory with $\Zln$-coefficients} to be represented by $\hKU\wedge \Zln$. The {\rm profinite \'etale K-theory} of a scheme $X$ over $k$, is hence defined by $$\hKU^i_{\et}(X):=\hKU^i(\hEt X)$$
and {\rm profinite \'etale K-theory with $\Zln$-coefficients} by 
$$\hKU^i_{\et}(X;\Zln):=\hKU^i(\hEt X;\Zln).$$   
\end{defn}

\begin{theorem}\label{KUahss}
Let $k$ be a field. For any smooth scheme $X$ over $k$ and any $\nu$ there is a strongly convergent spectral sequence $\{E_r^{p,q}\}$: 
$$E_2^{p,q} = H_{\et}^p(X;\Zln) \Longrightarrow \hKU_{\et}^{p+q}(X;\Zln).$$
\end{theorem} 
\begin{proof}
This is the spectral sequence of Theorem \ref{ahss} where we use the isomorphism 
$$H^{\ast}_{\et}(X;\Zln) \cong H^{\ast}(\hEt X;\Zln)$$ of Remark \ref{remarketalecohom} for the groups $(\hKU\Zln_{\geq 0})^q$ which are isomorphic to $\Zln$, where $\hKU\Zln_{\geq 0}$ denotes a connective covering. Since the \'etale cohomology groups $H^i_{\et}(X;\Zln)$ vanish for large $i$, see e.g. \cite{milne} VI, Theorem 1.1, the convergence also follows from Theorem \ref{ahss}.  
\end{proof}

Since $\hKU$ is an $\hOmega$-spectrum, we have for every profinite space $X$ 
$$\hKU^i(X;\Zln)\cong \Hom_{\hHhp}(X_+,\hKU_i\wedge \Zln).$$ 
Hence the functor $\compl:\pro-\Sh \to \hSh$ yields natural maps 
$$K^{\epsilon}(\{X_s\};\Zln) \to \hKU^{\epsilon}(\hat{X};\Zln)$$
for every $\{X_s\}\in \pro-\Sh$. Furthermore, by Corollary \ref{hatKU} and (\ref{torsequencepi}) we have $\pi_q(\widehat{KU}\wedge \Z/\ell^k)\cong \Zln$ for $q$ even and $\cong 0$ for $q$ odd. Now Proposition \ref{Z/ell-comparison} implies that we have, on the one hand, for every $\{X_s\}\in \pro-\Sh$ an isomorphism of the $E_2$-terms of the corresponding Atiyah-Hirzebruch spectral sequences and, on the other hand, a natural map between the abutments. Both maps are compatible. Hence we have the diagram 
$$\begin{array}{rcccc}
E_2^p(\Zln) & = & H^p(\{X_s\};\Zln) & \Longrightarrow & K^{\ast}(\{X_s\};\Zln)\\
 &   & \downarrow \cong &  & \downarrow\\
E_2^p(\Zln) & = & H^p(X;\Zln) & \Longrightarrow & \hKU^{\ast}(\hat{X};\Zln).
\end{array}$$
The two sequences are in particular both convergent if $H^p(\{X_s\};\Zln)\cong H^p(\hat{X};\Zln)$ either equals $0$ for $p$ sufficiently large or is finite for all $p$. Since the \'etale types of schemes of finite type over a separably closed field satisfy the convergence conditions of the theorem, we get the following 
\begin{theorem}\label{competalektheory}
For every scheme $X$ of finite type over a separably closed field $k$ the
\'etale $K$-theory groups $K_{\et}^{\ast}(X;\Zln)$ of \cite{etaleK} are isomorphic to the profinite \'etale
$K$-theory groups $\hKU_{\et}^{\ast}(X;\Zln)$ defined above. $\Box$
\end{theorem} 
%

\subsection{Profinite \'etale cobordism} 

Following Definition \ref{profinetalecohom} we set
\begin{defn}\label{profinetaleMU}
1. We define the {\rm profinite \'etale cobordism} of a scheme $X$ of finite type over $k$, to be the cohomology theory represented by the profinitely completed spectrum $\hMU$ applied to $\hEt X$, i.e. $$\hMU^{\ast}_{\et}(X):=\hMU^{\ast}(\hEt X),$$  and {\rm profinite \'etale cobordism with $\Zln$-coefficients}  to be the cohomology theory represented by the profinitely completed spectrum $\hMU\Zln$ applied to $\hEt$, i.e. $$\hMU^{\ast}_{\et}(X;\Zln):=(\hMU\Zln)^{\ast}(\hEt X).$$ 
2. We also have relative \'etale cobordism groups $\hMU^{\ast}_{\et}(X,U)$, resp. $\hMU^{\ast}_{\et}(X,U;\Zln)$, for an open subscheme $U\subset X$ defined by
$$\hMU^{\ast}_{\et}(X,U):= \hMU^{\ast}(\hEt (X)/\hEt (U))$$ and 
$$\hMU^{\ast}_{\et}(X,U;\Zln):= \hMU^{\ast}(\hEt (X)/\hEt (U);\Zln).$$ 
3. We define the \'etale cobordism group
$$\hMU_{\et}^{\ast}(\Pro^{\infty}_k):= \lim_V \hMU_{\et}^{\ast}(\Pro(V))$$ 
and similarly for $\hMU\Zln_{\et}$ as in Definition \ref{profinetalecohom}.\\
4. We define {\rm $\ell$-adic \'etale cobordism} to be the limit 
$$\hMU_{\et}^{\ast}(X;\Z_{\ell}):= \lim_{\nu} \hMU_{\et}^{\ast}(X;\Zln).$$
\end{defn}

\begin{remark}\label{exteriorproduct}
1. The exterior product $MU_n \times MU_m \to MU_{m+n}$ on $MU$ together with the diagonal map $\Delta_{\hEt X}:\hEt X \to \hEt X \times \hEt X$ yield a commutative ring structure on $\hMU_{\et}^{\ast}(X)$ and $\hMU^{\ast}_{\et}(X;\Zln)$.\\
2. We have introduced the notion of $\ell$-adic cobordism, since it is the \'etale cohomology theory to which $\ell$-adic \'etale cohomology converges, see Theorem \ref{MUahss}. We will not give any further comment on this theory, since most of its properties follow from those of $\hMU\Zln$. \\
One should note that it is not true in general that $\hMU_{\et}^{\ast}(X;\Z_{\ell})=\hMU_{\et}^{\ast}(X)$. 
\end{remark}

The morphism of profinite spectra $\hMU \to H\Zln$ induced by the orientation also yields maps of profinite \'etale cohomology theories: 
\begin{prop}\label{hMUtoHZln}
Let $k$ be a field. For every $X$ in $\Sm/k$, there are unique induced morphisms of cohomology theories 
$$\hMU^{\ast}_{\et}(X) \to H_{\et}^{\ast}(X;\Zln)$$
and 
$$\hMU^{\ast}_{\et}(X;\Zln) \to H_{\et}^{\ast}(X;\Zln). \Box$$ 
\end{prop}

\begin{prop}\label{coeffofhMUet}
Let $R$ be a strict local henselian ring, i.e. a local henselian ring with separably closed residue field. Then $$\hMU^{\ast}_{\et}(\Spec R)\cong MU^{\ast}\otimes_{\Z} \Z_{\ell}.$$
In particular, for a separably closed field $k$ we get 
$$\hMU^{\ast}_{\et}(k)\cong MU^{\ast}\otimes_{\Z} \Z_{\ell}~\mathrm{and}~
\hMU^{\ast}_{\et}(k;\Zln)\cong MU^{\ast}\otimes_{\Z} \Zln.$$
\end{prop}
\begin{proof}
This follows from Corollary \ref{coeffofcompletedMU} and Example \ref{SpecR}. The last assertion follows from the exact sequence (\ref{torsequencecohom}), since $MU^{\ast}$ has no torsion.
\end{proof}

\begin{theorem}\label{MUahss}
Let $k$ be a field. For every scheme $X$ in $\Sm/k$, there is a convergent spectral sequence $\{E_r^{p,q}\}$ with 
$$E_2^{p,q}= H_{\et}^p(X;\Zln \otimes MU^q) \Longrightarrow \hMU_{\et}^{p+q}(X;\Zln).$$
In addition, we have also a convergent spectral sequence for $\ell$-adic cobordism
$$E_2^{p,q}=\lim_{\nu}H_{\et}^p(X;\Zln \otimes MU^q)= H_{\et}^p(X;\Z_{\ell} \otimes MU^q) \Longrightarrow \hMU_{\et}^{p+q}(X;\Z_{\ell}).$$
\end{theorem} 
\begin{proof}
This is the spectral sequence of Theorem \ref{ahss} together with the isomorphisms  
$$H^{\ast}_{\et}(X;\Zln \otimes MU^q) \cong H^{\ast}(\hEt X;\Zln \otimes MU^q)$$ 
of Remark \ref{remarketalecohom} for the finite group $M=\Zln \otimes MU^q=\Zln \otimes \hMU^q$. 
The condition on the cohomology of $X$ for strong convergence in Theorem \ref{ahss} is satisfied since $H^p_{\et}(X;\Zln \otimes MU^q)$ are finite groups for all $p$ by \cite{milne} VI, Corollary 5.5.\\
For the second sequence, the fact that $H^p_{\et}(X;\Zln \otimes MU^q)$ are finite groups for all $p$ implies that exact couples remain exact after taking limits and hence that $E_2^{p,q}=\lim_{\nu}E_2^{p,q}(\Zln) = H_{\et}^p(X;\Z_{\ell} \otimes MU^q)$ converges to $\lim_{\nu} \hMU_{\et}^{p+q}(X;\Zln)=\hMU_{\et}^{p+q}(X;\Z_{\ell})$.
\end{proof}

\begin{remark}\label{remarkoncoefficients}
1. Since $\Zln \otimes MU^t$ is a finitely generated free $\Zln$-module, i.e. isomorphic to a direct sum of finitely many copies of $\Zln$, and since \'etale cohomology is compatible with direct limits, the group $H_{\et}^s(X;\Zln\otimes MU^t)$ is isomorphic to $H_{\et}^s(X;\Zln) \otimes MU^t$. Hence in order to calculate the groups $H_{\et}^s(X;\Zln\otimes MU^t)$, we may reduce the problem to determine the groups $H_{\et}^s(X;\Zln)$.\\
2. By Theorem \ref{ahss}, we get a spectral sequence  
$$E_2^{p,q}=H^p(\hEt X;\hMU^q) \Longrightarrow \hMU_{\et}^{p+q}(X).$$
But the problem is that, in general, we do not have an isomorphism $H^p(\hEt X;\hMU^q) \cong H_{\et}^p(X;\hMU^q)$, since $\hMU^q$ is infinite. In addition, the \'etale cohomology with $\Z_{\ell}$-coefficients does not have good properties. 
\end{remark}

\begin{theorem}\label{complexvar}
Let $X$ be an algebraic variety over $\C$. Let $X(\C)$ be the topological space of complex points. For every $\nu$, there is an isomorphism $$\hMU_{\et}^{\ast}(X;\Zln)\cong MU^{\ast}(X(\C);\Zln).$$
\end{theorem}
\begin{proof}
By Corollary \ref{piofMU} we know that 
$$(\hMU\Zln)^q \cong (MU\Zln)^q \cong MU^q \otimes \Zln.$$ 
Since $\pi_q(MU)$ is a finitely generated abelian group, we get that $MU^q \otimes \Z/\ell^k$ is a finite group. We denote by $\mathrm{Sing}(X(\C))$ the singular simplicial set of the  topological space $X(\C)$. By Theorem 8.4 and Corollary 8.5 of \cite{fried}, by Remark \ref{remarketalecohom} and Corollary \ref{piofMU} we get a map of spectral sequences and an isomorphism $$H^p(\mathrm{Sing}(X(\C));(MU\Zln)^q)\cong H^p(\hEt X;(\hMU\Zln)^q)$$ between singular and profinite \'etale cohomology for all $\nu$, $p$ and $q$.  Since $X$ is finite dimensional, the corresponding complex, respectively \'etale Atiyah-Hirzebruch spectral sequences are convergent with isomorphic $E_2$-terms.
\end{proof}

Also the following statements are easy consequences of the existence of the spectral sequence and results on \'etale cohomology.
\begin{prop}\label{properbasechange}
Let $A$ be a strict local henselian ring and $S=\Spec A$. Let $f:X \to S$ be a proper morphism and let $X_0$ be the closed fibre of $f$. Then the induced map 
$$\hMU_{\et}^{\ast}(X;\Zln) \stackrel{\cong}{\longrightarrow} \hMU_{\et}^{\ast}(X_0;\Zln)$$
is an isomorphism on \'etale cobordism with $\Zln$-coefficients.
\end{prop}
\begin{proof}
The morphism $X_0 \to X$ induces a morphism of Atiyah-Hirzebruch spectral sequences. By Theorem 1.2 of Chapter IV in \cite{sga412} Arcata, it induces an isomorphism on \'etale cohomology with $\Zln$-coefficients and hence, by Remark \ref{remarkoncoefficients} an isomorphism on \'etale cobordism with $\Zln$-coefficients.
\end{proof}

\begin{prop}\label{fieldextension}
Let $K/k$ be an extension of separably closed fields with characteristic different from $\ell$. Let $X$ be a $k$-scheme. The induced map 
$$\hMU_{\et}^{\ast}(X;\Zln) \stackrel{\cong}{\longrightarrow} \hMU_{\et}^{\ast}(X_K;\Zln)$$
is an isomorphism on \'etale cobordism with $\Zln$-coefficients. 
\end{prop}
\begin{proof}
Again the natural morphism $X_K \to X$ induces a map of spectral sequences and the assertion then follows from the isomorphism on \'etale cohomology with $\Zln$-coefficients of Corollaire 3.3 of Chapter V in \cite{sga412} Arcata and Remark \ref{remarkoncoefficients}. 
\end{proof}

\subsection{Profinite \'etale cobordism is an oriented cohomology theory}

We prove that $\hMU_{\et}^{2\ast}(-)$ and $\hMU_{\et}^{2\ast}(-;\Zln)$ are oriented cohomology theories in the sense of \cite{lm} on the category $\Sm/k$ of smooth quasi-projective schemes of finite type over a separably closed field $k$ with char $k\neq \ell$.\\ 
We have to check the axioms given in \cite{lm}. The axioms of a cohomology theory have already been shown for every \'etale cohomology theory. In order to check the axioms of an oriented cohomology theory for $\hMU_{\et}$, our strategy is as follows. We define an orientation class for $\hMU_{\et}$ and use it to define the first Chern class of a line bundle. Using the Chern class of the canonical line bundle on the projective $n$-space, we construct an explicit isomorphism for the cobordism of $\Pro^n_k$. The method to prove this isomorphism relies on the Atiyah-Hirzebruch spectral sequence for profinite cobordism. We use the spectral sequence to reduce the computation of the cobordism of $\Pro^n$ to the computation of the \'etale cohomology of $\Pro^n$ with $\Zl$-coefficients. This reduction uses non-trivial results on profinite cobordism and is the key point for the proof that $\hMU^{\ast}_{\et}$, resp. $(\hMU\Zln)^{\ast}_{\et}$, yields an oriented cohomology theory on $\Sm/k$.\\ 
This computation and the K\"unneth theorem for $\hMU$, resp. $\hMU\Zln$, yield a projective bundle formula. Then we use Grothendieck's method to define Chern classes for all vector bundles.\\
The general theory of cohomology theories with Chern classes given in \cite{panin} then implies all the remaining properties.\\
But one should note that the fact that $\hMU_{\et}$ is an oriented cohomology theory on $\Sm/k$ is a nontrivial result and does not just follow from the fact that $MU$ is an oriented theory on topological spaces. This is due to the fact that we may not use an embedding of $k$ in $\C$. We need an orientation with respect to smooth schemes over any separably closed field.\\   
To deal with the twists in \'etale cohomology we choose a primitive $\ell^{\nu}$-th-root of unity which yields an isomorphism $\Zln \cong \Zln(1)$, see Remarks \ref{twisting} and \ref{MGLet}.

\begin{defn}
Let $\Pro_k^{\infty}$ be the infinite projective space in the category $\LU$ of spaces over $k$ and let $E$ be a commutative ring spectrum in $\hSp$. Since we have an isomorphism $\hEt \Pro_k^1 \cong S^2$ in $\hHh$ by Example \ref{p1}, we have a canonical map 
$$E^2(\hEt \Pro_k^{\infty}) \to E^2(\hEt \Pro_k^1) \stackrel{\cong}{\to} E^2(S^2) \stackrel{\cong}{\to} \pi_0(E).$$
An {\rm orientation} for $E$  with respect to $\Sm/k$ is a class $x_E \in \tilde{E}^2(\hEt \Pro^{\infty})$ that maps to $1$ under the above map.
\end{defn}
\begin{example}\label{orientation}
1. For $\hMU$, resp. $\hMU\Zln$, we choose the orientation $x_{\hMU} \in \hMU^2(\hEt \Pro_k^{\infty})$, resp. $x_{\hMU\Zln} \in \hMU^2(\hEt \Pro_k^{\infty};\Zln)$, as the image of the orientation $x_{MGL}$ under $\hEt$. This means that we define the map $x_{\hMU}:\LhEt(\Pro^{\infty}) \to \hMU \wedge S^2$ in the preimage of $1 \in \pi_0(\hMU)$ such that the following diagram 
$$\begin{array}{rcc}
\hEt(\Pro^{\infty}) & \stackrel{\tilde{x}_{\hMU}}{\longrightarrow} & \hMU \wedge S^2\\
\LhEt(x_{MGL})\downarrow & \nearrow \phi & \\
\LhEt (MGL) \wedge \Pro^1
\end{array}$$
is commutative, where $\phi$ is the morphism defined in the next section using an \'etale Thom class.\\ Note that $x_{\hMU} \in \hMU^2(\hEt \Pro_k^{\infty})$ also defines an element  $\tilde{x}_{\hMU} \in \hMU_{\et}^2(\Pro_k^{\infty})$ via the canonical map $\hMU^{\ast}(\hEt \Pro^{\infty}) \to \hMU^{\ast}_{\et}(\Pro^{\infty})$, recall Definition \ref{profinetalecohom} and Remark \ref{cohomofPinfty}.\\
If $k=\C$, this orientation $x_{\hMU}$ is mapped to the classical complex orientation of $MU$ via the isomorphism of Proposition \ref{complexvar}. This is due to the fact that $x_{MGL}$ corresponds to the complex orientation of $MU$ under the complex realization functor, \cite{milnor}.\\
2. For $H\Zln$ an orientation $x_{H\Zln}$ is given by a generator of $H^2(\hEt \Pro_k^{\infty};\Zln) \cong \Zln$. This last isomorphism is deduced from 
$$H^2(\hEt \Pro_k^{\infty};\Zln) \cong H^2(\Et \Pro_k^{\infty};\Zln)\cong \lim_n H^2(\Et \Pro_k^n;\Zln) \cong \Zln,$$
where we use the fact that $\Et$ commutes with colimits.
\end{example}

The importance of an orientation of $\hMU$, resp. $\hMU\Zln$, is that it enables us to define the first Chern class of a line bundle.\\ 
A line bundle $L/X$ over a scheme $X$ of dimension $n$ over $k$ corresponds to a morphism $X \to \Pro_k^{n+1}$ unique up to $\A^1$-weak equivalences. Since $\LhEt$ agrees with $\hEt$ on smooth schemes by Theorem \ref{LhEt}, this also defines a unique map $\lambda:\hEt X \to \hEt \Pro_k^{n+1}$ in $\hHh$.\\ 
It is clear from the definition that the orientation $x_{\hMU}$ restricted to $\hEt \Pro_k^{n+1}$  corresponds to a map $x_{\hMU}:\hat{\Sigma}^{\infty}(\hEt \Pro_k^{n+1}) \to \hMU \wedge S^2$. The composition of these two maps is the first  Chern class of $L/X$. 
\begin{defn}\label{Chernclass}
A line bundle defines via composition a morphism in $\hSHh$
$$c_1(L): \hat{\Sigma}^{\infty}(\hEt X) \stackrel{\lambda}{\to} \hat{\Sigma}^{\infty}(\hEt \Pro_k^{n+1}) \stackrel{x_{\hMU}}{\to} \hMU \wedge S^2,$$
i.e. an element $c_1(L) \in \hMU_{\et}^2(X)$ which we call the {\rm first Chern class of $L$}.\\
Similarly, we define the first Chern class of a line bundle in $\hMU_{\et}^2(X;\Zln)$.
\end{defn}

Now we have to show that this Chern class for line bundles yields a Chern structure on \'etale cobordism. Therefore, we have to calculate the \'etale cobordism of the projective $n$-space $\Pro^n_k$ over $k$. The method we use is indirectly based on the Atiyah-Hirzebruch spectral sequence for profinite cobordism. We use the Atiyah-Hirzebruch spectral sequence to reduce the computation of $\hMU_{\et}^{\ast}(\Pro^n)$ to the well-known case of \'etale cohomology $H_{\et}^{\ast}(\Pro^n;\Zl)$ with $\Zl$- coefficients. A crucial point is the isomorphism 
\begin{equation}\label{hMUmod}
\hMU_{\et}^{\ast}(X)/(\ell \cdot \hMU^{\ast \leq0}) \cong H^{\ast}_{\et}(X;\Zl)
\end{equation}
for every $X$ in $\Sm/k$.

\begin{prop}\label{hMUofPn}
Let $k$ be a separably closed field. Let $\Pro^n_k$ be the projective space of dimension $n$ over $k$. Let $\Oh(1) \to \Pro^n_k$ be the canonical quotient line bundle. 
1.  Set $\xi :=c_1(\Oh(1)) \in \hMU_{\et}^2(\Pro^n_k)$. The \'etale cobordism of $\Pro^n_k$ is given by the isomorphism
$$\hMU_{\et}^{\ast}(\Pro^n_k)\cong \hMU^{\ast}[u]/(u^{n+1}), ~ \xi \mapsto u.$$ 
2. Set $\xi :=c_1(\Oh(1))\in \hMU_{\et}^2(\Pro^n_k;\Zln)$. The \'etale cobordism with $\Zln$-coefficients of $\Pro^n_k$ is given by the isomorphism
$$\hMU_{\et}^{\ast}(\Pro^n_k;\Zln)\cong (\hMU\Zln)^{\ast})[u]/(u^{n+1}), ~ \xi \mapsto u.$$ 
\end{prop}
\begin{proof}
We deduce both results from general facts about filtered $MU^{\ast}$-modules in \cite{dehon}. The crucial point is that $H^{\ast}(\hEt \Pro^n_k;\Zl^m)\cong H^{\ast}_{\et}(\Pro^n_k;\Zl^m)$ surjects onto $H^{\ast}(\hEt \Pro^n_k;\Zl)$ for every $m$, cf. \cite{milne} VI, 5.6. Hence $\hEt \Pro^n_k$ is a profinite space of finite dimension without $\ell$-torsion, see the discussion of Theorem \ref{KunnethforhMU}. \\
By Proposition 2.1.8 of \cite{dehon} respectively by definition, both $MU^{\ast}$-modules $M:=\hMU_{\et}^{\ast}(\Pro^n_k)$ and $N:=\hMU^{\ast}[u]/(u^{n+1})$ are objects of the category $\hat{\Lh}$ of free $MU^{\ast}$-modules that are complete with respect to the $\ell$-adic filtration $f^{\ast}$ of \cite{dehon}, \S 2.1.\\
By Corollaire 2.1.4 of \cite{dehon}, it suffices to show that $M/f^1 \cong N/f^1$ as $\Zl$-vector spaces. But since $\hEt \Pro^n_k$ is a finite dimensional profinite space, the isomorphism $M/f^1 \cong N/f^1$ follows from the isomorphism (\ref{hMUmod}), proved in Proposition 2.1.8 of \cite{dehon} together with Remark \ref{remarketalecohom}, and the isomorphism of \'etale cohomology $H^{\ast}_{\et}(\Pro^n_k;\Zl) \cong \Zl[u]/(u^{n+1})$. This proves the assertion for $\hMU$. One should remark that the proof of Proposition 2.1.8 of \cite{dehon} is based on the Atiyah-Hirzebruch spectral sequence for $\hMU$.\\ 
The second assertion follows from the first one using the fact that for a space $X$ without $\ell$-torsion, we have $\hMU^{\ast}(X;\Zln) \cong \hMU^{\ast}(X)\otimes \Zln$, see \cite{adams} III, Proposition 6.6.\\
We could deduce the second assertion also directly from the Atiyah-Hirzebruch spectral sequence for $\hMU\Zln_{\et}$ as in \cite{adams} I, Lemma 2.5, using the determination of $H_{\et}^{\ast}(\Pro^n_k;\Zln)$, since all differentials vanish.
\end{proof}

\begin{prop}\label{KunnethforPn}
For the projective space $\Pro^n_k$ over a separably closed field $k$ and every $X \in \Sm/k$, there are K\"unneth isomorphisms for the \'etale cobordism
$$\hMU_{\et}^{\ast}(X) \otimes_{\hMU^{\ast}} \hMU_{\et}^{\ast}(\Pro^n_k)
\stackrel{\cong}{\longrightarrow}\hMU_{\et}^{\ast}(X \times \Pro^n_k)$$ and 
$$\hMU_{\et}^{\ast}(X;\Zln) \otimes_{(\hMU\Zln)^{\ast}} \hMU_{\et}^{\ast}(\Pro^n_k;\Zln)
\stackrel{\cong}{\longrightarrow}\hMU_{\et}^{\ast}(X \times \Pro^n_k;\Zln).$$
\end{prop}
\begin{proof}
Again, the crucial point is that $H^{\ast}(\hEt \Pro^n_k;\Zl^m)\cong H^{\ast}_{\et}(\Pro^n_k;\Zl^m)$ surjects onto $H^{\ast}(\hEt \Pro^n_k;\Zl)$ for every $m$. Hence $\hEt \Pro^n_k$ is a profinite space without $\ell$-torsion. Over the separably closed field $k$, $H_{\et}^i(X;\Zln)$ and $H_{\et}^i(\Pro^n_k;\Zl)$ vanish for large $i$ and hence the profinite spaces $\hEt X$ and $\hEt \Pro^n_k$ are finite dimensional spaces in $\hSh$. This enables us to deduce the proposition from Theorem \ref{KunnethforhMU}.
\end{proof}

\begin{theorem}\label{PBF}{\rm Projective Bundle Formula}\\
Let $E \to X$ be a rank $n$ vector bundle over $X$ in $\Sm/k$, $\Oh(1) \to \Pro(E)$ the canonical quotient bundle with zero section $s:\Pro(E) \to \Oh(1)$. Set $\xi:= c_1(\Oh(1)) \in \hMU_{\et}^2(\Pro(E))$. Then $\hMU_{\et}^{\ast}(\Pro(E))$ is a free $\hMU_{\et}^{\ast}(X)$-module with basis $(1, \xi, \ldots , \xi^{n-1}).$\\
Similarly, $\hMU_{\et}^{\ast}(\Pro(E);\Zln)$ is a free $\hMU_{\et}^{\ast}(X;\Zln)$-module with basis $(1, \xi, \ldots , \xi^{n-1}),$ where $\xi:= c_1(\Oh(1)) \in \hMU_{\et}^2(\Pro(E);\Zln)$.
\end{theorem}
\begin{proof}
We consider $\hMU_{\et}$. We prove this theorem first for the case of a trivial bundle on $X$. We have to show that the map
$\hMU_{\et}^{\ast}(X \times \Pro^n_k) \stackrel{\cong}{\longrightarrow} \hMU_{\et}^{\ast}(X)[u]/(u^n)$, sending $\xi$ to $u$, is an isomorphism. By Proposition \ref{hEtofproduct}, we know that the canonical map $\hEt (X) \times \hEt (\Pro^n_k) \to \hEt (X\times \Pro^n_k)$ is a weak equivalence in $\hSh$. Proposition \ref{hMUofPn} yields an isomorphism 
$$\hMU_{\et}^{\ast}(X)[u] /(u^n)\stackrel{\cong}{\to} \hMU_{\et}^{\ast}(X)\otimes_{\hMU^{\ast}} \hMU_{\et}^{\ast}(\Pro^n_k).$$ 
Hence it suffices to prove that the projections induce an isomorphism 
$$\hMU_{\et}^{\ast}(X)\otimes_{\hMU^{\ast}} \hMU_{\et}^{\ast}(\Pro^n_k) \stackrel{\cong}{\to} 
\hMU_{\et}^{\ast}(X \times \Pro^n_k).$$
This follows from the previous proposition and finishes the proof of the theorem in the case that $E$ is a trivial bundle on $X$.\\
A similar argument using the second assertion of Proposition \ref{hMUofPn} proves the case of a trivial bundle for $\hMU\Zln_{\et}$.\\
For the general case, since $E$ is by definition locally trivial for the Zariski topology on $X$, it suffices to show that the theorem holds for $X$ if it holds for open subsets $X_0$, $X_1$ and $X_0 \cap X_1$, with $X=X_0 \cup X_1$. This point now follows from comparing the two Mayer-Vietoris-sequences  
$$\begin{array}{c}\cdots \to \hMU_{\et}^{i}(\Pro(E_0)) \oplus \hMU_{\et}^{i}(\Pro(E_1))\\
\to \hMU_{\et}^{i}(\Pro(E_{01})) \to \hMU_{\et}^{i+1}(\Pro(E))\to \cdots \end{array}$$
$$\begin{array}{c}\hMU^{\ast} [u]/(u^n)\otimes (\cdots \to \hMU_{\et}^{i}(X_0) \oplus \hMU_{\et}^{i}(X_1)\\
\to \hMU_{\et}^{i}(X_{01}) \to \hMU_{\et}^{i+1}(X)\to \cdots)\end{array}$$
where we have written $X_{01}$ for $X_0 \cap X_1$, $E_0$, $E_1$ and $E_{01}$ for the restrictions of $E$ to $X_0$, $X_1$ and $X_{01}$.\\
Again, a similar argument shows the general case for $\hMU\Zln_{\et}$.
\end{proof}

We use Grothendieck's idea to introduce higher Chern classes for vector bundles. 
\begin{defn}\label{chernclasses}
Let $E$ be a vector bundle of rank $n$ on $X$ in $\Sm/k$ and let $\xi:=c_1(\Oh(1))$ be as above. Then we define the {\rm $i$th Chern class of $E$} to be the element $c_i(E) \in \hMU_{\et}^{2i}(X)$ such that $c_0(E)=1$, $c_i(E)=0$ for $i>n$ and $\Sigma_{i=0}^n (-1)^i c_i(E)\xi^{n-i}=0$. By the projective bundle formula of Theorem \ref{PBF} these elements are unique.\\
Similarly, we define $c_i(E) \in \hMU_{\et}^{2i}(X;\Zln)$ by Theorem \ref{PBF}.
\end{defn}

\begin{lemma}\label{nilpotencyofc1}
1.  For each line bundle over a smooth $X$, the class $c_1(L) \in \hMU_{\et}^2(X)$, resp. $c_1(L) \in \hMU_{\et}^2(X;\Zln)$, is nilpotent.\\
2. The assignment of Chern classes is functorial with respect to the pullback, i.e. $f^{\ast}(c_i(E))=c_i(f^{\ast}E)$ for a bundle $E$ on $X$ and a morphism $f:Y\to X$ in $\Sm/k$.
\end{lemma}
\begin{proof}
1. This may be proved exactly as in \cite{panin}, Lemma 3.6.4.\\
2. This follows immediately from Theorem \ref{PBF}.
\end{proof}

\begin{theorem}\label{transfer}
The \'etale cobordism theories satisfy the axioms of an oriented ring cohomology theory on $\Sm/k$ of \cite{panin}, Definition 3.1.1. This implies that for every projective morphism $f:Y \to X$ of codimension $d$ there are natural transfer maps 
$$f_{\ast}:\hMU_{\et}^{\ast}(Y) \to \hMU_{\et}^{\ast+2d}(X)$$ and
$$f_{\ast}:\hMU_{\et}^{\ast}(Y;\Zln) \to \hMU_{\et}^{\ast+2d}(X;\Zln)$$
which satisfy all the axioms of a push-forward in an oriented cohomology theory of \cite{lm}. 
\end{theorem}
\begin{proof}
We have proved in Theorem \ref{Ecohom} that $\hMU_{\et}(-)$ and $\hMU_{\et}(-;\Zln)$ satisfy the axioms of a cohomology theory in the sense of \cite{panin}, Definition 2.0.1. By Remark \ref{exteriorproduct}, we know that $\hMU_{\et}(-)$, resp.  $\hMU_{\et}(-;\Zln)$, is also a ring cohomology theory.\\ 
Secondly, the orientation on $\hMU_{\et}(-)$, resp. $\hMU_{\et}^{\ast}(-;\Zln)$, and Theorem \ref{PBF} show that we have in fact an oriented theory with Chern classes in the sense of \cite{panin}, Definition 3.1.1, using Theorem 3.2.4 of \cite{panin}.\\
Now there is a general procedure, described in 4.2 and 4.3 of \cite{panin}, to define natural transfer maps for every projective morphism in $\Sm/k$. Note that since we require transfer maps only for projective maps, we only have to construct transfers for closed immersions of smooth schemes and the projection of the projective space $\Pro^n_X \to X$ over $X$. In Theorem 4.7.1 of \cite{panin} it is proved that the push-forward is independent of the decomposition of $f$.\\ 
Furthermore, also the compatibility of pullbacks and push-forwards for transversal morphisms is proved in Theorem 4.7.1 of \cite{panin}.\\ 
That $f_{\ast}$ shifts the degree by $2d$ comes from the fact that Chern classes live in even degrees and the fact that the construction of $f_{\ast}$ depends on Chern classes.\\ 
This finishes the proof of the theorem.
\end{proof}

We state a general consequence of this extra structure on \'etale cobordism. 
\begin{cor}\label{Gysinsequence}
For every closed embedding of smooth varieties $i:Z \hookrightarrow X$ with open complement $j:U\hookrightarrow X$ there are exact sequences
$$\hMU_{\et}^{\ast}(Z) \stackrel{i_{\ast}}{\longrightarrow} \hMU_{\et}^{\ast}(X) \stackrel{j^{\ast}}{\longrightarrow} \hMU_{\et}^{\ast}(U)$$ and
$$\hMU_{\et}^{\ast}(Z;\Zln) \stackrel{i_{\ast}}{\longrightarrow} \hMU_{\et}^{\ast}(X;\Zln) \stackrel{j^{\ast}}{\longrightarrow} \hMU_{\et}^{\ast}(U;\Zln).$$
We call them the {\rm Gysin sequences} for $\hMU^{\ast}_{\et}(-)$, resp. $\hMU^{\ast}_{\et}(-;\Zln)$.
\end{cor}
\begin{proof}
This is proved in \cite{panin} 4.4.2.
\end{proof}

As another corollary of the last theorem we get
\begin{cor}
Let $k$ be a separably closed field and $\ell$ a prime number different from the characteristic of $k$. Profinite \'etale cobordism $\hMU^{2\ast}_{\et}(-)$, resp. $\hMU^{2\ast}_{\et}(-;\Zln)$, is an oriented cohomology theory on $\Sm/k$ in the sense of \cite{lm}. $\Box$
\end{cor}

Using Theorem \ref{PBF} we can identify $\hMU_{\et}^{2\ast}(\Pro^n)$ with the ring $\hMU^{\ast}[u]/(u^{n+1})$, identifying $c_1(\Oh(1))$ with $u$, and in the same way $$\hMU_{\et}^{\ast}(\Pro^n\times \Pro^n)=\hMU^{\ast}[u,v]/(u^{n+1},v^{n+1}).$$ As Quillen pointed out in \cite{quillen3}, Proposition 2.7, we can associate to the theory $\hMU_{\et}$, resp. $\hMU\Zln_{\et}$, a formal group law, i.e.  a power series $F(u,v) \in \hMU^{\ast}[[u,v]]$, $F(u,v)=\Sigma_{i,j\geq 0}a_{ij}u^iv^j$ satisfying some axioms with $a_{ij} \in \hMU^{2-2i-2j}$, to an oriented cohomology theory such that
$$c_1(L \otimes M)=F(c_1(L),c_1(M))$$ 
for any two line bundles $L$, $M$ over $X$.\\
Let $\Lee_{\ast}$ be the Lazard ring. There is Quillen's famous theorem that $\Lee_{\ast} \to \pi_{2\ast}(MU)$ is an isomorphism.\\
We have already determined the coefficients of $\hMU_{\et}$ over a separably closed field to be $\hMU^{\ast} \cong MU^{\ast}\otimes_{\Z} \Z_{\ell}$, resp. $(\hMU\Zln)^{\ast} \cong MU^{\ast}\otimes_{\Z} \Zln$, see Proposition \ref{coeffofhMUet}. But we have to check that the natural morphism $\Phi:\Lee^{\ast} \to \hMU^{\ast}_{\et}$ induced by the formal group law on $\hMU_{\et}$ agrees with this isomorphism.
\begin{prop}\label{FGLforhMU}
Let $k$ be a separably closed field. The map $\Phi_{\ell}:\Lee^{\ast} \otimes \Z_{\ell} \to \hMU^{2\ast}_{\et}(k)$ induced by the natural map of formal group laws is an isomorphism.\\
Similarly, the induced map $\Phi/\ell^{\nu} :\Lee^{\ast} \otimes \Zln \to \hMU^{2\ast}_{\et}(k;\Zln)$ is an isomorphism.
\end{prop}
\begin{proof}
We consider the induced map $\Phi_{\ell}:\Lee^{\ast} \otimes \Z_{\ell} \to \hMU^{2\ast}_{\et}(k)$. It is clear that this map $\Phi_{\ell}$ is a morphism in the category $\hat{\Lh}$ of \cite{dehon}, \S 2.1, that we already used in the proof of Proposition \ref{hMUofPn}. Furthermore, it is clear that this map induces an isomorphism $\Zl = \Lee^{\ast} \otimes \Z_{\ell}/f^1 \stackrel{\cong}{\to} \hMU^{2\ast}_{\et}(k)/f^1 \cong \Zl$ modulo the filtration $f^{\ast}$ of \cite{dehon}, \S 2.1. By \cite{dehon}, Corollaire 2.1.4, this implies that the $\Phi_{\ell}$ is an isomorphism.\\
The second assertion follows from the first one by the coefficient exact sequence.     
\end{proof}

\begin{remark}
The proof of the last proposition is similar to the proof of $\Phi:\Lee^{\ast} \cong MGL^{2\ast,\ast}(k)$ by Hopkins and Morel \cite{homo}. The idea is to construct a spectrum $MGL/(x_1,x_2,\ldots)$, where the $x_i$'s are the generators of $\Lee^{\ast}\cong \Z[x_1, x_2,\ldots]$. The crucial point for the proof of the surjectivity of the map $\Phi:\Lee^{\ast} \to MGL^{2\ast,\ast}(k)$ is to show that the map of spectra
$$MGL/(x_1, x_2, \ldots ) \to H\Z$$
is an isomorphism, where $H\Z$ denotes the motivic Eilenberg-MacLane spectrum. In the proof of Proposition \ref{FGLforhMU}, this corresponds to the isomorphism $\Zl = \Lee^{\ast} \otimes \Z_{\ell}/f^1 \stackrel{\cong}{\to} \hMU^{2\ast}_{\et}(k)/f^1 \cong H^{\ast}(\hEt k;\Zl)=\Zl$. \\
Then one shows by induction on $n$ that $\Lee^n \to MGL^{2n,n}(k)$ is surjective.\\
The injectivity of $\Phi$ can be shown as for $\Lee^{\ast} \to \Omega^{\ast}(k)$ in \cite{lm}.
\end{remark}

\begin{remark}\label{twisting}
The proof of Theorem \ref{transfer} is based on two results. On the one hand, we used Panin's discussion and results on oriented cohomology theories to deduce the missing structures from the existence of an orientation and the projective bundle bundle formula. On the other hand, we had to prove the projective bundle formula. For this purpose, we used the Atiyah-Hirzebruch spectral sequence to reduce the problem to the projective bundle formula for \'etale cohomology. \\
Although we worked over a separably closed base field, \'etale cohomology $\oplus_q H_{\et}^{\ast}(-;\mu_{\ell^{\nu}}^{\otimes q})$ is an oriented theory over any field $k$ with char $k \neq \ell$. The problem for our argument arises with the twist of the coefficients. Nevertheless, it seems to be very likely that the above proof also works over any field with characteristic prime to $\ell$. Maybe one needs to modify the definition of \'etale cobordism for this purpose.\\
One could start with a naive version of twisted coefficients for $\hMU_{\et}$ as follows. For any field $k$, let $\mu_{\ell^{\nu}}$ be the group of $\ell^{\nu}$th roots of unity in $k$. Since we can form a Moore spectrum $MG$ for every finite abelian group, we may consider the twisted Moore spectrum $M(\mu_{\ell^{\nu}}^{\otimes q})=M\Zln(q)$. Since $\mu_{\ell^{\nu}}$ is finite, we can consider the profinite completion $\hMU\Zln(q)$ of $\mathrm{MU}\wedge M\Zln(q)$. 
We could define profinite \'etale cobordism with twisted coefficients $\hMU_{\et}^{\ast}(-;\Zln(\ast))$ to be the profinite cohomology represented by the profinite spectrum $\hMU\Zln(\ast)$ applied to the functor $\hEt$ on schemes over $k$.\\ 
If $k$ is a separably closed field and if $X \in \Sm/k$, we have a canonical isomorphism 
$$\hMU_{\et}^{\ast}(X;\Zln(r)) \cong \hMU_{\et}^{\ast}(X;\Zln)\otimes \Zln(r)$$
for every $r$.\\ 
For, we have a canonical isomorphism for \'etale cohomology for every $r$
$$E_2^{p,q}=H_{\et}^p(X;MU^q \otimes \Zln(r)) \cong H_{\et}^p(X;MU^q \otimes \Zln) \otimes \Zln(r)=: \tilde{E}_2^{p,q}.$$ The two converging Atiyah-Hirzebruch spectral sequences imply the result.\\
But the geometry of $\hSHh$ does not seem to reflect this extra structure. It is more likely, that one has to find another way to define {\rm twisted \'etale cobordism} to handle this problem.
\end{remark}

\subsection{Etale Cobordism over finite fields}

In the case of a finite base field $k=\F_q$ of characteristic $p \neq \ell$ and with $q=p^r$ a power of $p$, we can also show that \'etale cobordism is an oriented cohomology theory on $\Sm/k$. This is a remarkable feature of this \'etale topological theory.\\
Again, we have to prove that $\hMU_{\et}^{\ast}(-)$ satisfies a projective bundle formula. We proceed in an analgue way as before. 

\begin{prop}\label{hEtkelltorsion} 
1. The profinite space $\hEt k$ is a space without $\ell$-torsion in $\hSh$.\\ 
2. The \'etale cobordism of the finite field $k$ is given by the following isomorphism
\[
\hMU^n_{\et}(k) = \left \{ \begin{array}{r@{\quad: \quad}l}
                                                  \hMU^n & n ~\mathrm{even} \\
                                                  \hMU^{n-1} & n ~\mathrm{odd}.
                                                 \end{array} \right. 
\]
We get an analogue result for $\hMU^{\ast}_{\et}(k;\Zl)$.
\end{prop}
\begin{proof}
1. The cohomology $H^i(\hEt k;\Z/\ell^m)=H^i(G_k;\Z/\ell^m)$ surjects onto $H^i(\hEt k;\Zl)$ for every $m$, since the Galois cohomology $H^i(G_k;\Z/\ell^m)$ of $k$ with coefficients in the trivial $G_k$-module $\Z/\ell^m$ equals $\Z/\ell^m$ for $i=0,1$ and vanishes for $i>1$ by \cite{serre}.\\
2. Let us denote the right hand side of the equation by $M^{\ast}$.
The isomorphism follows from the fact that $\hEt k$ has no $\ell$-torsion and hence 
$\hMU^{\ast}(\hEt k)/f^1 \to H^{\ast}(\hEt k;\Zl)$ is an isomorphism by Proposition2.1.8 of \cite{dehon}. Now $H^i(\hEt k;\Zl)$ equals $\Zl$ for $i=0,1$ and vanishes otherwise. Hence we have also an isomorphism $M^{\ast}/f^1 \cong H^{\ast}(\hEt k;\Zl)$ of free $MU^{\ast}$-algebras in $\hat{\Lh}$. By Corollaire 2.1.4 of \cite{dehon}, this implies the second assertion.
\end{proof}

\begin{prop}\label{hEtPnelltorsion}
The profinite space $\hEt \Pro_k^n$ is a space without $\ell$-torsion in $\hSh$ for every $n$.
\end{prop}
\begin{proof}
By the projective bundle formula for \'etale cohomology, the cohomology 
$H^{\ast}_{\et}(\Pro^n_k;\Z/\ell^m)$ of $\Pro^n_k$ is a free module over $H^{\ast}_{\et}(k;\Z/\ell^m)$ for every $m$. Since $H^{\ast}_{\et}(k;\Z/\ell^m)$ surjects onto $H^{\ast}_{\et}(k;\Zl)$ for every $m$, we deduce that the map $H^{\ast}(\hEt \Pro^n_k;\Z/\ell^m) \to H^{\ast}(\hEt \Pro^n_k;\Zl)$ is surjective for every $m$, too.
\end{proof}

\begin{prop}\label{hMUofPnfinite}
Let $k$ be as above and let $\Pro^n_k$ be the projective space of dimension $n$ over $k$.\\ 
1. The \'etale cobordism of $\Pro^n_k$ is given by the isomorphism
$$\hMU_{\et}^{\ast}(\Pro^n_k)\cong \hMU^{\ast}_{\et}(k)[u]/(u^{n+1}).$$
2. The \'etale cobordism with $\Zl$-coefficients of $\Pro^n_k$ is given by the isomorphism
$$\hMU_{\et}^{\ast}(\Pro^n_k;\Zl)\cong \hMU^{\ast}_{\et}(k;\Zl)[u]/(u^{n+1}).$$
\end{prop}
\begin{proof}
By Proposition \ref{hEtPnelltorsion} and Remark \ref{hEtkelltorsion} above, $\hEt \Pro^n_k$ and $\hEt k$  are profinite spaces without $\ell$-torsion. Now one can continue as above.
\end{proof}

\begin{prop}\label{weforPn}
For every $X \in \Sm/k$ the canonical map 
$\hEt (\Pro^n_k \times_k X) \to \hEt \Pro^n_k \times_{\hEt k} \hEt X$ is a weak equivalence in $\hSh$.
\end{prop}
\begin{proof}
There is a relative K\"unneth spectral sequence for the category $\hSh/Z$ of objects over $Z$
$$E_2=\Tor_{H^{\ast}(Z;\Zl)}(H^{\ast}(Y;\Zl), H^{\ast}(X;\Zl))$$ 
converging to $H^{\ast}(Y\times_Z X;\Zl)$. One can prove this result using the corresponding K\"unneth spectral sequence for simplicial finite sets and using the fact that colimits respect tensor products, see e.g. \cite{mcc} Theorem 8.34.\\
Since $H^{\ast}(\hEt \Pro^n_k;\Zl)$ is a free $H^{\ast}(\hEt k;\Zl)$-algebra, the Tor-term in the spectral sequence for $Y=\hEt \Pro^n_k$ and $Z=\hEt k$ is trivial. Hence in this special case, we deduce an isomorphism
$$H^{\ast}(\hEt \Pro^n_k;\Zl)\otimes_{H^{\ast}(\hEt k;\Zl)} H^{\ast}(\hEt X;\Zl) 
\stackrel{\cong}{\longrightarrow} H^{\ast}(\hEt \Pro^n_k \times_{\hEt k} \hEt X;\Zl)$$
for every $X \in \Sm/k$. Since the left hand side is isomorphic to $H^{\ast}(\hEt(\Pro^n_k \times_k X);\Zl)$ by the projective bundle formula for \'etale cohomology, we have proved the assertion of the proposition.
\end{proof}

The next step is a relative K\"unneth formula for profinite cobordism. 

\begin{prop}\label{relativeKunneth}
Let $Z \in \hSh$ and let $X$ and $Y$ be objects in the category of profinite spaces over $Z$. We suppose that $Y$ and $Z$ are finite dimensional spaces without $\ell$-torsion such that both 
$H^{\ast}(Y;\Zl)$ is a free $H^{\ast}(Z;\Zl)$-algebra and 
$\hMU^{\ast}(Y)$ is a free $\hMU^{\ast}(Z)$-algebra, respectively $\hMU^{\ast}(Y;\Zl)$ is a free $\hMU^{\ast}(Z;\Zl)$-algebra. Furthermore, we suppose that the $\Zl$-cohomology of $X$ is finite in all degrees.\\
Then the canonical maps
$$\hMU^{\ast}(Y)\otimes_{\hMU^{\ast}(Z)} \hMU^{\ast}(X)
\stackrel{\cong}{\longrightarrow} \hMU^{\ast}(Y \wedge_Z X)$$
and 
$$\hMU^{\ast}(Y;\Zl)\otimes_{\hMU^{\ast}(Z;\Zl)} \hMU^{\ast}(X;\Zl)
\stackrel{\cong}{\longrightarrow} \hMU^{\ast}(Y \wedge_Z X;\Zl)$$
are isomorphisms.
\end{prop}
\begin{proof}
I am very grateful to Francois-Xavier Dehon for an explanation of this point. We will prove the assertion following his argument.\\
We suppose first that both $X$ and $Y$ are finite dimensional without $\ell$-torsion. As described in the proof of Proposition 2.1.8 in \cite{dehon}, the spectral sequence for $MU$ degenerates at $E_2$. Since $X$, $Y$ and $Z$ have no $\ell$-torsion, the same is true for $X\wedge_Z Y$ as one shows inductively via the long exact sequence associated to the exact sequence of coefficients 
$0 \to \Zl \to \Z/\ell^{n+1} \to \Z/\ell^n \to 0$. Hence the classical relative K\"unneth formula for $H^{\ast}(X\wedge_Z Y; \hMU^{\ast})$, see e.g. \cite{smith} using the fact that $H^{\ast}(Y;\Zl)$ is a free $H^{\ast}(Z;\Zl)$-algebra, impies the relative K\"unneth formula for $\hMU$.\\
For the case that $X$ has $\ell$-torsion, we use induction on the skeletal filtration of $X$. Let $X^n$ be the $n$-th skeleton of $X$. If $\hMU^{\ast}(Y)\otimes_{\hMU^{\ast}(Z)}\hMU^{\ast}(W) \to \hMU^{\ast}(Y \wedge_Z W)$ is an isomorphism for $W=X^n$ and for $W=X^{n+1}/X^n$, then it is also an isomorphism for $W=X^{n+1}$. This follows from the long exact sequence associated to the cofibre sequence $X^n \to X^{n+1} \to X^n/X^{n+1}$ and the exactness of the tensor product with  $\hMU^{\ast}(Y)$ over $\hMU^{\ast}(Z)$. The exactness of $\hMU^{\ast}(Y) \otimes_{\hMU^{\ast}(Z)} -$ is due to the fact that $Y$ and $Z$ have no $\ell$-torsion and $\hMU^{\ast}(Y)$ is a free $\hMU^{\ast}(Z)$-algebra.\\
Finally, if $X$ has arbitrary dimension, we argue as above and take the limit over the skeletons of $X$. Then we deduce the result by taking into account that the cohomology groups of $X$ are finite dimensional $\Zl$-vector spaces by hypotheses.\\
Exactly the same argument holds for the case of $\hMU^{\ast}(-;\Zl)$ instead of $\hMU^{\ast}(-)$.
\end{proof}

\begin{cor}\label{relativeKunnethforPn}
For the projective $n$-space and every smooth scheme $X$ over $k$ the canonical maps for \'etale cobordism
$$\hMU^{\ast}_{\et}(\Pro^n_k)\otimes_{\hMU^{\ast}_{\et}(k)} \hMU^{\ast}_{\et}(X)
\stackrel{\cong}{\longrightarrow} \hMU^{\ast}_{\et}(\Pro^n_k \times_k X)$$
and 
$$\hMU^{\ast}_{\et}(\Pro^n_k;\Zl)\otimes_{\hMU^{\ast}_{\et}(k;\Zl)} \hMU^{\ast}_{\et}(X;\Zl)
\stackrel{\cong}{\longrightarrow} \hMU^{\ast}_{\et}(\Pro^n_k \times_k X;\Zl)$$
are isomorphisms.
\end{cor}
\begin{proof}
This follows immediately, since $\hMU^{\ast}_{\et}(\Pro^n_k)$ is a free $\hMU^{\ast}_{\et}(k)$-algebra by Proposition \ref{hMUofPnfinite}.
\end{proof}

Let $\Oh(1) \to \Pro^n_k$ be the canonical quotient line bundle. In fact, the isomorphism (\ref{hMUmodfinite}) is proved in \cite{dehon} by the following argument. Since $\hEt \Pro^n_k$ has no $\ell$-torsion the Atiyah-Hirzebruch spectral sequence implies that the orientation map 
$\hMU^{\ast}(\hEt \Pro^n_k) \to H^{\ast}(\hEt k;\Zl)$ factors through the isomorphism
$$\hMU^{\ast}(\hEt \Pro^n_k;\Zl) \longrightarrow \oplus_s H^s(\hEt \Pro^n_k;\Zl)\otimes \hMU^{\ast -s},$$
see the proof of Proposition 2.1.9 of \cite{dehon}. This implies that the element $\xi_H =c_1(\Oh(1)) \in  H^2(\hEt \Pro^n_k;\Zl)$ induces an element $\xi_{\hMU} \in \hMU^2(\hEt \Pro^n_k)$, 
and $\xi_{\hMU} \in \hMU^2(\hEt \Pro^n_k;\Zl)$ respectively, and that $\xi_{\hMU}$ is the image of $u$ under the isomorphism of Proposition \ref{hMUofPnfinite} above.\\
Hence we can reformulate the above assertion as
\begin{equation}\label{PFBfinite}
\hMU_{\et}^{\ast}(\Pro^n_k)\cong \oplus_{i=0}^n \hMU^{\ast}_{\et}(k)\xi^i
\end{equation}
and similarly with $\Zl$-coefficients.\\
The element $\xi_{\hMU} \in \hMU^2(\hEt \Pro^n_k)$ is given by a morphism in $\hSHh$
$$\xi_{\hMU}: \hat{\Sigma}^{\infty}(\hEt \Pro^n_k) \longrightarrow \hMU \wedge S^2.$$

Now let $E \to X$ be a vector bundle over $X$ in $\Sm/k$. Let $\Pro(E)$ be the associated projective bundle and let $\Oh(1)$ be the canonical quotient line bundle. This line bundle determines a morphism 
$\Pro(E) \to \Pro^N_k$ for some sufficiently large $N$. Together with the morphism $\xi_{\hMU}$ we get an element $\xi \in \hMU^2_{\et}(\Pro(E))$.\\
Now one can continue as above to show that cobordism with $\Zl$-coefficients over a finite field is an oriented cohomology theory.

\section{Algebraic versus Profinite Etale Cobordism}

Let $\ell$ be a fixed prime number different from the characteristic of the base field $k$. Using the results of the previous sections, we compare \'etale cobordism with the algebraic cobordism theories of Levine/Morel and Voevodsky over a separably closed field.\\ 
Furthermore, we construct a map from the $MGL$-theory to profinite \'etale cobordism via the results of the previous section on the stable \'etale realization functor. The construction is not as straightforward as one would like, since we do not know if $\hEt MGL$ is isomorphic to $\hMU$ in $\hSHh$. Nevertheless, the stable \'etale realization of Theorem \ref{P1real} is the main technical advantage of the use of $\hEt$ instead of $\Et$ and is a necessary ingredient for the construction of the map from algebraic to profinite \'etale cobordism. The other key point is an \'etale Thom class for vector bundles. This Thom class yields the map $\hEt MGL \to \hMU$.\\
At the end, we conjecture that our map from algebraic cobordism to profinite \'etale cobordism is an isomorphism for $\Zln$-coefficients after inverting a Bott element. We explain a strategy to prove this conjecture and give arguments for the truth of the conjecture.  

\subsection{Comparison with $\Omega^{\ast}$}

We consider the algebraic cobordism theories $\Omega^{\ast}(-)$ and $\Omega^{\ast}(-;\Zln)$. As a corollary of Theorem \ref{transfer} we get

\begin{theorem}\label{OmegatohMU}
By universality of $\Omega^{\ast}$, the structure of an oriented cohomology theory on $\hMU_{\et}$ and $\hMU\Zln_{\et}$ yields canonical morphisms
$$\theta:\Omega^{\ast}(X) \to \hMU_{\et}^{2\ast}(X)$$ and  
$$\theta:\Omega^{\ast}(X;\Zln) \to \hMU_{\et}^{2\ast}(X;\Zln)$$  
defined by sending a generator $[f:Y\to X] \in \Omega^{\ast}(X)$, $f:Y \to X$ a projective morphism between smooth schemes, to the element $f_{\ast}(1_Y) \in \hMU_{\et}^{2\ast}(X)$; similarly with $\Zln$-coefficients.
\end{theorem}

If $k$ is a field of characteristic zero, Theorem 4 of \cite{lm} states that $\Omega^{\ast}(k) \cong \Lee^{\ast}$. But it is conjectured that this isomorphism holds for every field, see Conjecture 2 of \cite{lm}.  
We refer to this conjecture as 
\begin{assump}\label{assumption}
For the field $k$, the canonical morphism $\Lee^{\ast} \stackrel{\cong}{\to} \Omega^{\ast}(k)$ is an isomorphism.
\end{assump}

If we suppose that this isomorphism holds for $k$, we may conclude 
\begin{prop}\label{compofcoeff}
For every separably closed field $k$, the morphisms 
$$\theta:\Omega^{\ast}(k;\Zln) \to \hMU_{\et}^{2\ast}(k;\Zln)$$ and 
$$\theta \otimes \Z_{\ell}:\Omega^{\ast}(k) \otimes_{\Z} \Z_{\ell} \to \hMU_{\et}^{2\ast}(k)$$
are surjective.\\
If we suppose in addition that Assumption \ref{assumption} is true for $k$,  
then $\theta$ becomes an isomorphism on coefficient rings 
$$\theta:\Omega^{\ast}(k;\Zln) \cong \hMU_{\et}^{2\ast}(k;\Zln)$$ and
$$\theta:\Omega^{\ast}(k) \otimes_{\Z} \Z_{\ell} \cong \hMU_{\et}^{2\ast}(k).$$
\end{prop}
\begin{proof}
For the surjectivity, consider the canonical map $\Phi:\Lee^{\ast} \otimes \Zln \to \Omega^{\ast}(k;\Zln)$. Levine and Morel have shown that it is injective for every field, see \cite{lm}, Corollary 12.3. When we compose this map with $\theta: \Omega^{\ast}(k;\Zln) \to \hMU_{\et}^{2\ast}(k;\Zln)$, we get the canonical map 
$$\Lee^{\ast} \otimes \Zln \stackrel{\Phi}{\to} \Omega^{\ast}(k;\Zln) \stackrel{\theta}{\to} \hMU_{\et}^{2\ast}(k;\Zln).$$
Since this map is unique, it must be the canonical isomorphism. Hence $\theta: \Omega^{\ast}(k;\Zln) \to \hMU_{\et}^{2\ast}(k;\Zln)$ is surjective.\\
A similar argument works for $\theta: \Omega^{\ast}(k) \otimes \Z_{\ell} \to \hMU_{\et}^{2\ast}(k)$.\\
Note that for the surjectivity assertion we did not need to suppose that Assumption \ref{assumption} is true.\\
The last part follows from the assumption and  Proposition \ref{coeffofhMUet}.\\
\end{proof}

\subsection{Comparison with $MGL^{\ast,\ast}$}

Now we dedicate our attention to the comparison with the algebraic cobordism theory represented by the $MGL$-spectrum in the stable motivic $\A^1$-homotopy category. The technical advantage of $\hEt$  compared to $\Et$ now becomes apparent since we could construct a stable realization of $\Pro^1$-spectra in the previous section. We will use this functor to define a map from $MGL$-theory to profinite \'etale cobordism. First we define Thom classes which will enable us to define canonical maps of oriented cohomology theories. We will do this in essentially the same way in which one shows that $MGL$ is the universal oriented $\Pro^1$-spectrum.\\
Let $k$ be a separably closed field. Let $V$ be a vector bundle of rank $d$ over $X$ in $\Sm/k$. We recall that the Thom space $\Thom (V) \in \LU$ of $V$ is defined to be the quotient $$\Thom(V) =\Thom(V/X)=V/(V-i(X))$$ where $i:X \to V$ denotes the zero section of $V$. We reformulate a lemma from $\A^1$-homotopy theory.
\begin{prop}
Let $V$ be a vector bundle over $X$ and $\Pro(V) \to \Pro(V \oplus \Oh)$ be the closed embedding at infinity. Then the canonical morphism of pointed sheaves: $\Pro(V\oplus \Oh)/\Pro(V) \to \Thom(V)$ induces a weak equivalence in $\hShp$ via $\hEt$.  
\end{prop} 
\begin{proof}
This is the same proof as for Proposition 3.2.17 of \cite{mv} where we use the fact that $\hEt$ preserves $\A^1$-weak equivalences between smooth schemes, see Theorem \ref{LhEt}, and commutes with quotients, i.e $\hEt (X/A) = \hEt X/\hEt A$. 
\end{proof}

\begin{cor}
Let $T=\A^1/(\A^1-0)$ in $\LU$. Then we have an isomorphism in $\hHhp$ of pointed profinite spaces $\hEt(T) \cong S^2$.
\end{cor}
\begin{proof}
This follows from the previous proposition using the fact that $T=\Thom(\Oh)$ and $\hEt \Pro^1\cong S^2$ in $\hHh$ by Example \ref{p1}. 
\end{proof} 
  
Now we define the Thom class of $V$ in $\hMU_{\et}^{2d}(\Thom(V))$. From the isomorphism $\Pro(V\oplus \Oh)/\Pro(V) \cong \Thom(V)$ of the previous proposition we deduce an exact sequence induced by the cofiber sequence
$$\hMU_{\et}^{\ast}(\Thom(V)) \to \hMU_{\et}^{\ast}(\Pro(V\oplus \Oh)) \to \hMU_{\et}^{\ast}(\Pro(V)).$$
Using the projective bundle formula of Theorem \ref{PBF} this sequence is isomorphic to the exact sequence 
$$\hMU_{\et}^{\ast}(\Thom(V)) \to \hMU^{\ast}[1,u,\ldots,u^d] \to \hMU^{\ast}[1,u,\ldots,u^{d-1}].$$  
The element $u^d -c_1(V)u^{d-1}+\ldots+(-1)^dc_d(V) \in \hMU^{\ast}[1,u,\ldots,u^d]$ is sent to $u^d -c_1(V)u^{d-1}+\ldots+(-1)^dc_d(V) \in \hMU^{\ast}[1,u,\ldots,u^{d-1}]$ which is $0$ by the definition of Chern classes via the projective bundle formula. By exactness, there is an element $\thomet(V) \in \hMU_{\et}^{2d}(\Thom(V))$ that is sent to the above element $u^d+\ldots+(-1)^d c_d(V) \in \hMU^{\ast}[1,u,\ldots,u^d]$ which is homogeneous of degree $2d$ since $u$ has degree $2$. We call $\thomet(V)$ an {\em \'etale Thom class of $V$}. Hence the Thom class of $V$ is a morphism 
$$\thomet(V):\hat{\Sigma}^{\infty}(\hEt \Thom(V)) \to \hMU \wedge S^{2d}$$ in $\hSHh$.\\
In the same way, we define an \'etale Thom class for $\hMU\Zln_{\et}$.\\
We apply this argument to the tautological $n$-bundle $\gamma_n$ over the infinite Grassmannian. Since $MGL_n =\Thom(\gamma_n)$ we get via $\hEt$ a map of profinite spaces
$$\hEt(MGL_n) \to \hMU_{2n}.$$
In the proof of Theorem \ref{P1real} we have seen that $\hEt$ yields a map of spectra on the level of homotopy categories. Hence the above map yields a morphism in $\hSHh_2$   
\begin{equation}\label{hEtMGLtohMU}
\phi:\LhEt(MGL) \longrightarrow \hMU
\end{equation}
where $\LhEt$ denotes the total left derived functor of $\hEt$ on $\Pro^1$-spectra.\\
Since smooth schemes $X \in \Sm/k$ are cofibrant objects, $\LhEt$ and $\hEt$ agree on $\Sm/k$ by Theorem \ref{LhEt}. Hence this morphism in $\hSHh_2$ yields a natural morphism of cohomology theories
$$(\LhEt MGL)^{\ast}(\hEt X) \to \hMU_{\et}^{\ast}(X)$$
for every scheme $X$ in $\Sm/k$. It is clear that this construction is natural in $X$.\\ 
Recall the definition of algebraic cobordism for a scheme $X$ in $\Sm/k$
$$MGL^{p,q}(X):=\Hom_{\SHhP}(\Sigma^{\infty}_{\Pro^1}(X),MGL[p-2q]\wedge (\Pro^1)^{\wedge q}),$$
where $E[n]:=E\wedge (S^1_s)^{\wedge n}$, $S^1_s$ denoting the simplicial circle. We have seen in Section 6 that $S^1_s \wedge -$ commutes with $\hEt$. Hence the elements in degree $p,q$ are sent to elements in degree $p$ via $\phi$. We summarize this discussion in the following

\begin{theorem}\label{MGLcomparison}
The obvious induced map of Theorem \ref{P1real} yields a natural map for every $X$ in $\Sm/k$
\begin{equation}\label{MGLtohMU}
\phi: MGL^{\ast,\ast}(X) \to \hMU_{\et}^{\ast}(X),
\end{equation}
which is $\phi: MGL^{p,q}(X) \to \hMU_{\et}^{p}(X)$ on the $p,q$-level.\\
In the same way we get a map for $MGL$ with $\Zln$-coefficients
\begin{equation}\label{MGLZlntohMU}
\phi: MGL^{\ast,\ast}(X;\Zln) \to \hMU_{\et}^{\ast}(X;\Zln).
\end{equation}
\end{theorem}

\begin{remark}\label{MGLtoHZln}
In the same way as above we get a morphism of profinite spectra $\hEt MGL \to H\Zln$ which, for every $X$ in $\Sm/k$, gives rise to maps from algebraic cobordism to \'etale cohomology
$$MGL^{\ast,\ast}(X) \to H_{\et}^{\ast}(X;\Zln)$$ and 
$$MGL^{\ast,\ast}(X;\Zln) \to H_{\et}^{\ast}(X;\Zln).$$
\end{remark}

Recall that a morphism of oriented cohomology theories is a natural transformation that commutes with Chern classes and transfer maps. In fact, Theorem 4.1.4 of \cite{panin} states that we have to check only the compatibility with one of these structures. It is clear from the construction that we get the following 
\begin{theorem}\label{orientedmap}
The natural transformations $\phi: MGL^{\ast,\ast}(-) \to \hMU_{\et}^{\ast}(-)$ and $\phi: MGL^{\ast,\ast}(-;\Zln) \to \hMU_{\et}^{\ast}(-;\Zln)$ are morphisms of oriented cohomology theories on $\Sm/k$. 
\end{theorem}
\begin{proof}
We have to check that $\phi$ respects Chern classes in both theories. Then Theorem 4.1.4 of \cite{panin} implies that the transfer maps in both theories agree. But the construction of Chern classes in $\hMU_{\et}$ is the exact analog of the construction of Chern classes in $MGL$ only involving the functor $\hEt$ at the appropriate places, i.e. the orientation $x_{\hMU}$ of $\hMU_{\et}$ is just $\hEt (x_{MGL})$, i.e. the image of the orientation $x_{MGL}$ of $MGL$ under $\hEt$. Hence we have $\phi(c^{MGL}_1(L))=c^{\hMU}_1(L)$ in $\hMU_{\et}^2(X)$. Then Theorem 4.1.4 of \cite{panin} implies that $\phi$ respects the oriented structures and pushforwards, i.e. $\phi$ is a morphism of oriented cohomology theories. \\
The same argument holds for $\Zln$-coefficients. 
\end{proof}

\begin{remark}
We could have deduced from this theorem that $\hMU_{\et}$ carries in fact the formal group law induced by $\Lee^{\ast} \to \Lee^{\ast} \otimes \Z_{\ell}$. Since Hopkins and Morel proved in \cite{homo} that $\Phi: \Lee^{\ast} \to MGL^{2\ast,\ast}(k)$ is an isomorphism, the above theorem implies that the formal group law on $\hMU_{\et}$ is the one induced by $MGL^{2\ast,\ast}$.
\end{remark}

As a corollary we get the following
\begin{theorem}\label{triangle}
Let $k$ be a separably closed field. The canonical map of oriented cohomology theories $\theta:\Omega^{\ast}(-) \to \hMU_{\et}^{2\ast}(-)$ factors through $MGL^{2\ast,\ast}(-)$, i.e. the following diagram
$$\begin{array}{ccc}
\Omega^{\ast}(X) & \stackrel{\theta_{MGL}}{\longrightarrow} &  MGL^{2\ast,\ast}(X)\\
\theta_{\hMU} \searrow &  & \swarrow \phi \\
 & \hMU_{\et}^{2\ast}(X) & 
\end{array}$$
is commutative for every $X$ in $\Sm/k$.\\
Furthermore, the diagram
$$\begin{array}{ccc}
\Omega^{\ast}(X;\Zln) & \stackrel{\theta_{MGL}}{\longrightarrow} &  MGL^{2\ast,\ast}(X;\Zln)\\
\theta_{\hMU} \searrow &  & \swarrow \phi \\
 & \hMU_{\et}^{2\ast}(X;\Zln) & 
\end{array}$$
is commutative for every $X$ in $\Sm/k$.
\end{theorem}
\begin{proof}
By universality of $\Omega^{\ast}$, $\theta_{\hMU}$ is the unique morphism of oriented cohomology theories on $\Sm/k$ from $\Omega^{\ast}(-)$ to $\hMU_{\et}^{2\ast}(-)$. Since the composition of $\theta_{MGL}:\Omega^{\ast}(-) \to MGL^{2\ast,\ast}(-)$ with $\phi$ is also a morphism of oriented cohomology theories by Theorem \ref{orientedmap}, the triangle commutes.\\
The same argument holds for $\Zln$-coefficients.
\end{proof}

\begin{remark}\label{MGLet}
1. Hopkins and Morel have proved that $\Omega^{\ast}(X) \to MGL^{2\ast,\ast}(X)$ is surjective for every $X$ in $\Sm/k$ over every field of characteristic zero in \cite{homo}.\\
2. Let $MGL_{\et}$ be the $\Pro^1$-Thom-spectrum of algebraic cobordism in $\SHhP_{\et}$, i.e. instead of using the Nisnevich topology, we consider $MGL$ and smooth schemes as sheaves for the \'etale topology.\\
The map $\phi$ actually factors through $MGL_{\et}^{\ast,\ast}(X)$, since $\hEt$ factors through the \'etale version of $\A^1$-homotopy theory, see Remark \ref{factorization}.\\ 
3. Suppose that one has constructed an Atiyah-Hirzebruch spectral sequence from \'etale cohomology with coefficients in $\Zln \otimes MU^{\ast}$ to $MGL_{\et}(-;\Zln)$. It is clear from the construction of $\phi$ and the isomorphism $H_{\et}^{\ast}(X;\Zln \otimes MU^{\ast}) \stackrel{\cong}{\to} H^{\ast}(\hEt X;\Zln \otimes MU^{\ast})$ that this spectral sequence would imply an isomorphism 
$$\phi_{\et}:MGL_{\et}^{\ast,\ast}(X;\Zln) \stackrel{\cong}{\longrightarrow} \hMU_{\et}^{\ast}(X;\Zln)$$
whenever both spectral sequences converge.\\
Such an Atiyah-Hirzebruch spectral sequence for motivic cohomology and Nisnevich $MGL$ has been constructed by Hopkins and Morel \cite{homo}. We will discuss an application of it in the following section.\\
4. We have seen that this definition of profinite \'etale cobordism is not good enough to reflect questions of the $q$-twist of the $MGL$-theory, see also Remark \ref{twisting}. This point should be clarified in future work. 
\end{remark}

\subsection{A conjecture}

At the end of this section we would like to explain and give arguments for a conjecture on the relation of  algebraic and profinite \'etale cobordism after inverting a Bott element.\\
First we recall a theorem of Levine for cohomology. Let $H^p(X;\Z/n(q))$ denote the motivic cohomology of a smooth scheme $X$ over a field $k$. For $\Spec k$ there is an isomorphism 
$H^0(\Spec k;\Z/n(1)) \cong \mu_n(k)$ with the group of $n$-th roots of unity in $k$. We suppose that $k$ contains an $n$-th root of unity $\zeta$, we have a corresponding {\em motivic Bott element} $\beta_n \in  
H^0(\Spec k;\Z/n(1))$. \\
Levine has shown in \cite{bott} that motivic $\Z/n$-cohomology of a smooth scheme over $k$ agrees with \'etale $\Z/n$-cohomology after inverting the Bott element. More precisely, we form the bigraded ring 
$$H^{\ast}(X;\Z/n(\ast))[\beta_n^{-1}]:= (\oplus_{p,q}H^p(X;\Z/n(q)))[\beta_n^{-1}].$$ 
There is a cycle class map $cl^{\ast,\ast}:H^{\ast}(X;\Z/n(\ast)) \to H^{\ast}_{\et}(X;\mu_n^{\otimes \ast})$ with $cl^{0,1}(\beta_n)=\zeta$. Since $k$ contains $\zeta$, the cup product with $\zeta$ is an isomorphism on \'etale cohomology. Hence we get an induced map
\begin{equation}\label{motbottmap}
cl^{\ast,\ast}:H^{\ast}(X;\Z/n(\ast))[\beta_n^{-1}] \to H^{\ast}_{\et}(X;\mu_n^{\otimes \ast}).
\end{equation} 
This map is an isomorphism for every smooth scheme $X$ over $k$ if $n$ is odd or if $4|n$ or if $k$ contains a square root of $-1$, cf. Theorem 1.1 of \cite{bott}.\\
Furthermore, Hopkins and Morel \cite{homo} have constructed an Atiyah-Hirzebruch spectral sequence for algebraic cobordism
\begin{equation}\label{MGLAHSS}
E_2^{\ast,\ast,\ast}=H^{\ast;\ast}(X,MU^{\ast}) \Longrightarrow MGL^{\ast,\ast}(X),
\end{equation}
where $H^{\ast,\ast}(X;MU^{\ast})$ denotes motivic cohomology with coefficients in the ring $MU^{\ast}$. We still suppose that $k$ contains an $n$-th root of unity. We consider the spectral sequence for $\Spec k$. Since $E_2^{p,1,q}(\Spec k)$ is concentrated in degrees $p=0$ and $p=1$ we deduce an isomorphism $MGL^{0,1}(k)\cong k^{\times}$. For $\Z/n$-coefficients, the exact sequence for coefficients implies that we get an isomorphism $$MGL^{0,1}(k;\Z/n)\cong \mu_n(k)$$ and, via the spectral sequence, the motivic Bott element defined above, more precisely $\beta_n\cdot 1_{MU/n}$, is sent to an induced Bott element $\beta_n \in MGL^{0,1}(k;\Z/n)$. \\
Let us now suppose that $k$ is separably closed, $n=\ell^{\nu}$, and that the characteristic of $k$ is not equal to $\ell$. In particular, $k$ contains an $\ell^{\nu}$-th root of unity $\zeta$. It  defines an element $\zeta \in H_{\et}^0(\Spec k;\mu_{\ell^{\nu}})$. The element $\zeta \cdot 1_{MU/\ell^{\nu}}$ induces via the Atiyah-Hirzebruch spectral sequence an element $\zeta_{\hMU} \in \hMU_{\et}^0(\Spec k;\Zln)$ or equivalently a unit element in $\hMU_{\et}^0(\Spec k;\Zln)\cong \Zln$. The multiplication with $\zeta$ yields an isomorphism of spectral sequences and hence multiplication with $\zeta_{\hMU}$ is an isomorphism on \'etale cobordism. Furthermore, the map $\phi$ of (\ref{MGLtohMU}) maps the element $\beta_{\ell^{\nu}}$ to the element $\zeta_{\hMU}$. This implies that $\phi$ also yields a map $\phi: MGL^{\ast,\ast}(X;\Zln)[\beta ^{-1}] \to \hMU_{\et}^{\ast}(X;\Zln)$.

\begin{conjecture}\label{Bottconjecture}
Let $X$ be a smooth scheme of finite type over a separably closed field $k$ of characteristic different from $\ell$. Suppose that $\ell$ is odd or that $\ell^{\nu}\geq 4$. Let $\beta \in MGL^{0,1}(k;\Zln)$ be the Bott element defined above. We conjecture that the morphism
$$\phi: MGL^{\ast,\ast}(X;\Zln)[\beta ^{-1}] \to \hMU_{\et}^{\ast}(X;\Zln)$$
is an isomorphism.
\end{conjecture}

In fact, one should conjecture that $\phi$ can be defined for every field and yields an isomorphism over every field that satisfies the hypothesis. For our construction of $\phi$ in the previous section, we have used that $k$ was separably closed. But it seems to be very likely, that $\phi$ does not depend on that fact. On the other hand, we have argued in the previous section, see Remarks \ref{twisting} and \ref{MGLet}, that $\hMU_{\et}$ may not be the best theory over a non-separably closed field. We intend to ameliorate the theory for this purpose.\\ 
Nevertheless, the proof of the conjecture over a separably closed field would already be a big step towards the proof of a more general conjecture. Thomason \cite{thomason} and Levine \cite{bott} proved that algebraic K-theory, respectively motivic cohomology, satisfies \'etale descent after inverting the Bott element. Then they showed that the statement holds over an algebraically closed field and reduced the general case via descent to this special case.\\
Hence for algebraic cobordism one could try the same, i.e. one could try to prove that $MGL^{\ast,\ast}(-;\Zln)[\beta ^{-1}]$ satisfies \'etale descent and then reduce to the special case of a separably closed field. 

\section{Calculations and Applications}

Most of the following calculations of profinite \'etale cobordism are due to the Atiyah-Hirzebruch spectral sequence. We consider some examples of schemes whose \'etale cohomology groups are known and deduce from this the \'etale cobordism groups with $\Zln$-coefficients. But we start with the naturally arrising question how the absolute Galois group of a field acts on \'etale cobordism.
\subsection{The Galois action on \'etale cobordism}
The comparison with $\Omega^{\ast}(-;\Zln)$ enables us to study the action of the Galois group of a field on $\hMU_{\et}^{\ast}$. \\
Let $k$ be a field of characteristic $p \neq \ell$ and let $\ok$ be a separable closure of $k$. If $X$ is a scheme over $k$, there is a natural action on $\Xok=X \otimes_k \ok$ of the Galois group $G_k:=\Gal(\ok/k)$ of $k$. The definition of $\hEt$ immediately implies an action of $G_k$ on $\hEt \Xok$. 
\begin{defn}\label{galois}
Let $X$ be defined over a field $k$. We define the action of the Galois group $G_k:=\Gal(\ok/k)$ on $\hMU_{\et}^{\ast}(\Xok)$, resp. $\hMU_{\et}^{\ast}(\Xok;\Zln)$, to be the one induced by the action of $G_k$ on $\hEt \Xok$.
\end{defn}
In order to study the Galois action on $\hMU^{\ast}_{\et}(\ok;\Zln)$, we use Proposition \ref{compofcoeff}. We deduce the following result from \cite{lm}, \S 3.2. Let $H_{m,n} \hookrightarrow \Pro^n \times \Pro^m$ denote a smooth closed subscheme defined by the vanishing of a transverse section of the line bundle $p_1^{\ast}\Oh(1) \otimes p_2^{\ast}\Oh(1)$. The smooth projective $\ok$-schemes $H_{m,n}$ are known as Milnor hypersurfaces.\\ 
From the definition of the map $\Omega^{\ast}(\ok) \to \hMU^{2\ast}_{\et}(\ok;\Zln)$, sending the class $[Y \to \Spec \ok]$ to $\pi_{\ast}(1_Y)$, we deduce the following
\begin{prop}
Let $\ok$ be a separably closed field.
Then $\hMU_{\et}^{\ast}(\ok)$, resp. $\hMU_{\et}^{\ast}(\ok;\Zln)$, is generated by the classes of projective $n$-spaces and the classes of smooth Milnor hypersurfaces, i.e by the elements $\pi_{\ast}(1_{\Pro^n_{\ok}})$ for the projections $\pi:\Pro^n_{\ok} \to \Spec \ok$ and the elements $\pi_{\ast}(1_{H_{m,n}})$ with projection $\pi:H_{m,n} \to \Spec \ok$.  
\end{prop}
\begin{proof} 
This follows from Proposition \ref{compofcoeff} and the fact that $\Omega^{\ast}(\ok)$ is generated by the classes $[\Pro^n_{\ok}]$ and $[H_{m,n}]$ for all $m$, $n$, see \cite{lm}, Remark 3.7. 
\end{proof}

Hence, in order to determine the action of the absolute Galois group on the group $\hMU_{\et}(\ok)$, resp. $\hMU_{\et}(\ok;\Zln)$, it suffices to know the action of the Galois group on the projective space $\pi:\Pro^n_{\ok} \to \ok$ and the hypersurface $H_{m,n} \to \ok$.\\ 
By \cite{milne} VI, 5.6, the \'etale cohomology of $Y$, for $Y=\Pro^n_{\ok}$ or $Y=H_{m,n}$, is given by $H_{\et}^i(Y;\Zl) \cong \Zl(-i/2)$ if $i$ is even, $0 \leq i \leq 2n$, resp. $0\leq i \leq 2(m+n-1)$, and vanishes otherwise, where the twist is given by the cyclotomic character. Consider $\sigma \in G_k$. Since $k_s$ is separably closed, this implies that the induced morphism $\sigma: Y \to Y^{\sigma}$ induces an isomorphism on \'etale cohomology with $\Zl$-coefficients, where $Y^{\sigma}$ is the twisted $\ok$-scheme with structure morphism $\sigma \circ \pi$. Hence $\sigma$ induces an isomorphism in $\hHh$: 
$$\sigma:\hEt Y^n_{\ok} \stackrel{\cong}{\to} \hEt Y^{\sigma}$$
for $Y=\Pro^n_{\ok}$ or $Y=H_{m,n}$.
This implies the following 
\begin{theorem}\label{Galoisaction}
The action of $G_k$ on $\hMU^{\ast}_{\et}(\ok;\Zln)$ is trivial.
\end{theorem}
\begin{remark}
It is clear from the construction of the Atiyah-Hirzebruch spectral sequence in Theorem \ref{ahss} that the corresponding spectral sequence for $\hMU\Zln_{\et}$ is equivariant under the Galois action. This can be seen as follows. By definition, the simplicial structure of $\hEt X$ is given by the connected components $\pi(U_i)$ of the \'etale hypercoverings $\ldots \to U_1 \to U_0 \to X$. Hence the action of $G_k$ on $\hEt X$ is in fact an action on every simplicial dimension. In particular, this implies that the action of $G_k$ preserves the skeletal filtration on $\hEt \Xok$ and induces an action on every exact couple. Since both actions on $H_{\et}^{\ast}(X;\Zln)$ and on $\hMU_{\et}^{\ast}(\Xok;\Zln)$ are induced by the one on $\hEt \Xok$, the assertion follows.
\end{remark}
 \subsection{Etale Cobordism for local fields}
Next, we calculate the profinite \'etale cobordism groups of a local field. 
\begin{theorem}\label{hMUofQp}
Let $k$ be a local field, i.e. either a finite extension of the field $\Q_p$ or a finite extension of the field of formal power series $\F((t)))$ over a finite field of characteristic $p$. We assume $p\neq \ell$ and we denote by $q=p^f$ the number of elements in the residue class field of $k$. Let $\nu_0 \leq \nu$ be such that $\ell^{\nu_0}=(q-1,\ell^{\nu})$ is the greatest common divisor of $q-1$ and $\ell^{\nu}$. For all $n$, the profinite \'etale cobordism groups with $\Zln$-coefficients of $k$ are given by 
\[
\hMU^n_{\et}(k;\Zln) = \left \{ \begin{array}{r@{\quad: \quad}l}
                                          (\Z/\ell^{\nu_0}\otimes MU^{n-2})\oplus(\Zln\otimes MU^n) & n ~\mathrm{even} \\
                                          (\Zln \oplus \Z/\ell^{\nu_0})\otimes MU^{n-1} & n ~\mathrm{odd}.
                                          \end{array} \right. 
\] 
\end{theorem} 
\begin{proof}
We prove this theorem again via the Atiyah-Hirzebruch spectral sequence for \'etale cobordism using the identification of Galois and \'etale cohomology. Since $\Zln \otimes MU^t$ is a finitely generated free $\Zln$-module with a trivial $G_k$-operation, we may identify the Galois cohomology groups $H^i(k;\Zln \otimes MU^t)$ with $H^i(k;\Zln)\otimes MU^t$. Hence we only have to consider the cohomology groups $H^i(k;\Zln)$. These groups vanish for $i>2$, see e.g. \cite{nsw}, Chapter VII \S 1.\\
Since the $G_k$-operation on $\Zln$ is trivial, we have 
$$H^0(k;\Zln)=\Zln.$$ 
By Theorem 7.2.9 of Chapter VII, \S 2 in \cite{nsw}, there is a duality for the finite $G_k$-modules $A:=\Zln$ and $A':=\Hom(A,\mu_{\ell^{\nu}})$ between the finite groups $H^i(k;A)$ and $H^{2-i}(k;A')$ for $0\leq i \leq 2$ given by the cup-product pairing
\begin{equation}\label{Galois-duality}
H^i(k;A) \times H^{2-i}(k;A') \longrightarrow H^2(k;\mu_{\ell^{\nu}})=\Zln.\end{equation}
For $i=1$, by Theorem 7.1.8 in \cite{nsw}, Chapter VII \S 1, we know $H^1(k;\mu_{\ell^{\nu}})=k^{\times}/k^{\times \ell^{\nu}}$. This yields by (\ref{Galois-duality}) an isomorphism $H^1(k;\Zln) \cong \Hom(H^1(k;\mu_{\ell^{\nu}}),\Zln)=\Hom(k^{\times}/k^{\times \ell^{\nu}}, \Zln)$. By local class field theory, \cite{neukirch}, Chapter II \S 5, Proposition 5.7, using the fact that $p\neq \ell$, there is an isomorphism $k^{\times}/k^{\times \ell^{\nu}} \cong \Zln \oplus \Z/\ell^{\nu_0}$, where $\ell^{\nu_0}$ with $\nu_0 \leq \nu$ is the greatest common divisor $(q-1, \ell^{\nu})$. Hence we get an isomorphism 
$$H^1(k;\Zln)\cong \Zln \oplus \Z/\ell^{\nu_0}.$$
For $i=2$, by Theorem 7.2.9 in \cite{nsw}, Chapter VII \S 1, we know $H^2(k;\mu_{\ell^{\nu}})=\Zln$. Pairing (\ref{Galois-duality}) yields an isomorphism $H^2(k;\Zln)\cong \Hom(H^0(k;\mu_{\ell^{\nu}}), \Zln)$. But $H^0(k;\mu_{\ell^{\nu}})$ is equal to the group $\mu(k)$of $\ell^{\nu}$th roots of unity in $k$, which is isomorphic to $\Z/\ell^{\nu_0}$, with $\nu_0 = (q-1,\ell^{\nu})$, by \cite{neukirch}, Chapter II \S 5, Proposition 5.3. Hence we get an isomorphism
$$H^2(k;\Zln)\cong \Z/\ell^{\nu_0}.$$
Hence the associated Atiyah-Hirzebruch spectral sequence from the Galois cohomology of $k$ to $\hMU_{\et}^{\ast}(k;\Zln)$ is still easy to describe. It collapses at the $E_2$-stage since $E_2^{s,t}\neq 0$ if and only if $s=0,1$ or $2$ and $t$ even, and we know the occuring terms by the above arguments.\\ 
If $n$ is odd, $\hMU_{\et}^n(k;\Zln)$ is isomorphic to $E_2^{1,n-1}=H^1(G_k;(\hMU\Zln)^{n-1})\cong (\Zln \oplus \Z/\ell^{\nu_0})\otimes MU^{n-1}$.\\ 
If $n$ is even, there are two terms occurring: $E_2^{0,n}=H^0(G_k;(\hMU\Zln)^n)\cong \Zln\otimes MU^n$ and $E_2^{2,n-2}=H^2(G_k;(\hMU\Zln)^{n-2})\cong \Z/\ell^{\nu_0} \otimes MU^{n-2}$. Since there only two non-vanishing terms, these two groups are related by an exact sequence of $\Zln$-modules
$$0 \to E_2^{2,n-2} \to \hMU_{\et}^n(k;\Zln) \to E_2^{0,n} \to 0.$$
But since $E_2^{0,n}\cong \Zln \otimes MU^n$ is a free $\Zln$-module, this sequence splits. Hence, if $n$ is even, we get an isomorphism 
$$\hMU_{\et}^n(k;\Zln) \cong \Z/\ell^{\nu_0} \otimes MU^{n-2} \oplus \Zln \otimes MU^n.$$
\end{proof}
\subsection{Etale Cobordism for smooth curves}
The Atiyah-Hirzebruch spectral sequence and the knowledge on the \'etale cohomology allows us to determine the \'etale cobordism of  a smooth curve. 
\begin{theorem}\label{ecofsmoothcurve}
Let $k$ be a separably closed field of characteristic $p \neq \ell$ and let $X$ be a connected projective smooth curve of genus $g$ over $k$. Then the profinite \'etale cobordism of $X$ with $\Zln$-coeffiecients is given by  
\[
\hMU^n_{\et}(X;\Zln) \cong \left \{ \begin{array}{r@{\quad: \quad}l}
                                                  (\Zln \otimes MU^n) \oplus (\Zln \otimes MU^{n-2}) & n ~\mathrm{even} \\
                                                 \oplus_{i=1}^{i=2g} (\Zln \otimes MU^{n-1}) & n ~\mathrm{odd}.
                                                 \end{array} \right. 
\] 
\end{theorem}
\begin{proof}
By \cite{sga412}, Arcata, III, Corollaire 3.5, we know that the cohomology groups $H_{\et}^s(X;\Zln)$ vanish for $s >2$ and are free over $\Zln$ of rank $1, 2g, 1$ for $s=0,1,2$. The $E_2^{s,t}$-terms of the AHSS do not vanish only for $s=0,1,2$ and $t\leq 0$, $t$ even. Hence the spectral sequence collapses at the $E_2$-stage and since the occuring terms are all free $\Zln$-modules, the $n$th \'etale cobordism group $\hMU_{\et}^n(X;\Zln)$ is the direct sum of the $E_2^{s,t}$-terms with $s+t=n$. 
\end{proof}

\newpage

\begin{appendix}

\section{Existence of left Bousfield localization of fibrantly generated model categories}

The aim of the first part of this appendix is the proof of the localization Theorem \ref{leftloc} for the special situation of a left proper fibrantly generated model category $\Mh$, e.g. $\hSh$. The dual version of this theorem is the left localization theorem of \cite{hirsch}, but our result is more general. We only assume that the given model structure is fibrantly generated whereas in \cite{hirsch} it is used that the model structure is cellular, which is stronger than cofibrantly generated. The outline and most of the lemmas used to prove the theorem are just dual version of the ones in \cite{hirsch}. We will give the proof only of those statements that are not exact dual assertions. In particular, we point out at which places we use left properness and fibrantly generated and how we avoid the dual of cellularity.\\
We use this theorem in the next appendix to prove that there is a stable modle structure on the category of spectra for every fibrantly generated left proper model category. This enables us to stabilize the category of profinite spaces.\\
We assume the reader to be familiar with the notions of model categories, as explained e.g. in \cite{hirsch} or \cite{hoveybook}. 

\begin{convention}\label{convention}
We make the convention that a model category has by definition all finite limits and finite colimits as in \cite{homalg}. If they need to have more than finite limits or finite colimits we will indicate this. In our applications the categories have all small limits and all results we need only use the existence of small limits.
\end{convention}

In this subsection we prove an analogue of Theorem 4.1.1 of \cite{hirsch}. We recall the definition of local objects and local equivalences.
\begin{defn}\label{local}
Let $\Mh$ be a simplicial model category. Let $K$ be a set of objects in $\Mh$.
\begin{enumerate}
\item A map $f:A\to B$ is called a {\rm $K$-local equivalence} if for every element $X$ of $K$ the induced map of simplicial mapping spaces $f^{\ast}:\Map(B,X)\to \Map(A,X)$ is a weak equivalence of simplicial sets.
\item Let $\Ch$ denote the class of $K$-local equivalences. An object $X$ is called {\rm $K$-local} if it is $\Ch$-local,  i.e. if $X$ is fibrant and for every $K$-local equivalence $f:A\to B$ the induced
map $f^{\ast}:\Map(B,X) \to \Map(A,X)$ is a weak equivalence.
\end{enumerate}
\end{defn}

\begin{theorem}\label{leftloc}
Let $\Mh$ be a left proper fibrantly generated simplicial model category with all small limits. Let $K$ be a set of fibrant objects in $\Mh$ and let $\Ch$ be the class of $K$-local equivalences. 
\begin{enumerate}
\item The left Bousfield localization of $\Mh$ with respect to $\Ch$ exists, i.e.
there is a model category structure $\LCM$ on the underlying category $\Mh$ in which
\begin{enumerate}
\item the class of weak equivalences of $\LCM$ equals the class of $K$-local equivalences of $\Mh$,
\item the class of cofibrations of $\LCM$ equals the class of cofibrations of $\Mh$,
\item the class of fibrations of $\LCM$ is the class of maps with the right lifting property with respect
to those maps that are both cofibrations and $K$-local equivalences.
\end{enumerate}
\item The fibrant objects of $\LCM$ are the $K$-local objects of $\Mh$.
\item $\LCM$ is left proper. It is fibrantly generated if every object of $\Mh$
is cofibrant.
\item The simplicial structure of $\Mh $ gives $\LCM$ the structure of a simplicial model category.
\end{enumerate}
\end{theorem}

\begin{remark}
1. Although {\rm fibrantly generated} is the dual notion of {\rm cofibrantly generated} the proof of theorem \ref{leftloc} is not the dual of Theorem 4.1.1 of \cite{hirsch} for left Bousfield localization. But we can use a dual proof of Theorem 5.1.1 of \cite{hirsch} for right Bousfield localization where the cofibrations and fibrations change their role plus an idea of \cite{prospectra} in order to avoid the dual of Proposition 5.2.3 the proof of which uses the cellularization. The fact that we can construct a {\rm left} localization is due to the left properness. \\
2. The theorem does not depend on the simplicial structure of $\Mh$. In fact, it
is also true for a non-simplicial $\Mh$ but we would have to use the machinery of
homotopy function complexes instead of just simplicial mapping spaces.\\
3. Note that no result that we need below makes use of the fact that more than finite colimits exist in $\Mh$. Only in order to verify the factorization model axiom we need the existence of small limits since we use a variant of the cosmall object argument.
\end{remark}

We split the proof into some lemmas that we collect from \cite{hirsch} and \cite{prospectra}. Let $P$ be the set of generating fibrations in $\Mh$ and let $Q$ denote the set of generating trivial fibrations. Let $K$ be a fixed set of fibrant objects in $\Mh$. We define the two sets of maps
$$\Lambda (K):= \{ f:X^{\Delta^n} \to X^{\partial\Delta^n}, X\in K, n\geq0\}$$
and
$$\overline{\Lambda(K)} := \Lambda (K) \cup Q.$$

\begin{lemma}\label{prop524}
Let $\Mh$ be fibrantly generated model category and let $K$ be a set of objects of $\Mh$. If $B$ is a cofibrant object of $\Mh$ then a map $g:A\to B$ is $\overline{\Lambda(K)}$-projective if and only if it is both a cofibration and a $K$-local equivalence.
\end{lemma}
\begin{proof}
The proof is dual to the one of Proposition 5.2.4 of \cite{hirsch}.
\end{proof}

\begin{lemma}\label{prop525}
Let $\Mh$ be fibrantly generated model category and let $K$ be a set of objects of $\Mh$. Then a map that is a transfinite tower of pullbacks of maps in $\overline{\Lambda(K)}$ is a $K$-local fibration.
\end{lemma}
\begin{proof}
The proof is dual to the one of Proposition 5.2.5 of \cite{hirsch}. Note that it does not make use of the cellularity of $\Mh$.
\end{proof}

\begin{lemma}\label{lemma532}
Let $\Mh$ be a model category and $K$ a set of objects. A map $g:X\to Y$ is both a $K$-local fibration and a $K$-local weak equivalence if and only if it is a trivial fibration.
\end{lemma}
\begin{proof}
This is a dual of Lemma 5.3.2 of \cite{hirsch} whose proof does not depend on the cellularity or right properness of $\Mh$. Hence the dual proof of Lemma 5.3.2 gives a proof for $\Mh$ left proper and fibrantly generated.
\end{proof}

\begin{lemma}\label{prop533}
Let $\Mh$ be a model category and $K$ a set of objects. Then there is a functorial factorization for every map $g:X\to Y$ in $\Mh$ as $X\stackrel{q}{\to} W \stackrel{p}{\to} Y$ in which $q$ is a $K$-local cofibration and $p$ is both a $K$-local weak equivalence and a $K$-local fibration.
\end{lemma}
\begin{proof}
This follows from Lemma \ref{lemma532} and the functorial factorization into a cofibration followed by a trivial fibration.
\end{proof}

\begin{lemma}\label{lemma534}
Let $\Mh$ be a left proper fibrantly generated model category and $K$ a set of objects. If $g:X\to Y$ is a weak equivalence, $h:Y\to Z$ is a fibration and the composition $hg:X\to Z$ is a $K$-local fibration, then $h$ is a $K$-local fibration.
\end{lemma}

\begin{proof}
This is the dual proof to Lemma 5.3.4 in \cite{hirsch}, but we point out where we use left properness and how the diagram looks like.\\
If $i:A\to B$ is a map in $K-W \cap K-\mathrm{Cof} = W \cap \mathrm{Cof}$ (last equality holds by Lemma \ref{lemma532}), then one can choose by Proposition 8.1.23 of \cite{hirsch} a cofibrant approximation
$\tilde{i}: \tilde{A} \to \tilde{B}$ to $i$ such that $\tilde{i} \in$ Cof. By Propositions 3.1.5 and 3.2.3 of \cite{hirsch} $\tilde{i}$ is a $K$-local equivalence. Since $\Mh$ is left proper, Proposition 13.2.1 of \cite{hirsch} implies that it is sufficient to show that $h$ has the right lifting property with respect to $\tilde{i}$.\\
Suppose we have a diagram
$$\begin{array}{rcl}
 &  & X\\
 &  & \downarrow g \\
\tilde{A} & \stackrel{s}{\to} & Y \\
\tilde{i}\downarrow &  & \downarrow h \\
\tilde{B} & \stackrel{t}{\to} & Z.
\end{array}$$
Now we continue as in the proof of Lemma 5.3.4 in \cite{hirsch} by constructing lifts $j:\tilde{A} \to X$ and $k:\tilde{B} \to X$, consider the composition $u:=gk:\tilde{B}\to Y$ which has the property $ui \simeq s$ in the category of objects over $Z$ and conclude by the homotopy extension property of cofibrations in this case (Proposition 7.3.10 of \cite{hirsch}) that there is a map $v:\tilde{B}\to Y$ such that $v\simeq u$, $v\tilde{i}=s$ and hence $hv=t$.
\end{proof}

In \cite{hirsch} the cellularity of the given model category is used for the existence of factorizations. We have to avoid this point. In order to do so, we need a dual version of Lemma 2.5 of \cite{prospectra}.

\begin{lemma}\label{lemma25}
Let $\Mh$ be a left proper fibrantly generated simplicial model category and let $K$ be a set of fibrant objects. Then every map $g:A\to B$ has a factorization into a cofibration $s:A\to W$ that has the left lifting property with respect to all maps $X^{\Delta^n} \to X^{\partial\Delta^n}$, $X \in K$, $n\geq 0$
followed by a $K$-fibration.
\end{lemma}

\begin{proof}
This is a variation of the small object argument, dual to the proof of Lemma 2.5 of \cite{prospectra}. We prove this lemma by transfinite induction of length $\kappa$ where $\kappa$ is a regular cardinal such that the codomains of the generating trivial fibrations $Q$ and the objects of $K$ are $\kappa$-cosmall relative to the fibrations.\\
Let $J_o$ be the set of all squares
$$\begin{array}{rcc}
A & \to & X^{\Delta^n} \\
g \downarrow &  & \downarrow \\
B & \to &  X^{\partial\Delta^n}
\end{array}$$
for all objects $X\in K$. Define $Z_0$ to be the pullback
$\prod_{J_0}X^{\Delta^n} \times_{\prod_{J_0}X^{\partial\Delta^n}} B$ and let $j_0:A\to Z_0$, $q_0:Z_0\to B$ be the obvious maps. Now we factor $j_0$ into a cofibration $i_0:A\to W_0$ followed by a
trivial fibration $p_0:W_0 \to Z_0$. This is the first step of the induction.\\
If $\beta$ is a limit ordinal, set $W_{\beta}:=\lim_{\alpha < \beta} W_{\alpha}$ and set
$i_{\beta}:= \colim_{\alpha < \beta} i_{\alpha}$.\\
If $\beta$ is a successor ordinal, define $J_{\beta}$ to be the set of all squares
$$\begin{array}{rcc}
A & \to & X^{\Delta^n} \\
g \downarrow &  & \downarrow \\
W_{\beta -1} & \to &  X^{\partial\Delta^n}
\end{array}$$
for all $X\in K$. Define $Z_{\beta}$ to be the pullback $\prod_{J_{\beta}}X^{\Delta^n} \times_{\prod_{J_{\beta}}X^{\partial\Delta^n}} W_{\beta -1}$
and let $j_{\beta}:A\to Z_{\beta}$, $q_{\beta}:Z_{\beta}\to Z_{\beta}$ be the obvious maps.
Now we factor $j_{\beta}$ into a cofibration $i_{\beta}:A\to W_{\beta}$ followed by a
trivial fibration $p_{\beta}:W_{\beta} \to Z_{\beta}$.\\
By transfinite induction we get a $\kappa$-tower 
$$\ldots \to W_{\beta} \to \ldots  \to W_1 \to W_0 \to B.$$
Note that there are compatible maps $i_{\beta}:A\to W_{\beta}$. Since filtered limits exist in $\Mh$ we may continue with the following definition. Let $W$ be $\lim_{\beta} W_{\beta}$ and let $i : A \to W$ be
$\lim_{\beta} i_{\beta}$. By construction, each $i_{\beta}$ is a cofibration for every successor ordinal $\beta$. Since $\Mh$ is fibrantly generated and $\kappa$ is such that the codomains of $Q$ and the objects of $K$ are $\kappa$-cosmall relative to the fibrations we get that $i$ is a cofibration. Furthermore one can show that $i$ has the left lifting property with respect to all maps $X^{\Delta^n} \to X^{\partial\Delta^n}$, $X \in K$, $n\geq 0$ by the cosmall object argument as follows:\\
Suppose we have a diagram
$$\begin{array}{rcl}
A & \stackrel{k}{\to} & X^{\Delta^n} \\
i \downarrow &  & \downarrow f\\
W & \stackrel{h}{\to} &  X^{\partial\Delta^n}
\end{array}$$
for some $X\in K$ and $n\geq0$. By the theory of simplicial model categories if $X$ is a fibrant object then the map $f:X^{\Delta^n} \to X^{\partial\Delta^n}$ is a fibration. Since the objects in $K$ are $\kappa$-cosmall relative to the fibrations, there is a $\beta < \kappa$ such that $h$ is the composite $W \to W_{\beta} \stackrel{h_{\beta}}{\to} X^{\partial\Delta^n}$. By construction, there is a map $k_{\beta}:W_{\beta +1} \to X^{\Delta^n}$ such that $fk_{\beta}=h_{\beta}j$ and $k=k_{\beta}i_{\beta +1}$ where
$j$ is the map $W_{\beta +1} \to W_{\beta}$. The composition $W \to W_{\beta +1} \stackrel{k_{\beta}}{\to} X^{\Delta^n}$ is the required lift in the diagram.\\
It remains to show that $p:W \to Y$ is a $K$-fibration. Since $K$-fib is defined by a right lifting property and $p$ is a transfinite tower of maps $W_{\beta +1} \to W_{\beta}$ it suffices to show that each $W_{\beta +1} \to W_{\beta}$ is a $K$-fibration. The map $W_{\beta +1} \to W_{\beta}$ is the composition $W_{\beta +1} \to Z_{\beta +1} \to W_{\beta}$ where the first map is a trivial fibration and hence also a map in $K-W \cap K-\fib$. The second map is a product of maps of the form $X^{\Delta^n} \to X^{\partial\Delta^n}$ each of which is a $K$-fibration by Lemma \ref{prop524}. Hence $Z_{\beta +1} \to W_{\beta}$ is a $K$-fibration.
\end{proof}

\begin{lemma}\label{prop535}
Let $\Mh$ be a left proper fibrantly generated model category and let $K$ be a set of objects. Then there is functorial factorization for every map $g:X\to Y$ in $\Mh$ as $X\stackrel{q}{\to} W \stackrel{p}{\to} Y$ in which $q$ is both a $K$-local weak equivalence and a $K$-local cofibration and $p$ is a $K$-local fibration.
\end{lemma}
\begin{proof}
Again we dualize the proof of Proposition 5.3.5 of \cite{hirsch} point out how the diagrams must look
like and where we use the left properness. Let $j:\tilde{X}\to X$ be a fibrant cofibrant approximation to $X$. Lemma \ref{lemma25} implies that there is a functorial factorization of the composition $gj$
as $\tilde{X} \stackrel{s}{\to} \tilde{W} \stackrel{r}{\to} X$ such that $s$ is a cofibration and in $\overline{\Lambda (K)}$-projective and $r$ is a $K$-local fibration. Let $Z:=X\amalg_{\tilde{X}}\tilde{W}$. We factor the induced map $Z\to Y$ into $Z \stackrel{u}{\to} W \stackrel{p}{\to} Y$ where $u$ is a trivial cofibration and $p$ is a fibration. This yields the following diagram 
 
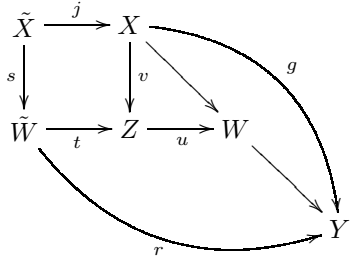
\[
\xymatrix{
\tilde{X} \ar[r]^j \ar[d]_s & X \ar[d]^v \ar[dr] \ar@/^2pc/[ddrr]^g \\
\tilde{W} \ar[r]_t \ar@/_2pc/[drrr]_r & Z \ar[r]_u & W \ar[dr] \\
 & & & Y 
}
\]
Let $q:=uv$. Now since $j$ is a weak equivalence and $s$ is a cofibration the left properness of $\Mh$ implies that $t$ is also a weak equivalence. Hence $ut$ is a weak equivalence and $s$ is a cofibrant approximation to $q$. Lemma \ref{prop524} implies that $s$ is a $K$-local equivalence.\\
By Lemma \ref{prop524}, since $s$ is in $\overline{\Lambda (K)}$-proj it is in particular a $K$-local equivalence. Hence $q\in K-W$ since $j \in K-W$ and $qj=ts \in K-W$ (two out of three holds for $K-W$).
Now since $q=uv$ is the composition of two cofibrations and the $K$-local cofibrations equal the cofibrations in $\Mh$, we have $q\in K-W \cap K-Cof$. By Lemma \ref{prop525} this q is a weak equivalence and a cofibration.\\
Since $u$, $t$ and hence also $ut$ are weak equivalences and $r=put$ is by construction a $K$-local fibration, we conclude by Lemma \ref{lemma534} that $p$ is a $K$-local fibration.
\end{proof}

\begin{proof} {\bf of Theorem \ref{leftloc}}
1. Axiom M1 is satisfied since it is satisfied in $\Mh$. Axiom M2 follows from Proposition 3.2.3 of
\cite{hirsch}. Axiom M3 follows from Proposition 3.2.4 and Lemma 7.2.8 of \cite{hirsch}.
The right lifting property of $K$-local fibrations is given by definition.
The left lifting property of $K$-local cofibrations follows from Lemma \ref{lemma532}.
The factorization of a map into a $K$-local cofibration followed by a trivial $K$-local fibration
follows from Lemma \ref{prop533} and the factorization into a trivial $K$-local cofibration followed
by a $K$-local fibration follows from Lemma \ref{prop535}.\\
2. This follows from Proposition 3.4.1 of \cite{hirsch}.\\
3. The left properness follows from Proposition 3.4.4 of \cite{hirsch}.
If every object in $\Mh$ is cofibrant, then Lemma \ref{lemma25} and lemma \ref{prop525}
imply that $\LCM$ is fibrantly generated with generating fibrations
$\overline{\Lambda(K)}$ and generating trivial fibrations $Q$.\\
4. This is the analogue proof as for part 4 of Theorem 5.1.1 in \cite{hirsch}.
\end{proof}

\section{Stable model structure on spectra}

The well known stabilization of the category of simplicial spectra of \cite{bousfried}
uses the fact that $\Sh$ is proper in an essential way in order to construct functorial factorizations. 
Since the category $\hSh$ is only left but not right proper, we may not use their localization method.\\
Hovey \cite{hovey} has pointed out a general way to stabilize a left proper cellular model category with respect to a left Quillen endofunctor $T$. Since we do have left properness but not the property of being
cofibrantly generated we have to modify this construction.\\
The aim is to construct a stable structure on $\Ch$ with respect to an endofunctor via the well known construction of spectra on $\Ch$. We reformulate the results of \cite{hovey} for the situation of a left proper fibrantly generated model category $\Ch$, e.g. $\hSh$. The projective and stable model structure is the same as in \cite{hovey}, but we have to give different proofs for the key points of the construction. We also have to choose different generating fibrations. But the idea for the construction and the proof that the structure in fact becomes stable are due to \cite{hovey}.\\
In the first part, we construct a stable structure on the usual classical notion of spectra in analogy to the one of \cite{bousfried}. \\
In the second part, we show, following \cite{hovey}, that this works also for symmetric spectra on a symmetric monoidal model structure that is left proper and fibrantly generated. 

\subsection{Spectra}

In this section we assume always that $\Ch$ is a left proper fibrantly generated simplicial model category with all small limits. Let $T:\Ch \to \Ch$ be a left Quillen endofunctor on $\Ch$. Let $U:\Ch \to \Ch$ be its right adjoint.

\begin{defn}
A {\rm spectrum} $X$ is a sequence $(X_n)_{n\geq 0}$ of objects of $\Ch$ together
with structure maps $\sigma :TX_n \to X_{n+1}$ for all $n\geq 0$. A {\rm map of spectra}
$f:X\to Y$ is a collection of maps $f_n:X_n \to Y_n$ commuting with the structure maps.
We denote the category of spectra by $\SpC$.
\end{defn}

We begin by defining an intermediate strict structure, see \cite{bousfried} or \cite{hovey}. But note that we are in a dual situation of a fibrantly generated model structure. Consider the following functors:

\begin{defn}\label{FnRn}
Given $n\geq 0$, the {\rm evaluation functor} $\Ev_n:\SpC \to \Ch$ takes $X$ to $X_n$. It has a left adjoint $F_n:\Ch \to \SpC$ defined by $(F_nA)_m=T^{m-n}A$ if $m\geq n$ and $(F_nA)_m=\ast$ otherwise. The structure maps are the obvious ones.\\
The evaluation functor also has a right adjoint $R_n:\Ch \to \SpC$ defined by $(R_nA)_i=U^{n-i} A$ if $i\leq n$ and $(R_nA)_i=\ast$ otherwise. The structure map $TU^{n-i}A \to U^{n-i} A$ is adjoint to the identity map of $U^{n-i}$ when $i<n$.
\end{defn}

Note that $Ev_n$ commutes with limits and colimits, where limits and colimits are constructed levelwise in $\SpC$. If $G: I \to \SpC$ is a functor from a small category to $\SpC$, then, for limits, the structure maps are defined as the composites 
$$T(\lim \Ev_n \circ G)  \to \lim (T\circ \Ev_n \circ G) \stackrel{\lim(\sigma \circ G)}{\to} \lim \Ev_{n+1} \circ G.$$ 
For colimits one uses the fact that $T$ is left adjoint and hence preserves colimits to get the obvious structure maps.

\begin{defn}\label{defstrict}
A map $f$ in $\SpC$ is a {\rm projective weak equivalence} (resp. {\rm projective fibration}) if each map $f_n$ is a weak equivalence (resp. fibration). A map $i$ is a {\rm projective cofibration} if it has the left lifting property with respect to all projective trivial fibrations. 
\end{defn}

%
We will show that the projective structure is in fact a fibrantly generated model structure. We denote the set of generating fibrations of $\Ch$ by $P$ and the set of generating trivial fibrations by $Q$. Inspired by Proposition \ref{projcofs}, we set 
$$\tP:= 
\{g_n: R_nR \to R_nS\times_{R_{n-1}US}R_{n-1}UR, ~ \mathrm{for~all} ~ f:R\to S ~\mathrm{in} ~ P\}$$ 
and 
$$\tQ:=
\{g_n: R_nR \to R_nS\times_{R_{n-1}US}R_{n-1}UR, ~ \mathrm{for~all} ~ f:R\to S ~\mathrm{in} ~ Q\}$$ 
where $R_n$ is the right adjoint of the evaluation functor $Ev_n$ and $g_n$ is the map induced by the commutative diagram
\begin{equation}\label{diagRn}
\begin{array}{ccc}
R_nR& \to & R_nS \\
\downarrow &  & \downarrow \\
R_{n-1}UR & \to & R_{n-1}US.
\end{array}
\end{equation}
We will show that the maps in $\tP$ are the generating fibrations and the maps in $\tQ$ are the generating trivial fibrations for the projective model structure on $\SpC$. 

\begin{lemma}\label{levelfibs}
If $f:X \to Y$ is a fibration (resp. trivial fibration) in $\Ch$, then the maps 
$g_n: R_nX \to R_nY\times_{R_{n-1}UY}R_{n-1}UX$ are fibrations (resp. trivial fibrations) in $\SpC$.
\end{lemma}
\begin{proof}
Since $R_n$ is defined via the right Quillen functor $U$ which preserves fibrations and trivial fibrations, it is clear that the map $(R_nf)_i=U^{n-i}(f)$ on the $i$-th level is a fibration (resp. trivial fibration). In  diagram (\ref{diagRn}), the maps on the $i$-th level are either identities or equal to $U^{n-i}(f)$. Hence the induced map $g_n$ is also a fibration (resp. trivial fibration).
\end{proof}

The definition of $\tP$ and $\tQ$ becomes clear after the following 
\begin{prop}\label{projcofs}
A map  $i:A \to B$ in $\SpC$ is a projective (trivial) cofibration if and only if $i_0:A_0 \to B_0$ and the induced maps  $j_n:A_n \amalg_{TA_{n-1}}TB_{n-1}$ for $n\geq 1$ are (trivial) cofibrations in $\Ch$. 
\end{prop}
\begin{proof}
We prove the cofibration case, the case of trivial cofibrations is similar. We use the standard idea of \cite{hovey}, Proposition 1.14. Let $i:A\to B$ be a map with the left lifting property with respect to projective trivial fibrations. Since projective fibrations are defined levelwise, it is clear that $i_0:A_0 \to B_0$ is a cofibration. Let $f:X\to Y$ be a trivial fibration in $\Ch$. We have to show that there is a lift $B_n \to X$ for any commutative diagram
$$\begin{array}{ccc}
A_n \amalg_{TA_{n-1}}TB_{n-1}& \to & X \\
j_n \downarrow &  & \downarrow f \\
B_n & \longrightarrow & Y.
\end{array}$$
By adjunction, this diagram has a lift if and only if the induced diagram
$$\begin{array}{ccc}
A& \longrightarrow & R_nX \\
i \downarrow &  & \downarrow g_n \\
B & \to & R_nY\times_{R_{n-1}UY}R_{n-1}UX
\end{array}$$
has a lift. Now $R_n$ is defined via the right Quillen functor $U$, see Definition \ref{FnRn}. This implies that all maps in Diagram (\ref{diagRn}) are projective trivial fibrations and hence the induced map is one, too. So $g_n$ is a trivial fibration if $f$ is one, see Lemma \ref{levelfibs}. Hence there is a lift in the last diagram.\\
Now suppose that  $i_0$ and $j_n$ are cofibrations in $\Ch$ for all $n>0$. Let 
$$\begin{array}{ccc}
A& \stackrel{a}{\to} & X \\
i \downarrow &  & \downarrow f \\
B & \stackrel{b}{\to} & Y
\end{array}$$
be a commutative diagram in $\SpC$ and let $f$ be a trivial fibration. We construct a lift $h_n:B_n \to X_n$, compatible with the structure maps, by induction on $n$. There is a lift $h_0$ since $i_0$ is a cofibration  and $f_0$ is a trivial fibration in $\Ch$. Suppose $h_j$ is constructed for all $j<n$. Then we can consider the diagram
$$\begin{array}{ccc}
A_n\amalg_{TA_{n-1}}TB_{n-1}& \stackrel{(f_n,\sigma \circ Th_{n-1})}{\longrightarrow} & X \\
j_n \downarrow &  & \downarrow f_n \\
B_n & \stackrel{b_n}{\longrightarrow} & Y_n.
\end{array}$$
The lift $h_n:B_n \to X_n$ of this diagram is the desired map. 
\end{proof}

\begin{theorem}\label{thmstrict}
Let $\Ch$ be a left proper fibrantly generated simplicial model category with all small limits.
The projective weak equivalences, projective fibrations and projective cofibrations define a left proper
fibrantly generated simplicial model structure on $\SpC$ with set of generating fibrations $\tP$ and set of generating trivial fibrations $\tQ$.
\end{theorem}

We will prove this theorem in a series of lemmas. The existence of finite colimits and limits is clear since they are defined levelwise. The two-out-of-three- and the retract-axiom are immediate as well. 
The crucial point is the factorization axiom for which we use the cosmall object argument.
For the notations $\tP-\proj$, $\tP-\fib$, $\tQ-\proj$, $\tQ-\fib$ and the cosmall object argument, which is the exact dual of the small object argument, we refer the reader to the book of Hovey \cite{hoveybook}: 
in $\tP-\proj$ are all the maps having the left lifting property with respect to all maps in $\tP$; in $\tP-\fib$ are all the maps having the right lifting property with respect to maps in $\tP-\proj$.

\begin{prop}\label{tP-proj}
The projective cofibrations are exactly the maps which have the left lifting property with respect to all maps in $\tQ$. The projective trivial cofibrations are exactly the maps which have the left lifting property with respect to all maps in $\tP$. 
\end{prop}
\begin{proof}
By Proposition \ref{projcofs}, $i:A\to B$ is a projective cofibration if and only if $i_0$ and $j_n$ are cofibrations in $\Ch$ for all $n>0$. But since $\Ch$ is a fibrantly generated model category with $Q$ the set of generating trivial fibrations, $i_0$ and $j_n$ are cofibrations if and only if they have the left lifting property with respect to all maps in $Q$, i.e. if and only if $i_0$ and $j_n$ are in $Q-\proj$. 
We have seen that a lift in a diagram
$$\begin{array}{ccc}
A_n \amalg_{TA_{n-1}}TB_{n-1}& \to & R \\
j_n \downarrow &  & \downarrow f \\
B_n & \longrightarrow & S
\end{array}$$ 
is by adjunction equivalent to a lift in the diagram
$$\begin{array}{ccc}
A& \longrightarrow & R_nR \\
i \downarrow &  & \downarrow g_n \\
B & \to & R_nS\times_{R_{n-1}US}R_{n-1}UR.
\end{array}$$
Hence $i_0$ and $j_n$ are in $Q-\proj$ if and only if the map $i$ is an element of $\tQ-\proj$. This completes the proof in the cofibration case. The trivial cofibration case is similar. 
\end{proof}

\begin{cor}\label{tP-fib=fib}
The projective fibrations are exactly the maps in $\tP-\fib$. The projective trivial cofibrations are exactly the maps in $\tQ-\fib$.
\end{cor}

\begin{lemma}\label{cosmall}
If $Z$ is a cosmall object in $\Ch$ relative to the fibrations, then $R_nZ$, $n\geq 0$, is cosmall relative to the level fibrations in $\SpC$. Similarly, if $Z$ is a cosmall object in $\Ch$ relative to the trivial fibrations, then $R_nZ$, $n\geq 0$, is cosmall relative to the level trivial fibrations in $\SpC$.
In addition, the codomains of the maps in $\tP$ (resp. $\tQ$) are cosmall relative to the level fibrations (resp. trivial fibrations) in $\SpC$.
\end{lemma}
\begin{proof}
This is clear since $Ev_n$ commutes with all limits.
\end{proof} 

We can now finish the proof of Theorem \ref{thmstrict}. It remains to show the factorization axiom.
Note that by hypothesis, $\Ch$ and hence $\SpC$ have all small limits. This enables us to apply the cosmall object argument.  By Lemma \ref{cosmall} and since we have shown in Corollary \ref{tP-fib=fib} that the class of projective fibrations (resp. trivial fibrations) is equal to $\tP-\fib$ (resp. $\tQ-\fib$), we can apply the cosmall object argument to the class $\tP$ (resp. $\tQ$). In the first case, this yields a factorization into a map in $\tP-\proj$ followed by a map in $\tP-\cocell$, where $\tP-\cocell$ consists of maps which are transfinite compositions of pullbacks of maps in $\tP$. But we have shown that $\tP-\proj$ is equal to the class of projective trivial cofibrations and $\tP-\cocell \subseteq \tP-\fib$, \cite{hoveybook} dual to Lemma 2.1.10. Since we have shown that $\tP-\fib$ is equal to the class of projective fibrations, this is the first desired factorization.\\
The cosmall object argument applied to $\tQ$ yields a factorization into a map in $\tQ-\proj$ followed by a map in $\tQ-\cocell$. Since we know that $\tQ-\proj$ is equal to the class of projective cofibrations and $\tQ-\cocell \subseteq \tQ-\fib$ and since $\tQ-\fib$ is equal to the class of projective trivial fibrations, we have the second factorization.\\
This proves that the projective structure is in fact a model structure on $\SpC$. But we have also shown that this structure satisfies the axioms of a fibrantly generated model category, see \cite{hoveybook}, Definition 2.1.17. \\
The left properness is clear since coproducts are defined levelwise in $\SpC$.
This completes the proof of Theorem \ref{thmstrict}. $\Box$  

It remains to modify this structure in order to get a stable structure, i.e. one in which the prolongation of $T$ is a Quillen equivalence. We will do this by applying the localization Theorem \ref{leftloc} to the projective model structure on spectra. We want the stable weak equivalences to be the maps that induce isomorphisms on all generalized cohomology theories. A generalized cohomology theory is represented by the analogue of an $\Omega$-spectrum. Since the $K$-local objects are the fibrant objects in the localized model structure, we have to choose the set $K$ to consist exactly of these analogues of $\Omega$-spectra.

\begin{defn}
A spectrum $E\in \SpC$ is defined to be an {\rm $\Omega$-spectrum} if each $E_n$  is fibrant and the adjoint structure maps $E_n \to UE_{n+1}$ are weak equivalences for all $n\geq 0$.
\end{defn}

By Corollary 9.7.5 in \cite{hirsch}, since each $E_n$ is fibrant and since the right Quillen functor $U$ preserves fibrations, we have that $E_n \to UE_{n+1}$ is a weak equivalence in $\Ch$ if and only if the induced map $\Map (A,E_n) \to \Map (A, UE_{n+1})$ is a weak equivalence in $S$ for every cofibrant object $A$. By adjunction this is equivalent to $\Map (F_nA,E) \to \Map (F_{n+1}TA, E)$ being a weak equivalence in $S$ for every cofibrant object $A$. Note that by \cite{hirsch}, Corollary 9.7.5, this is what the weak equivalences have to be in a model structure. \\
So in order to get the $\Omega$-spectra as stable fibrant objects and the maps inducing isomorphisms on cohomology theories as stable equivalences we have to choose the maps $F_{n+1}TA \to F_nA$ adjoint to the identity map of $TA$ to be the $K$-local equivalences. Hence we define the set $K$ to consist of all $\Omega$-spectra.\\
Using the fact that the projective model structure on $\SpC$ is fibrantly generated, left proper and simplicial, we may deduce the following 
\begin{theorem}\label{stablestructure}
Let $\Ch$ be a left proper fibrantly generated simplicial model category with all small limits. Let $K$ be the set of $\Omega$-spectra. There is a {\rm stable model structure} on $\SpC$ which is defined to be the $K$-localized model structure $L_S\SpC$ of Theorem \ref{leftloc} where $S$ is the class of all $K$-local equivalences.
\end{theorem}
\begin{remark}
In particular, the {\rm stable equivalences} are the $K$-local equivalences; the {\rm stable cofibrations} are the projective cofibrations; the {\rm stable fibrations} are the maps that have the right lifting property with respect to all stable trivial cofibrations; the stable fibrant objects are the $\Omega$-spectra.
\end{remark}

Define the {\em prolongation of $T$} to be the functor $T:\SpC \to \SpC$, defined by $(TX)_n=TX_n$ with structure maps $T(TX_n)\stackrel{T\sigma}{\to}TX_{n+1}$ where $\sigma$ is the structure map of $X$. The adjoint $\Omega$ is defined in the analog way. Just as in \cite{hovey} one can prove that $T$ is a Quillen equivalence on the stable structure.

\begin{theorem}\label{Tstable}
Let $\Ch$ be a left proper fibrantly generated simplicial model category with a left Quillen endofunctor $T$. Then the prolongation $T:\SpC \to \SpC$ is a Quillen equivalence with respect to the stable structure.
\end{theorem}
\begin{proof}
This is the proof of Theorem 3.9 in \cite{hovey}.
\end{proof}

\subsection{Symmetric spectra}

We can transfer these results also to symmetric spectra on a left proper fibrantly generated simplicial model category $\Ch$ with a symmetric monoidal structure. The ideas and definitions are due to \cite{hovey}. We extend his results to our unusual and slightly more general situation. For the definition of these notions, we refer the reader to \cite{hoveybook}.\\
For the whole section, let $\Ch$ be a left proper fibrantly generated simplicial closed symmetric monoidal model category with all small limits. For simplicity, we suppose that its unit object $S$ is cofibrant. Since any endofunctor $T:\Ch \to \Ch$ is of the form $T(L)=L\otimes K$ for $K=T(S)$, we will consider in this section the functor $- \otimes K$ instead of $T-$. The right adjoint of $T$ is then given by $L\to L^K$. All the definitions and ideas for the construction of symmetric spectra are taken from \cite{hovey}.

\begin{defn}
A {\rm symmetric sequence} in $\Ch$ is a sequence $X_0, X_1, \ldots , X_n, \ldots$ of objects in $\Ch$ with an action of $\Sigma_n$ on $X_n$, where $\Sigma_n$ is the symmetric group of $n$ letters. A map of symmetric sequences is a sequence of $\Sigma_n$-equivariant maps $X_n \to Y_n$. We denote the resulting category by $\Ch^{\Sigma}$.\\
We define the object $\SymK$ in $\Ch^{\Sigma}$ to be the symmetric sequence $(S, K, K\otimes K, \ldots , K^{\otimes n}, \ldots )$ where $\Sigma_n$ acts on $K^{\otimes n}$ by permutation, using the commutativity and associativity isomorphisms.   
\end{defn}

\begin{defn}\label{defsymmspectra}
A {\rm symmetric spectrum} $X$ is a sequence of $\Sigma_n$-objects $X_n$ in $\Ch$ with $\Sigma_n$-equivariant structure maps $X_n\otimes K \to X_{n+1}$, such that the composite
$$X_n \otimes K^{\otimes p} \to X_{n+1}\otimes K^{\otimes p-1} \to \ldots \to X_{n+p}$$
is $\Sigma_n \times \Sigma_p$-equivariant for all $n, p \geq 0$. A map of symmetric spectra is a collection of $\Sigma_n$-equivariant maps $X_n \to Y_n$ compatible with the structure maps. 
We denote the category of symmetric spectra on $\Ch$ by $\SpsC$.\\
We define the symmetric spectrum $\SymKb$ in $\SpsC$ to be the initial object $0$ in degree $0$ and $K^{\otimes n}$ in degree $n$ for $n>0$ with the obvious structure maps.
\end{defn}

One can show as in \cite{symmspectra} that $\SpsC$ is again a closed symmetric monoidal category with $\SymK$ as the unit. We denote the monoidal structure by $X \wedge Y = X \otimes_{\SymK}Y$, and the closed structure by $\Hom_{\SymK}(X,Y)$, where we consider $X$ and $Y$ as $\SymK$-modules. 

\begin{defn}\label{symmFnRn}
1. Given $n\geq 0$, the {\rm evaluation functor} $\Ev_n:\SpsC \to \Ch$ takes $X$ to $X_n$.\\
2. The evaluation functor has a left adjoint $F_n:\Ch \to \SpsC$, defined by $F_nA=\tilde{F}_nA \otimes \SymK$, where $\tilde{F}_nA$ is the symmetric sequence $(0, \ldots , 0, \Sigma_n \times A, 0,Ê\ldots)$.\\
3. The evaluation functor has a right adjoint $R_n:\Ch \to \SpsC$, defined by $R_nA=\Hom(\SymK, \tilde{R}_nA)$, where $\tilde{R}_nA$ is the symmetric sequence that is the terminal object in dimensions other than $n$ and is the cofree $\Sigma_n$-object $\Hom_{\Ch}(\Sigma_n,A)$, i.e. $n!$-product of $A$ in $\Ch$ together with the $\Sigma_n$-action defined by $(\rho f)(\rho')=f(\rho'\rho)$ for $\rho \in \Sigma_n$ and $f\in \Hom_{\Ch}(\Sigma_n,A)$.  
\end{defn}

\begin{defn}\label{symmdefstrict}
A map $f$ in $\SpsC$ is a {\rm projective weak equivalence} (resp. {\rm projective fibration}) if each map $f_n$ is a weak equivalence (resp. fibration). A map $i$ is a {\rm projective cofibration} if it has the left lifting property with respect to all projective trivial fibrations. 
\end{defn}

Again we will show that the projective structure is in fact a fibrantly generated model structure. We recall that the set of generating fibrations of $\Ch$ is denoted by $P$ and the set of generating trivial fibrations by $Q$. Inspired by the forthcoming proposition we set 
\begin{eqnarray*}
\tP & := & 
\{q_n: R_nX \to R_nY\times_{R_n\Hom_K(K^{\otimes n}, Y)} R_n\Hom_K(K^{\otimes n}, X),\\
 & & \hspace*{2cm} ~ \mathrm{for~all} ~ p:X \to Y ~\mathrm{in} ~ P\}\end{eqnarray*}
and 
\begin{eqnarray*}
\tQ & := & 
\{q_n: R_nX \to R_nY\times_{R_n\Hom_K(K^{\otimes n}, Y)} R_n\Hom_K(K^{\otimes n}, X),\\
& & \hspace*{2cm} ~ \mathrm{for~all} ~ p:X \to Y ~\mathrm{in} ~ Q\}\end{eqnarray*}
where $R_n$ is the right adjoint of the evaluation functor $Ev_n$ and $q_n$ is the induced map
in the commutative diagram
\begin{equation}\label{symmdiagRn}
\begin{array}{ccc}
R_nX& \to & R_nY \\
\downarrow &  & \downarrow \\
R_n(\Hom_K(K^{\otimes n},X)) & \to & R_n(\Hom_K(K^{\otimes n},Y)).
\end{array}
\end{equation}
We will show that the maps in $\tP$ are the generating fibrations and the maps in $\tQ$ are the generating trivial fibrations for the projective model structure on $\SpC$. First we describe the projective cofibrations in $\SpsC$. This is slightly more complicated than in the case of usual spectra. The reason for this comes from the fact that we have to ensure that all maps are equivariant with respect to the actions of the groups $\Sigma_n$. The following definition and the idea how to solve this problem is taken from \cite{hovey}, Definition 8.4.\\
The following proposition is well-known, see \cite{hovey}, Proposition 8.5, or \cite{symmspectra}, Proposition 5.2.2. But we will give a different proof similar to the one of Proposition \ref{projcofs}, which enables us to understand the definition of the generating fibrations and trivial fibrations in the projective model structure on $\SpsC$. It avoids the assumption that we had already shown that the projective structure is in fact a model structure. Furthermore, we do not have to assume that $\Ch$ is cofibrantly generated, which is not the case in our situation. 

\begin{defn}\label{latchingspace}
We define the {$n$-th latching space} of $A \in \SpsC$ by $L_nA:=\Ev_n(A \wedge \SymKb)$. The obvious map $\SymKb \to \SymK$ induces a map of spectra $i:A \wedge \SymKb \to A\wedge \SymK$ and a $\Sigma_n$-equivariant map $\Ev_n(i):L_nA \to A_n$.  
\end{defn}
Note that the latching space is a $\Sigma_n$-object in $\Ch$. There is a model structure on the category of $\Sigma_n$-objects of $\Ch$, where the fibrations and weak equivalences are the underlying ones.\\
We recall that for given maps $f:A\to B$ and $g:C\to D$ in $\Ch$ one defines the pushout product $f \Box g$ of $f$ and $g$ to be the map 
$$f\Box g:(A\otimes D)\amalg_{A\otimes C}(B\otimes C) \to B \otimes D$$ 
induced by the obvious commutative diagram.  
\begin{prop}\label{symmprojcofs}
A map $f:A \to B$ in $\SpsC$ is a projective cofibration if and only if the induced map $\Ev_n(f \Box i) :A_n \amalg_{L_nA} L_nB \to B_n$ is a $\Sigma_n$-cofibration in $\Ch$ for all $n$. Similarly, $f$ is a projective trivial cofibration if and only if $\Ev_n(f \Box i)$ is a trivial $\Sigma_n$-cofibration in $\Ch$ for all $n$. 
\end{prop} 
\begin{proof}
As usual we prove only the cofibration case, since the case of trivial cofibrations is similar. Suppose first that each map $\Ev_n(f \Box i)$ is a $\Sigma_n$-cofibration in $\Ch$. The following argument is exactly the one in the proof of Proposition 8.5 in \cite{hovey}. \\
Suppose that $f$ is a projective cofibration of symmetric spectra. We have to show that each map $\Ev_n(f \Box i)$ is a $\Sigma_n$-cofibration in $\Ch$, i.e. each map $\Ev_n(f \Box i)$ has the left lifting property with respect to every trivial $\Sigma_n$-fibration $p:X \to Y$ in $\Ch$. We recall that this just means that $p$ is an underlying trivial fibration in $\Ch$. Suppose we have a commutative diagram 
\begin{equation}\label{diagsymm}\begin{array}{rcc}
A_n \amalg_{L_nA}L_nB & \to & X \\
\Ev_n(f \Box i) \downarrow &  & \downarrow p \\
B_n & \longrightarrow & Y.
\end{array}\end{equation}
By adjunction, this diagram has a lift if and only if the induced diagram
$$\begin{array}{ccc}
A & \longrightarrow & R_nX \\
f \downarrow &  & \downarrow  \\
B & \to & R_nY\times_{\Hom_{\SymK}(\SymKb, R_nY)} \Hom_{\SymK}(\SymKb, R_nX).
\end{array}$$
Using the fact, that $\Ev_n$ and $R_n$ is a pair of adjoint closed monoidal functors, we get an isomorphism of objects in $\SpsC$
\begin{equation}\label{closedadjunction}
\Hom_{\SymK}(A,R_nZ) \cong R_n \Hom_K(\Ev_n A,Z)
\end{equation}
where $\Hom_K$ and $\Hom_{\SymK}$ denote the closed structures of $\Ch$ and $\SpsC$, respectively. When we apply this isomorphism to the above diagram, we deduce that there is a lift in diagram(\ref{diagsymm}) if and only if there is a lift in the following diagram
\begin{equation}\label{adjointdiagsymm}
\begin{array}{ccc}
A & \longrightarrow & R_nX \\
f \downarrow &  & \downarrow q_n \\
B & \to & R_nY\times_{R_n\Hom_K(K^{\otimes n}, Y)} R_n\Hom_K(K^{\otimes n}, X).
\end{array}\end{equation}
Since $K$ is a cofibrant object in $\Ch$ and $\Ch$ being a closed monoidal model category, the functors $R_n$, $\Hom_K(K^{\otimes n},-)$ and pullbacks preserve fibrations and trivial fibrations, the map $q_n$ is a projective trivial fibration in $\SpsC$. Hence there is a lift in the last diagram, when $f$ is a projective cofibration.\\
Now, if each map $\Ev_n(f \Box i)$ is a $\Sigma_n$-cofibration in $\Ch$, one can construct by induction on $n$ a lift $B \to X$ in the commutative diagram
$$\begin{array}{ccc}
A& \stackrel{a}{\to} & X \\
f \downarrow &  & \downarrow p \\
B & \stackrel{b}{\to} & Y
\end{array}$$
where $p$ is a level trivial fibration. The argument is exactly the one given in the first part of the proof of Proposition 8.5 in \cite{hovey}.    
\end{proof}

\begin{prop}\label{symmtP-proj}
The projective cofibrations are exactly the maps which have the left lifting property with respect to all maps in $\tQ$. The projective trivial cofibrations are exactly the maps which have the left lifting property with respect to all maps in $\tP$. 
\end{prop}
\begin{proof}
By Proposition \ref{symmprojcofs} $f:A\to B$ is a projective cofibration if and only if $\Ev_n(f \Box i)$ are $\Sigma_n$-cofibrations in $\Ch$ for all $n$. Since $\Ch$ is a fibrantly generated model category with $Q$ the set of generating trivial fibrations, $\Ev_n(f \Box i)$ are $\Sigma_n$-cofibrations if and only if they have the left lifting property with respect to all maps in $Q$, i.e. if and only if $\Ev_n(f \Box i)$ are in $Q-\proj$. 
We have seen that a lift in a diagram
$$\begin{array}{rcc}
A_n \amalg_{L_nA}L_nB& \to & X \\
\Ev_n(f \Box i) \downarrow &  & \downarrow p \\
B_n & \longrightarrow & Y
\end{array}$$
is by adjunction equivalent to a lift in the diagram
$$\begin{array}{ccc}
A & \longrightarrow & R_nX \\
f \downarrow &  & \downarrow q_n \\
B & \to & R_nY\times_{R_n\Hom_K(K^{\otimes n}, Y)} R_n\Hom_K(K^{\otimes n}, X).
\end{array}$$
Hence the maps $\Ev_n(f \Box i)$ are in $Q-\proj$ if and only if the map $f$ is an element of $\tQ-\proj$. This completes the proof in the case of a cofibration. The case of a trivial cofibration is similar. 
\end{proof}

\begin{cor}\label{symmtP-fib=fib}
The projective fibrations are exactly the maps in $\tP-\fib$. The projective trivial cofibrations are exactly the maps in $\tQ-\fib$.
\end{cor}

As in the case of ordinary spectra we deduce the following 
\begin{theorem}\label{symmthmstrict}
Let $\Ch$ be a left proper fibrantly generated simplicial closed symmetric monoidal model category with all small limits. The projective weak equivalences, projective fibrations and projective cofibrations define a left proper fibrantly generated simplicial closed symmetric monoidal model structure on $\SpsC$ with set of generating fibrations $\tP$ and set of generating trivial fibrations $\tQ$ and unit object $\SymK$.
\end{theorem}
\begin{proof}
The model structure part may be deduced exactly as in Theorem \ref{thmstrict}. The only thing we have to explain is the closed symmetric monoidal structure on $\SpsC$. We have to show that given a projective cofibration $f:A \to B$ and a projectivfe fibration $g:W\to Z$ the induced map $\Hom_{\SymK , \Box}(f,g)$
$$\Hom_{\SymK}(B,W) \to \Hom_{\SymK}(B,Z)\times_{\Hom_{\SymK}(A,Z)}\Hom_{\SymK}(A,W)$$ 
is a fibration which is trivial if $f$ or $g$ is trivial. By the adjoint statement of Corollary 4.2.5 of \cite{hoveybook}, it suffices to show this if $g$ is a generating fibration or generating trivial fibration, respectively. So we suppose $g=q_n: R_nX \to R_nY\times_{R_n\Hom_K(K^{\otimes n}, Y)} R_n\Hom_K(K^{\otimes n}, X)$ in $\tP$ or $\tQ$ as above. Then the map $\Hom_{\SymK , \Box}(f,q_n)$ becomes
$$\Hom_{\SymK}(B,R_nX) \to \Hom_{\SymK}(B,R_n(X,Y))\times_{\Hom_{\SymK}(A,R_n(X,Y))}\Hom_{\SymK}(A,R_nX)$$
where we denote $R_nY\times_{R_n\Hom_K(K^{\otimes n}, Y)} R_n\Hom_K(K^{\otimes n}, X)$ by $R_n(X,Y)$. Using adjuntion (\ref{closedadjunction}) and the fact that the diagrams (\ref{diagsymm}) and (\ref{adjointdiagsymm}) are adjoint to each other, we see that the map $\Hom_{\SymK , \Box}(f,q_n)$ is in fact equal to the map $R_n \Hom_{K,\Box}(\Ev_n(f\Box i), p)$ 
$$R_n\Hom_K(B_n,X) \to R_n\Hom_K(B_n,Y)\times_{R_n\Hom_K(A_n \amalg_{L_nA}L_nB,Y)}R_n\Hom_K(A_n \amalg_{L_nA}L_nB,X)$$
where $p:X \to Y$ is the map in $\Ch$ corresponding to $q_n$, as above. 
Since $R_n$ is a right Quillen functor and since $\Ch$ is a closed symmetric monoidal model category, we get that the map $R_n \Hom_{K,\Box}(\Ev_n(f\Box i), p)$ is a fibration which is trivial if $\Ev_n(f \Box i)$ or $p$ is trivial. By Proposition \ref{symmprojcofs}, the map $f$ is a (trivial) projective cofibration if and only if the maps $\Ev_n(f \Box i)$ are (trivial) cofibrations for all $n$. Since $q_n$ is a (trivial) fibration if and only if $p$ is a (trivial) fibration, we are done.
\end{proof}

In order to obtain a stable model structure on $\SpsC$ we have to localize the projective structure at an appropriate set of fibrant objects.
\begin{defn}\label{Omega-spectra}
A symmetric spectrum $X \in \SpsC$ is a {\rm symmetric $\Omega$-spectrum} if $X$ is level fibrant and the adjoint $X_n \to X_{n+1}^K$ of the structure map of $X$ is a weak equivalence for all $n$. \\
A symmetric Spectrum $E$ is called {\rm injective} if it has the extension property with respect to every monomorphism $f$ of symmetric spectra that is also a level equivalence, i.e. for every diagram in $\SpsC$
$$\begin{array}{rcc}
X & \stackrel{g}{\longrightarrow} & E \\
f \downarrow &  & \\
Y &  & 
\end{array}$$
where $f$ is a level equivalence and a monomorphism there is a map $h:Y \to E$ such that $g =hf$.
\end{defn}

We want the symmetric $\Omega$-spectra to be the stable fibrant objects in $\SpsC$. 
As in the previous section we have to choose the set maps 
$$\Sh := \{ \zeta^A_n:F_{n+1}(A\otimes K) \to F_nA; A~\mathrm{cofibrant~in~}\Ch \}$$ to be the $S$-local equivalences, where the map $\zeta^A_n$ is adjoint to the map 
$$A\otimes K \to \Ev_{n+1}F_nA=\Sigma_{n+1}\times (A\otimes K)$$ corresponding to the identity of $\Sigma_{n+1}$.  With the set $S$ consisting of all symmetric $\Omega$-spectra we obtain the desired set $\Sh$ of  $S$-local equivalences.

\begin{defn}\label{symmstablestructure}
Let $\Ch$ be a left proper fibrantly generated simplicial model category with all small limits. Let $S$ be the set of all symmetric $\Omega$-spectra. We define the {\rm stable model structure} on $\SpsC$ to be the $S$-localized model structure $L_{\Sh}\SpsC$ of Theorem \ref{leftloc} where $\Sh$ is the class of all $S$-local equivalences.
\end{defn}
\begin{remark}
In particular, the {\rm stable equivalences} are the $S$-local equivalences; the {\rm stable cofibrations} are the projective cofibrations; the {\rm stable fibrations} are the maps that have the right lifting property with respect to all stable trivial cofibrations; the {\rm stable fibrant objects} are the symmetric $\Omega$-spectra.
\end{remark}

\begin{theorem}
The functor $-\otimes K$ is a Quillen equivalence with respect to the stable structure of $\SpsC$.
\end{theorem}
\begin{proof} The proof is given by \cite{hovey}, Theorem 8.10.\end{proof}

Now we specialize to the case of $\Ch =\hSh$, which is a simplicial fibrantly generated closed symmetric monoidal model category. We can prove the following theeorem only for this special case. In order to prove the following theorem we use the ideas of the proof of Theorem 5.3.7 of \cite{symmspectra}. in particular, we use the fact that we can build the cofiber $Ci$ of a map of spectra $i:X \to Y$.
In the following we denote the internal Hom-object in $\SpshSh$ by $\Hom_S(X,Y)$.

\begin{prop}\label{prop3.1.4}
Let $f:X \to Y$ be a map of symmetric profinite spectra. The following conditions are equivalent:\\
a) $E^0f$ is an isomorphism for every injective $\Omega$-spectrum $E$.\\
b) $\Map_{\SpshSh}(f,E)$ is a weak equivalence for every injective $\Omega$-spectrum $E$.\\
c) $\Hom_S(f,E)$ is a level equivalence for every injective $\Omega$-spectrum $E$.
\end{prop}
\begin{proof}
The proof is the one of Proposition 3.1.4 in \cite{symmspectra}. CHECK!!
\end{proof}

\begin{lemma}\label{lemma3.1.6}
Let $f:X \to Y$ be a map of symmetric profinite spectra. \\
1. If $E \in \SpshSh$ is an injective spectrum and $f$ is a level equivalence, then $E^0f$ is an isomorphism of sets.\\
2. If $f:X \to Y$ is a map of injective spectra, $f$ is a level equivalence if and only if $f$ is a simplicial homotopy equivalence.
\end{lemma}
\begin{proof}
The proof is the one of Lemma 3.1.6 in \cite{symmspectra}.  CHECK!!
\end{proof}

\begin{theorem}\label{stablemonoidalstructure}
The stable model structure on $\SpshSh$ is monoidal.
\end{theorem}
\begin{proof}
Since we have already proved that the projective model structure is monoidal and the stable cofibrations are exactly the projective cofibrations, we only have to show that for two stabel cofibrations $f$ and $g$ the map $f\Box g$ is also a stable equivalence if $f$ or $g$ is a stably trivial cofibration.\\
Suppose that $g$ is a stable equivalence. The other case is of course similar. A level cofibration $i:X \to Y$ is a stabel equivalence if and only if its cofiber $C_i=Y/X$ is stably trivial. We have already seen that $f \Box g$ is a projective and hence in particular a level cofibration. By commuting colimits, the cofiber of $f \Box g$ is the smash product $Cf \wedge Cg$ of the cofiber $Cf$ of $f$ and the cofiber $Cg$ of $g$. Let $E$ be a symmetric $\Omega$-spectrum. We show that $\Hom_S(Cf\wedge Cg,E$ is a level trivial spectrum, and thus $Cf \wedge Cg$ is stably trivial. \\
By adjointness, we have $\Hom_S(Cf \wedge Cg,E) \cong \Hom_S(Cg,\Hom_S(Cf,E))$. We will prove in the following lemma that $D:=\Hom_S(Cf,E)$ is an injective $\Omega$-spectrum. Together with the above adjuntion and the fact that $Cg$ is stably trivial by hypothesis, this proves that $Cf \wedge Cg$ is stably trivial. 

\begin{lemma}
Let $E$ be an injective spectrum and let $h$ be a projective cofibration. Then $\Hom_S(h,E)$ has the right lifting property with respect to all level monomorphisms which are also level equivalences.\\
In particular, if $A$ is a projective cofibrant spectrum, then $\Hom_S(A,E)$ is an injective spectrum.
\end{lemma}
\begin{proof}
Let $g$ and $h$ be level monomorphisms and let $g$ be a level equivalence. The $g \Box h$ is defined as a coequalizer of colimits and pushouts of the maps $g_p \Box h_q$. Since $\hShp$ is a monoidal model category, these maps are monomorphisms and weak equivalences. Since colimits and pushouts preserve trivial cofibrations when all objects are cofibrant, $g \Box h$ is still a level monomorphism and level equivalence. This implies that if $E$ is an injective spectrum, then the pair $(g \Box h,E)$ has the lifting property. By adjointness, we have an isomoprhism 
$\Hom_{\SpshSh}(g,\Hom_S(h,E)) \cong \Hom_{\SpshSh}(g \Box h,E)$. This proves that $\Hom_S(h,E)$ has the right lifting property with respect to all level trivial cofibrations.  
\end{proof}

\begin{lemma}
Let $E$ be an injective $\Omega$-spectrum and let $A$ be a projective cofibrant spectrum. Then $\Hom_S(A,E)$ is an $\Omega$-spectrum, too.
\end{lemma}
\begin{proof}
We write $D:=\Hom_S(A,E)$. By adjointness, we have an isomorphism
$$\Ev_n D\cong \Map_{\SpshSh}(A \wedge F_nS^0,E) \cong \Map_{\SpshSh}(A, \Hom_S(F_nS^0,E)).$$
Since $E$ is an $\Omega$-spectrum, 
$(F_nS^0 \wedge \lambda)^{\ast}: \Hom_S(F_nS^0,E) \to \Hom_S(F_{n+1}S^1,E)$ is a level equivalence, where $\lambda: F_1S^1 \to F_0S^0$ is the adjoint to the identity map 
$S^1 \to \Ev_1S=S^0$. Since $E$ is injective, both the source and thee target are also injective, and so this map is a simplicial homotopy equivalence by Lemma \ref{lemma3.1.6}. Hence 
$\Ev_n D \to (\Ev_{n+1}D)^{S^1}$ is still a level equivalence, so $D=\Hom_S(A,E)$ is an injective $\Omega$-spectrum.
\end{proof}
\end{proof}

\begin{cor}
The profinite completion is a monoidal left Quillen functor from symmetric spectra to profinite symmetric spectra.
\end{cor}

\end{appendix}


\bibliographystyle{amsplain}

Mathematisches Institut, WWU M\"unster, Einsteinstr. 62, 48149 M\"unster, Germany\\
E-mail address: gquick@math.uni-muenster.de

\end{document}